\documentclass[bj, numbers, noshowframe]{imsart}

\usepackage{amsfonts, amssymb, amsmath, amsthm, bbm, multirow}
\usepackage{caption, subcaption}

\theoremstyle{plain}

\newtheorem{theorem}{Theorem}
\newtheorem{lem}{Lemma}
\theoremstyle{remark}

\newtheorem{assumption}{Assumption}

\newcommand{\Expectation}{\mathbb{E}}
\newcommand{\Cov}{\mathbb{C}\mathsf{ov}}
\newcommand{\Var}{\mathbb{V}}

\newcommand{\Prob}{\mathbb{P}}
\newcommand{\Indicator}{\mathbbm{1}}

\newcommand{\precsimConc}{\precsim_{\mathtt{TC}}}

\newcommand{\hd}{\frac{1}{h^d}}
\newcommand{\n}[1]{%
\ifthenelse{\equal{#1}{}}{\frac{1}{n}}{\frac{1}{n^{#1}}}%
}

\newcommand{\K}{\mathscr{K}}

\newcommand{\Kclass}{\mathcal{K}}

\newcommand{\Trans}{\intercal}
\newcommand{\diff}{\mathrm{d}}

\renewcommand{\epsilon}{\varepsilon}

\newcommand{\bP}{\mathbf{P}}
\newcommand{\bQ}{\mathbf{Q}}
\newcommand{\bR}{\mathbf{R}}
\newcommand{\bS}{\mathbf{S}}
\newcommand{\bT}{\mathbf{T}}

\newcommand{\bX}{\mathbf{X}}
\newcommand{\bY}{\mathbf{Y}}

\newcommand{\bb}{\mathbf{b}}
\newcommand{\bc}{\mathbf{c}}

\newcommand{\be}{\mathbf{e}}

\newcommand{\bmm}{\mathbf{m}}
\newcommand{\bp}{\mathbf{p}}
\newcommand{\bq}{\mathbf{q}}

\newcommand{\bu}{\mathbf{u}}
\newcommand{\bv}{\mathbf{v}}

\newcommand{\bx}{\mathbf{x}}

\newcommand{\bz}{\mathbf{z}}

\newcommand{\pOrder}{\mathfrak{p}}
\newcommand{\qOrder}{\mathfrak{q}}

\newcommand{\bbeta}{\boldsymbol{\beta}}

\newcommand{\bnu}{\boldsymbol{\nu}}
\newcommand{\bgamma}{\boldsymbol{\gamma}}

\newcommand{\bOmega}{\boldsymbol{\Omega}}

\newcommand{\Ysupp}{\mathcal{Y}}
\newcommand{\Xsupp}{\mathcal{X}}

\newcommand{\FProc}{\mathbb{F}}
\newcommand{\BProc}{\mathbb{B}}
\newcommand{\TProc}{\mathbb{T}}
\newcommand{\SProc}{\mathbb{S}}
\newcommand{\GProc}{\mathbb{G}}

\newcommand{\inv}{{-1}}
\newcommand{\eval}{y, \bx}

\DeclareMathOperator*{\argmin}{argmin}

\DeclareMathAlphabet{\mathsf}{OT1}{LibertinusSans-LF}{m}{n}

\newcommand{\sumIN}{\sum_{i=1}^n}

\newcommand{\cval}{\mathtt{cv}}

\begin{document}

\begin{frontmatter}
\title{Boundary adaptive local polynomial conditional density estimators}
\runtitle{Local polynomial conditional density estimators}

\begin{aug}
\author[A]{\inits{MC}\fnms{Matias D.}~\snm{Cattaneo}\ead[label=e1]{cattaneo@princeton.edu}}
\author[A]{\inits{RC}\fnms{Rajita}~\snm{Chandak}\ead[label=e2]{rchandak@princeton.edu}}
\author[B]{\inits{MJ}\fnms{Michael}~\snm{Jansson}\ead[label=e3]{mjansson@berkeley.edu}}
\author[D]{\inits{XM}\fnms{Xinwei}~\snm{Ma}\ead[label=e4]{x1ma@ucsd.edu}\orcid{0000-0001-8827-9146}}
\address[A]{Department of Operations Research and Financial Engineering, Princeton University, Princeton NJ, United States\printead[presep={,\ }]{e1,e2}}

\address[B]{Department of Economics, UC Berkeley, Berkeley CA, United States\printead[presep={,\ }]{e3}}

\address[D]{Department of Economics, UC San Diego, La Jolla CA, United States\printead[presep={,\ }]{e4}}
\end{aug}

\begin{abstract}
We begin by introducing a class of conditional density estimators based on local polynomial techniques. The estimators are boundary adaptive and easy to implement. We then study the (pointwise and) uniform statistical properties of the estimators, offering characterizations of both probability concentration and distributional approximation. In particular, we establish 
uniform convergence rates in probability and valid Gaussian distributional approximations for the Studentized $t$-statistic process. We also discuss implementation issues such as consistent estimation of the covariance function for the Gaussian approximation, optimal integrated mean squared error bandwidth selection, and valid robust bias-corrected inference. We illustrate the applicability of our results by constructing valid confidence bands and hypothesis tests for both parametric specification and shape constraints, explicitly characterizing their approximation errors. A companion \texttt{R} software package implementing our main results is provided.
\end{abstract}


\end{frontmatter}

\section{Introduction}\label{section:Introduction}

Suppose that $(y_1,\bx_1^\Trans),(y_2,\bx_2^\Trans),\dots,(y_n,\bx_n^\Trans)$ is a random sample from a distribution supported on $\Ysupp\times\Xsupp$, where $\Ysupp\subset\mathbb{R}$ and $\Xsupp\subset\mathbb{R}^d$ are compact. Letting $F(y|\bx)$ be the conditional cumulative distribution function (CDF) of $y_i$ given $\bx_i$, important parameters of interest in statistics, econometrics, and many other data science disciplines, are the conditional probability density function (PDF) and derivatives thereof:
\[f^{(\vartheta)}(y|\bx) = \frac{\partial^{1+\vartheta}}{\partial y^{1+\vartheta}} F(y|\bx), \qquad \vartheta\in\{0,1,2,\dots\},\]
where, in particular, $f(y|\bx) = f^{(0)}(y|\bx)$ is the conditional density function of $y_i$ given $\bx_i$.

Estimation and inference methodology for (conditional) PDFs has a long tradition in statistics \citep[e.g.,][and references therein]{Wand-Jones_1995_Book,Wasserman_2006_Book,simonoff2012smoothing,scott2015multivariate}. Unfortunately, without specific modifications, smoothing methods employing kernel, series, or other local approximation techniques are invalid at or near boundary points of $\Ysupp\times\Xsupp$. To address this challenge, we introduce a boundary adaptive nonparametric estimator of $f^{(\vartheta)}(y|\bx)$ based on local polynomial techniques \citep{Fan-Gijbels_1996_Book} and provide an array of distributional approximation results that are valid (pointwise and) uniformly over $\Ysupp\times\Xsupp$. In particular, we obtain a uniformly valid stochastic linear representation for the estimator and develop uniform inference methods based on strong approximation techniques leading to, for example, asymptotically valid confidence bands with careful characterization of their associated approximation errors.

To motivate our proposed estimation approach, suppose we start from an estimator of the conditional CDF, $\widehat{F}(\cdot|\bx)$. Then, for $y\in\mathbb{R}$, a natural estimator of $f^{(\vartheta)}(y|\bx)$ is obtained via local polynomial regression:
\begin{align}\label{eq:local reg along y direction}
    \widehat{f}^{(\vartheta)}(y|\bx) = \be_{1+\vartheta}^\Trans\widehat{\bbeta}(y|\bx), \qquad
      \widehat{\bbeta}(y|\bx) = \argmin_{\bu\in\mathbb{R}^{\pOrder+1}} \sum_{i=1}^n \left(\widehat{F}(y_i|\bx) - \bp(y_i-y)^\Trans\bu \right)^2 K_h(y_i;y),
\end{align}
where $\pOrder\geq1+\vartheta$ is the order of the polynomial basis
$\bp(y)=(1,y/1!,y^2/2!,\dots,y^\pOrder/\pOrder!)^\Trans$, $\be_l$ is the
conformable $(1+l)$-th unit vector, and $K_h(y_i;y)=K((y_i-y)/h)/h$ for some
kernel function $K$ and some positive bandwidth $h$. Since
$F(y|\bx_i)=\Expectation[\Indicator(y_i\leq y)|\bx_i]$, we employ a $\qOrder$-th
order local polynomial regression of the indicator function, $\Indicator(y_i\leq
y)$, to form the conditional CDF estimator that will be plugged into \eqref{eq:local reg along y direction}:
\[\widehat{F}(y|\bx) = \be_{\mathbf{0}}^\Trans\widehat{\bgamma}(y|\bx), \qquad
  \widehat{\bgamma}(y|\bx) = \argmin_{\bv\in\mathbb{R}^{\qOrder_d+1}} \sum_{i=1}^n \left(\Indicator(y_i\leq y) - \bq(\bx_i-\bx)^\Trans\bv \right)^2 L_b(\bx_i;\bx).
\]
Here, using standard multi-index notation, $\bq(\bx)$ denotes the $(\qOrder_d+1)$-dimensional vector collecting the terms $\bx^{\bmm}/\bmm ! = x_1^{m_1} \cdots x_d^{m_d} / (m_1! \cdots m_d!)$ for $\bx = (x_1, \dots, x_d)^\Trans \in \mathbb{R}^d$, $\bmm = (m_1, \dots, m_d)^\Trans \in \mathbb{Z}_+^d$ with $|\bmm| = m_1 + \dots + m_d \leq \qOrder$, and $\qOrder_d = (d+\qOrder)!/(\qOrder!d!)-1$. We also let $L_b(\bx_i;\bx)=L((\bx_i-\bx)/b)/b^d$ be some (multivariate) kernel function $L$ and positive bandwidth $b$. Our proposed estimator can also be written in closed-form as
\begin{align}\label{equation:closed-form expression}
    \widehat{f}^{(\vartheta)}(y|\bx) = \be_{1+\vartheta}^\Trans\widehat{\bS}_y^{-1} \widehat{\bR}_{y,\bx} \widehat{\bS}_\bx^{-1}\be_\mathbf{0},
\end{align}
where the matrices are
\begin{align*}
    \widehat{\bS}_y &= \frac{1}{n}\sum_{i=1}^{n}\bp\Big(\frac{y_i-y}{h}\Big)\frac{1}{h}\bP\Big(\frac{y_i-y}{h}\Big)^\Trans,
    \qquad
    \widehat{\bS}_\bx     = \frac{1}{n}\sum_{i=1}^{n}\bq\Big(\frac{\bx_i-\bx}{b}\Big)\frac{1}{b^d}\bQ\Big(\frac{\bx_i-\bx}{b}\Big)^\Trans,\\
    \widehat{\bR}_{y,\bx} &= \frac{1}{n^2h^{1+\vartheta}}\sum_{j=1}^{n}\sum_{i=1}^{n}\frac{1}{h}\bP\Big(\frac{y_j-y}{h}\Big)\frac{1}{b^d}\bQ\Big(\frac{\bx_i-\bx}{b}\Big)^\Trans  \Indicator(y_i\leq y_j),
\end{align*}
with the definitions $\bP(y)=\bp(y)K(y)$ and $\bQ(\bx)=\bq(\bx)L(\bx)$, which absorb the kernel function into the basis. See Appendix \ref{Appendix: closed-form expression and comparison with FYT} for derivation.

By virtue of being based on a local polynomial smoothing approach, the estimator $\widehat{f}^{(\vartheta)}(y|\bx)$ is not only intuitive, but also boundary adaptive. Furthermore, $\widehat{f}^{(\vartheta)}(y|\bx)$ admits a simple closed-form representation as we have shown in \eqref{equation:closed-form expression}, making it easy to implement. These features follow directly from its construction: unlike classical kernel-based conditional density (derivative) estimators, which seek to approximate the conditional PDF indirectly (e.g., by constructing a ratio of two unconditional kernel-based density estimators), our proposed estimator applies local polynomial techniques directly to the conditional CDF estimator $\widehat{F}(y|\bx)$. In addition, our approach offers an easy way to construct higher-order kernels to reduce misspecification (or smoothing) bias via the choice of polynomial orders $\pOrder$ and $\qOrder$.

We present two main uniform results for our proposed estimator. First, we provide precise uniform probability concentration bounds associated with a stochastic linear representation of $\widehat{f}^{(\vartheta)}(y|\bx)$ (Lemma \ref{lem:Equivalent Kernel Representation} and Theorem \ref{thm:Probability Concentration}). In addition to being useful for the purposes of characterizing the distributional properties of the conditional density estimator itself, the first main result can be used to analyze multi-step estimation and inference procedures whenever $\widehat{f}^{(\vartheta)}(y|\bx)$ enters as a preliminary step. As a by-product of the development of the first main result, we obtain a related class of conditional density estimators based on local smoothing. This new approach will require the knowledge of the support $\mathcal{Y}$. On the other hand, it is immune to ``low'' density regions of $y_i$. For details, see Appendix \ref{appendix:local smoothing based estimator}.

Our second main result employs the stochastic linear representation of $\widehat{f}^{(\vartheta)}(y|\bx)$ to establish a valid strong approximation for the standardized $t$-statistic stochastic process based on $\widehat{f}^{(\vartheta)}(y|\bx)$ and indexed over $\Ysupp\times\Xsupp$ (Theorem \ref{thm:Strong Approximation S-hat}). This result is established using a powerful result due to \citet{Rio_1994_PTRF}, which in turn builds on the celebrated Hungarian construction \citep{KMT_1975_Zeitschrift}. The $t$-statistic stochastic processes based on kernel-based nonparametric estimators are not asymptotically tight and, as a consequence, do not converge weakly as a process indexed over $\Ysupp\times\Xsupp$ \citep{van-der-Vaart-Wellner_1996_Book,Gine-Nickl_2016_Book}. Nevertheless, using strong approximations to such processes, it is possible to deduce distributional approximations for functionals thereof by employing anti-concentration \citep{Chernozhukov-Chetverikov-Kato_2014a_AoS}. Combining these ideas, we obtain valid distributional approximations for the suprema of the $t$-statistic stochastic process (Theorem \ref{thm:Kolmogorov-Smirnov Distance: Suprema}) based on $\widehat{f}^{(\vartheta)}(y|\bx)$ with approximation rates that are faster than those currently available in the literature for the case of $d=1$ (e.g., Remark 3.1(ii) in \cite{Chernozhukov-Chetverikov-Kato_2014b_AoS}).

In addition to our two main uniform estimation and distributional results, we discuss several implementation results that are useful for practice. First, we present a covariance function estimator for the Gaussian approximation and prove its uniform consistency (Lemma \ref{lem:Covariance Estimation}). This result enables us to estimate the statistical uncertainty underlying the Gaussian approximation for a feasible version of the $t$-statistic process. Second, in Section \ref{section:Applications} we discuss optimal bandwidth selection based on an asymptotic approximation to the integrated mean squared error (IMSE) of the estimator $\widehat{f}^{(\vartheta)}(y|\bx)$. This result allows us to implement our proposed estimator using point estimation optimal data-driven bandwidth selection rules. Finally, we employ robust bias correction \citep{Calonico-Cattaneo-Farrell_2018_JASA,Calonico-Cattaneo-Farrell_2022_Bernoulli} to develop valid inference methods based on the Gaussian approximation when using the estimated covariance function and IMSE-optimal bandwidth rule.

We illustrate our theoretical and methodological results with three substantive applications in Section \ref{section:Applications}. To be specific, we construct valid confidence bands for the unknown conditional density function (and derivatives thereof) and we develop valid hypothesis testing procedures for parametric specification and shape constraints of $f^{(\vartheta)}(y|\bx)$, respectively. All these methods are data-driven and, in some cases, optimal in terms of probability and/or distributional concentration, possibly up to $\log(n)$ factors. Furthermore, thanks to the precise probability approximation errors we obtain via strong approximation and other exponential concentration methods, we are able to characterize precise coverage error and rejection probability error rates for all the feasible inference procedures considered.

Another advantage of our proposed estimation procedure \eqref{eq:local reg along y direction} is that it allows for incorporating additional constraints easily. For example, setting $\vartheta=0$ (PDF), it may be desirable to require that the estimator is non-negative and integrates to 1. In Section \ref{section:additional constraints}, we proposed a modified conditional PDF estimator which satisfies these two properties. To be precise, non-negativity can be imposed by solving a constrained version of \eqref{eq:local reg along y direction}, as the feature is local to the evaluation point. On the other hand, ensuring the estimator integrates to 1 requires imposing a global constraint, which we implement by minimizing the Kullback-Leibler divergence to ensure that the final estimator is a valid conditional density in finite samples. Interestingly, this modified conditional PDF estimator requires introducing a normalization factor that affects the strong approximation in nontrivial ways, leading to a different distributional Gaussian process approximation (Theorem \ref{thm:Strong Approximation S-hat, NN and I estimators}).

Proofs of the main results are given in the \hyperref[appn]{Appendix}. In the supplementary material \citep{Cattaneo-Chandak-Jansson-Ma_2023_OnlineSupp}, we consider a more general setup and offer additional technical and methodological results of potential independent interest, including: (i) boundary adaptive estimators for the CDF and its derivatives with respect to the conditioning variable $\bx$; (ii) theoretical properties of the local smoothing based conditional PDF and derivatives estimators; (iii) additional details on bandwidth selection; (iv) alternative covariance function estimators. Last but not least, we provide a general purpose \texttt{R} software package (\texttt{lpcde}) implementing the main results in this paper.

\subsection{Related literature}\label{subsection:Related Literature}

Our paper contributes to the literature on kernel-based conditional density estimation and inference. See \citet{Hall-Wolff-Yao_1999_JASA}, \citet{DeGooijer-Zerom_2003_SN} and \citet{Hall-Racine-Li_2004_JASA} for earlier reviews, and \citet{Wand-Jones_1995_Book}, \citet{Wasserman_2006_Book}, \citet{simonoff2012smoothing} and \citet{scott2015multivariate} for textbook introductions. Traditional methods for conditional density estimation typically employ ratios of unconditional kernel density estimators, nonlinear kernel-based derivative of distribution function estimators, or local polynomial estimators based on some preliminary density-like approximation. In the leading special case of $\vartheta=0$, the closest antecedent to our proposed conditional density estimator is the local polynomial conditional density estimator introduced by \citet{Fan-Yao-Tong_1996_Biometrika}, which is formed by a local polynomial regression of $K_h(y_i;y)$ on $\bx_i$. Their estimator is valid at the boundary of $\mathcal{X}$, but is generally inconsistent at the boundary of $\mathcal{Y}$. See Appendix \ref{Appendix: closed-form expression and comparison with FYT} for more discussion.

More generally, classical methods for conditional density estimation are not boundary adaptive without specific modifications, and in some cases do not have a closed-form representation. Boundary adaptivity could be achieved by employing boundary-corrected kernels in some cases, but such conditional density estimation methods do not appear to have been considered in the literature before. Our first contribution is to introduce a novel boundary adaptive, closed-form conditional density (derivative) estimator. Our proposed construction does not rely on boundary-corrected kernels explicitly, but it rather builds on the idea that automatic boundary-adaptive density estimators can be constructed using local polynomial methods to smooth out the (discontinuous) distribution function \citep{Cattaneo-Jansson-Ma_2020_JASA}.

We also consider estimation of conditional CDF, as the intercept in Equation \eqref{eq:local reg along y direction} is an estimator of $F(y|\bx)$, that is, $\be_{0}^\Trans\widehat{\bbeta}(y|\bx)$. In addition to being boundary adaptive, this CDF estimator is also continuous in $y$ and $\bx$. We discuss properties of this estimator (probability concentration, strong approximation, etc.) in the supplementary material. To compare, the conditional CDF estimator $\widehat{F}(y|\bx)$, which is constructed via a local polynomial regression of the indicators $\Indicator(y_i\leq y)$ on $\bx_i$, is generally discontinuous in $y$. Properties of $\widehat{F}(y|\bx)$, such as the uniform convergence rate, have been studied in the literature \citep{ferrigno2010uniform,einmahl2000empirical}.

\subsection{Notation and assumptions}\label{subsection:Notation and Assumptions}

To simplify the presentation, in the remainder of this paper we set $L$ to be the product kernel based on $K$: $L(\bx)=K(x_1)K(x_2)\cdots K(x_d)$ for a vector $\bx = (x_1,\dots,x_d)^\Trans$. We also employ the same bandwidth, $b=h$, in the construction of our proposed estimator, and assume $\qOrder=\pOrder - \vartheta - 1\geq 0$ throughout.

For two numbers $a$ and $b$, let $a\vee b = \max\{a,b\}$. Limits are taken with respect to the sample size tending to infinity  (i.e., $n\to\infty$). For two positive sequences $a_n$ and $b_n$, $a_n\precsim b_n$ means that $a_n/b_n$ is bounded and $a_n\precsim_\Prob b_n$ means that $a_n/b_n$ is bounded in probability. Constants that do not depend on the sample size or the bandwidth will be denoted by $\mathfrak{c}$, $\mathfrak{c}_1$, $\mathfrak{c}_2$, etc.

We introduce the notation $\precsimConc$, which not only provides an asymptotic order in probability, but also controls the tail probability ($\mathtt{TC}$): $a_n \precsimConc b_n$ implies that for any $\mathfrak{c}_1 > 0$, there exists some $\mathfrak{c}_2$ such that
\begin{align*}
    \limsup_{n\to\infty}\ n^{\mathfrak{c}_1}\ \Prob\big[a_n \geq \mathfrak{c}_2 b_n\big] < \infty.
\end{align*}
Finally, let $\bX = (\bx_1^\Trans,\dots,\bx_n^\Trans)^\Trans_{\phantom{n}}$ and $\bY = (y_1,\dots,y_n)^\Trans_{\phantom{n}}$ be the data matrices. We make the following assumptions on the joint distribution and the kernel function.

\begin{assumption}[DGP]\label{as:1-DGP}
(i) $(y_1,\bx_1^\Trans),\dots,(y_n,\bx_n^\Trans)$ is a random sample from an absolutely continuous distribution supported on $\mathcal{Y}\times\mathcal{X} = [0,1]^{1+d}$, and the joint Lebesgue density, $f(y,\bx)$, is continuous and bounded away from zero on $\mathcal{Y}\times\mathcal{X}$. (ii) $f^{(\pOrder)}(y|\bx)$ exists and is continuous. (iii) $\partial^{\bnu} f^{(\vartheta)}(y|\bx)/\partial \bx^{\bnu} $ exists and is continuous for all $|\bnu|=\pOrder-\vartheta$.
\end{assumption}

\begin{assumption}[Kernel]\label{as:2-Kernel}
$K$ is a symmetric, Lipschitz continuous PDF supported on $[-1,1]$.
\end{assumption}

Setting $\mathcal{Y}\times\mathcal{X}= [0,1]^{1+d}$ is a normalization without loss of generality: all our results generalize to the case that $\mathcal{Y}\times\mathcal{X}$ is a Cartesian product of closed intervals. Since our method is local in nature, all the pointwise properties (discussed in the supplementary material) continue to hold if the support $\mathcal{Y}\times\mathcal{X}$ is unbounded. Statements of uniform properties will also remain valid for compact subsets.

We also follow the literature to classify evaluation points as interior or (near) boundary (for example, Section 2.1.2 of \cite{cheng1994boundary}). To be precise, let $\mathrm{Cube}_h(y,\bx) = [y-h,y+h]\times [x_1-h,x_1+h]\times \dots \times [x_d-h,x_d+h]$ be the cube of length $2h$ centered at $(y,\bx)$. Then $(y,\bx)$ is interior if $\mathrm{Cube}_h(y,\bx)\subseteq \mathcal{Y}\times\mathcal{X}$. Otherwise it is called (near) boundary. This classification stems from properties of our estimator: as discussed in Appendix \ref{Appendix:properties of K-circle}, the equivalent kernel is compactly supported, meaning that the estimator only employs observations in an $h$-neighborhood of the evaluation point.

\section{Main results}\label{section:Main Results}

This section presents four main theoretical results. First, we provide a stochastic linearization of our estimator (Lemma \ref{lem:Equivalent Kernel Representation}). Based on this representation, we obtain a uniform probability concentration result for $\widehat{f}^{(\vartheta)}(y|\bx)$ (Theorem \ref{thm:Probability Concentration}). Next, we obtain valid strong approximation results for the standardized $t$-process based on $\widehat{f}^{(\vartheta)}(y|\bx)$ (Theorem \ref{thm:Strong Approximation S-hat}). Finally, we develop a feasible distributional approximation for the suprema of the Studentized $t$-process (Theorem \ref{thm:Kolmogorov-Smirnov Distance: Suprema}). We obtain a uniform consistency result for an estimator of the covariance function (Lemma \ref{lem:Covariance Estimation}) to establish Theorem~\ref{thm:Kolmogorov-Smirnov Distance: Suprema}.

\subsection{Stochastic linearization and uniform probability concentration}\label{subsection:Uniform Probability Concentration}

We first define the large-sample limits of the matrices $\widehat{\bS}_y$ and $\widehat{\bS}_\bx$:
\begin{align*}
  \bS_y       &= \int_\Ysupp \bp\Big(\frac{u-y}{h}\Big)\frac{1}{h}\bP\Big(\frac{u-y}{h}\Big)^\Trans \text{d}F_y(u)
                \qquad\text{and}\qquad
                \bS_\bx     = \int_\Xsupp \bq\Big(\frac{\bv-\bx}{h}\Big)\frac{1}{h^d}\bQ\Big(\frac{\bv-\bx}{h}\Big)^\Trans  \mathrm{d}F_\bx(\bv),
\end{align*}
with $F_y$ and $F_\bx$ denoting the CDFs of $y_i$ and $\bx_i$, respectively. The following uniform stochastic linear representation holds for $\widehat{f}^{(\vartheta)}(y|\bx)$.

\begin{lem}[Stochastic linearization]\label{lem:Equivalent Kernel Representation}
    Suppose Assumptions \ref{as:1-DGP} and \ref{as:2-Kernel} hold. If $nh^{1+d}/\log (n)\to\infty$ and $h\to 0$, then
    \[
    \sup_{y\in\Ysupp,\bx\in\Xsupp} \left| \widehat{f}^{(\vartheta)}(y|\bx) - f^{(\vartheta)}(y|\bx) - \bar{f}^{(\vartheta)}(y|\bx) \right|\ \precsimConc\ \mathtt{r_{SL}},\quad \mathtt{r_{SL}}=h^{\pOrder-\vartheta} + \frac{\log (n)}{\sqrt{n^2h^{1+2\vartheta+d + (2 \vee d)}}},
    \]
    where $\bar{f}^{(\vartheta)}(y|\bx) = n^{-1}\sum_{i=1}^n \K^\circ_{\vartheta,h}\big( y_i,\bx_i; y,\bx \big)$, and
    \[\K_{\vartheta,h}^\circ( a,\bb; y,\bx )
          = \frac{1}{h^{1+\vartheta}}\be_{1+\vartheta}^\Trans\bS_y^\inv \int_{\mathcal{Y}} \Big(\Indicator(a\leq u) - F(u|\bb ) \Big) \frac{1}{h}\bP\Big(\frac{u-y}{h}\Big)\diff F_y(u) \frac{1}{h^d}\bQ\left(\frac{\bb-\bx}{h}\right)^\Trans \bS_{\bx}^\inv \be_{\mathbf{0}}.
        \]
\end{lem}

The proof, given in Appendix \ref{Appendix:proof of stochastic linear rep},
involves showing that the matrices $\widehat{\bS}_y$, $\widehat{\bS}_\bx$ and
$\widehat{\bR}_{y, \bx}$ concentrate.
$\widehat{\bS}_y$ and $\widehat{\bS}_\bx$
concentrate in probability (and $\mathtt{TC}$ sense),
uniformly in $y$ and $\bx$ respectively, around $\bS_y$ and $\bS_{\bx}$. Characterizing the large-sample behavior of the matrix $\widehat{\bR}_{y,\bx}$
in \eqref{equation:closed-form expression} requires a little more care, but the
end result can be combined with the results for $\widehat{\bS}_y$ and
$\widehat{\bS}_\bx$ to obtain the uniform stochastic linear
representation for $\widehat{f}^{(\vartheta)}(y|\bx)$.

Lemma~\ref{lem:Equivalent Kernel Representation} implies that the properties of $\widehat{f}^{(\vartheta)}(y|\bx)$ are thus governed by the properties of the stochastic linear representation. In Appendix \ref{Appendix:properties of K-circle}, we first characterize the leading variance of $\bar{f}^{(\vartheta)}(y|\bx)$ (Lemma \ref{lem:variance of f-bar}). Define $\mathsf{V}_{\vartheta}(y,\bx) := \Var[ \bar{f}^{(\vartheta)}(y|\bx)]$, then
\begin{align}
    \nonumber    \mathsf{V}_{\vartheta}(y,\bx)  &= \frac{1}{nh^{1+d+2\vartheta}}f(y|\bx) \left(\be_{1+\vartheta}^\Trans\bS_{y}^{-1}\bT_y \bS_{y}^{-1}\be_{1+\vartheta}\right)\left(\be_{\mathbf{0}}^\Trans \bS_{\bx}^{-1} \bT_{\bx} \bS_{\bx}^{-1} \be_{\mathbf{0}}\right) + O\left(\frac{1}{nh^{d+2\vartheta}}\right),\\
    \nonumber    \text{where}\quad \bT_{y} &= \iint_{\Ysupp\times \Ysupp}   \frac{\min(u_1,u_2)-y}{h}\frac{1}{h^2}\bP\Big(\frac{u_1-y}{h}\Big)\bP\Big(\frac{u_2-y}{h}\Big)^\Trans \text{d}F_y(u_1)\text{d}F_y(u_2),\qquad\qquad\ \\
    \label{eq:leading variance}    \bT_{\bx} &= \int_{\Xsupp}   \frac{1}{h^d}\bQ\Big(\frac{\bv-\bx}{h}\Big)\bQ\Big(\frac{\bv-\bx}{h}\Big)^\Trans \mathrm{d}F_\bx(\bv).
\end{align}

Based on the stochastic linearization result in Lemma \ref{lem:Equivalent Kernel Representation} and the above leading variance characterization, we can obtain a pointwise (in $y$ and $\bx$) convergence rate of our estimator: $h^{\pOrder-\vartheta} + 1/\sqrt{nh^{1+d+2\vartheta}}$. In Theorem \ref{thm:Probability Concentration} below we will establish a uniform convergence rate and a probability concentration result.

Appendix \ref{Appendix:properties of K-circle} establishes additional important features of $\K^\circ_{\vartheta,h}$, such as boundedness and Lipschitz continuity which will play a crucial role in our strong approximation results. We also bound the uniform covering number for the class of functions formed by varying the evaluation point. This uniform covering number result takes into account the fact that the shape of $\K^\circ_{\vartheta,h}$ changes across different evaluation points. To this end, we provide in Appendix \ref{Appendix: covering number} a generic result on covering number calculation for function classes formed by kernels, which may be of independent interest. This result allows the kernel functions to take different shapes as well as to depend on a range of bandwidths --- the latter feature can be useful for establishing consistency and distributional approximation that are uniform in bandwidth (for example, \cite{einmahl2005uniform}). However, we do not further pursue along this uniform-in-bandwidth direction to avoid obscuring the main message of the paper.

The following theorem gives a uniform probability concentration result for our conditional density and derivative estimator. The proof is in Appendix \ref{Appendix:Proof uniform convergence rate}.

\begin{theorem}[Probability concentration]\label{thm:Probability Concentration}
    Suppose Assumptions \ref{as:1-DGP} and \ref{as:2-Kernel} hold. If $h\to 0$ and if $nh^{1+d}/\log (n)\to\infty$, then
    \begin{align*}
        \sup_{y\in \Ysupp, \bx \in \Xsupp}\left| \widehat{f}^{(\vartheta)}(y|\bx) - f^{(\vartheta)}(y|\bx) \right|\ \precsimConc\ \mathtt{r_{PC}},\quad \mathtt{r_{PC}}=h^{\pOrder-\vartheta} + \sqrt{\frac{\log(n)}{nh^{1+d+2\vartheta}}}.
    \end{align*}
\end{theorem}

The $h^{\pOrder-\vartheta}$ in Theorem~\ref{thm:Probability Concentration} stems from a bias term whose magnitude coincides with that of the pointwise bias at interior evaluation points. As a consequence, the theorem implies that the estimator is boundary adaptive.
The other term represents ``noise,'' whose magnitude is larger than its counterpart in Lemma \ref{lem:Equivalent Kernel Representation}, reflecting the fact that the estimation error  $\widehat{f}^{(\vartheta)}(y|\bx) - f^{(\vartheta)}(y|\bx)$ can be characterized by the bias and the randomness in $\bar{f}^{(\vartheta)}(y|\bx)$. By setting $h=(\log(n)/n)^{\frac{1}{1+d+2\pOrder}}$, it follows from the theorem that the estimator achieves the minimax optimal uniform convergence rate \cite{khas1979lower}, namely $\left({\log(n)}/{n}\right)^{\frac{\pOrder-\vartheta}{1+d+2\pOrder}}$.

\subsection{Strong approximation}\label{subsection:Strong Approximation}

We study the distributional properties of the standardized process $\widehat{\SProc}_{\vartheta}(y,\bx)$:
\begin{align}\label{eq: process-standardized}
    \widehat{\SProc}_{\vartheta}(y,\bx) = \frac{\widehat{f}^{(\vartheta)}(y|\bx) - f^{(\vartheta)}(y|\bx)}{\sqrt{\mathsf{V}_{\vartheta}(y,\bx)}}.
\end{align}

Using elementary tools, Theorem 2.1 in the supplementary material obtains a pointwise Gaussian approximation to $\widehat{\SProc}_{\vartheta}(y,\bx)$. However, the process $\widehat{\SProc}_{\vartheta}$ is not asymptotically tight and hence it does not converge weakly to a Gaussian process in $\ell^\infty(\Ysupp\times\Xsupp)$, the set of uniformly bounded real-valued functions on $\Ysupp\times\Xsupp$ equipped with the uniform norm \citep{van-der-Vaart-Wellner_1996_Book,Gine-Nickl_2016_Book}. To obtain a uniform distributional approximation, we use the result of \citet{Rio_1994_PTRF} and establish a strong approximation result for $(\widehat{\SProc}_{\vartheta}(y,\bx) : y\in\mathcal{Y},\bx\in\mathcal{X})$. To state the result, define the correlation function
\[
\rho_{\vartheta}(y,\bx,y',\bx') = {\mathsf{C}_{\vartheta}(y,\bx,y',\bx')} \Big/ { \sqrt{\mathsf{V}_{\vartheta}(y,\bx)\mathsf{V}_{\vartheta}(y',\bx')}},
\]
where $\mathsf{C}_{\vartheta}(y,\bx,y',\bx') = n^{-1} \mathbb{E}[\K^\circ_{\vartheta,h}( y_i,\bx_i; y,\bx)\K^\circ_{\vartheta,h}( y_i,\bx_i; y',\bx')]$.

\begin{theorem}[Strong approximation]\label{thm:Strong Approximation S-hat}
    Suppose Assumptions \ref{as:1-DGP} and \ref{as:2-Kernel} hold. If ${nh^{1+d+2\pOrder}}\to 0$ and if $nh^{1+d}/\log (n)\to\infty$, then there exist two stochastic processes, $\widehat{\SProc}^{\prime}_{\vartheta}$ and $\GProc_{\vartheta}$, in a possibly enlarged probability space, such that:
    \begin{enumerate}[label=(\roman*),noitemsep]
        \item $\widehat{\SProc}_{\vartheta}$ and $\widehat{\SProc}^{\prime}_{\vartheta}$ have the same distribution,
        \item $\GProc_{\vartheta}$ is a centered Gaussian process with unit variance and correlation $\rho_{\vartheta}$;
        \item the following holds:
        \begin{align*}
              \sup_{y \in \Ysupp, \bx \in \Xsupp}\Big|\widehat{\SProc}^{\prime}_{\vartheta}(y,\bx)-\GProc_{\vartheta}(y,\bx)\Big|\ \precsimConc\  \mathtt{r_{SA}},\quad \mathtt{r_{SA}} = \sqrt{nh^{1+d+2\pOrder}} +  \Big(\frac{\log^{1+d}(n)}{nh^{1+d}}\Big)^{\frac{1}{2+2d}}.
              \end{align*}
    \end{enumerate}

\end{theorem}

The theorem provides a Gaussian approximation for the entire stochastic process $\widehat{\SProc}_{\vartheta}$ rather than for a particular functional thereof. Later we will employ this result to approximate the distribution of the suprema of the process, based on which uniform confidence bands can be constructed.

\subsection{Variance-covariance estimation and suprema approximation}\label{subsection:Suprema Approximation}

Because both the process $\widehat{\SProc}_{\vartheta}$ and the correlation function $\rho_{\vartheta}$ depend on unknown features of the underlying data generating process (namely, the covariance function $\mathsf{C}_{\vartheta}$), Theorem \ref{thm:Strong Approximation S-hat}  in isolation cannot be used for inference. In this subsection we first propose an estimator of the covariance function, and then demonstrate how to obtain a feasible distributional approximation for the suprema of the Studentized $t$-process.

The covariance function $\mathsf{C}_{\vartheta}$ can be expressed as a functional of two unknowns: the conditional CDF of $y_i$ given $\bx_i$ and the marginal CDF of $y_i$. Replacing $F(y|\bx)$ and $F_y(y)$ with $\widehat{F}(y|\bx)$ and $\widehat{F}_y(y)=n^{-1}\sum_{i=1}^n \Indicator(y_i\leq y)$, respectively, we obtain the following plug-in covariance function estimator:
\begin{align*}
    \widehat{\mathsf{C}}_{\vartheta}(y,\bx,y',\bx') &= \frac{1}{n^2}\sumIN \widehat{\K}_{\vartheta,h}^\circ\Big( y_i,\bx_i; y,\bx \Big)\widehat{\K}_{\vartheta,h}^\circ\Big( y_i,\bx_i; y',\bx' \Big),
\end{align*}
where
\begin{align*}
    \widehat{\K}_{\vartheta,h}^\circ\Big( a,\bb; y,\bx \Big)
    &= \frac{1}{h^{1+\vartheta}}\be_{1+\vartheta}^\Trans\widehat{\bS}_y^\inv \bigg[\frac{1}{n}\sum_{j=1}^n \Big(\Indicator(a\leq y_j) - \widehat{F}(y_j|\bb ) \Big) \frac{1}{h}\bP\Big(\frac{y_j-y}{h}\Big)\bigg] \frac{1}{h^d}\bQ\left(\frac{\bb-\bx}{h}\right)^\Trans \widehat{\bS}_{\bx}^\inv \be_{\mathbf{0}}.
\end{align*}
The corresponding estimators of $\mathsf{V}_{\vartheta}$ and $\rho_{\vartheta}$ are given by $\widehat{\mathsf{V}}_{\vartheta}(y,\bx) = \widehat{\mathsf{C}}_{\vartheta}(y,\bx,y,\bx)$ and
\begin{align*}
\widehat{\rho}_{\vartheta}(y,\bx,y',\bx') &= \widehat{\mathsf{C}}_{\vartheta}(y,\bx,y',\bx') \Big/ \sqrt{\widehat{\mathsf{V}}_{\vartheta}(y,\bx)\widehat{\mathsf{V}}_{\vartheta}(y',\bx')}.
\end{align*}
Lemma~\ref{lem:Covariance Estimation} establishes a uniform probability
concentration result for $\widehat{\mathsf{V}}_{\vartheta}$ and $\widehat{\rho}_{\vartheta}$. We relegate the proof to the supplementary material as it is quite involved.

\begin{lem}[Covariance estimation]\label{lem:Covariance Estimation}
    Suppose Assumptions \ref{as:1-DGP} and \ref{as:2-Kernel} hold. If $h\to 0$ and if\\ $nh^{1+d}/\log (n)\to\infty$, then
    \begin{align*}
    &\sup_{y \in \mathcal{Y}, \bx \in \mathcal{X}}\left| \frac{\widehat{\mathsf{V}}_{\vartheta}(\eval)  - {\mathsf{V}}_{\vartheta}(\eval) }{{\mathsf{V}}_{\vartheta}(\eval) } \right|\ \precsimConc\ \mathtt{r_{VE}},\  \sup_{y,y'\in\Ysupp,\bx,\bx'\in\Xsupp} \Big| \widehat{\rho}_{\vartheta}(y,\bx,y',\bx') - {\rho}_{\vartheta}(y,\bx,y',\bx') \Big|\ \precsimConc\ \mathtt{r_{VE}},\\
    &\qquad\qquad\qquad\text{where }\mathtt{r_{VE}} = h^{\pOrder-\vartheta-\frac{1}{2}} + \sqrt{\frac{\log(n)}{nh^{1+d}}}.
    \end{align*}

\end{lem}

With a valid covariance (and variance) estimator, we replacing $\mathsf{V}_{\vartheta}(y,\bx)$ with $\widehat{\mathsf{V}}_{\vartheta}(y,\bx)$ in \eqref{eq: process-standardized} to obtain the Studentized $t$-process,
\begin{align*}
    \widehat{\TProc}_{\vartheta}(y,\bx) &= \frac{\widehat{f}^{(\vartheta)}(y|\bx) - f^{(\vartheta)}(y|\bx)}{\sqrt{\widehat{\mathsf{V}}_{\vartheta}(y,\bx)}}.
\end{align*}
By Theorem \ref{thm:Strong Approximation S-hat} and Lemma \ref{lem:Covariance Estimation}, the law of $(\widehat{\TProc}_{\vartheta}(y,\bx) : y\in\mathcal{Y},\bx\in\mathcal{X})$ can be approximated by that of a centered Gaussian process with unit variance and correlation function $\rho_{\vartheta}$, where the latter is estimated by $\widehat{\rho}_{\vartheta}$. As a consequence, functionals of $\widehat{\TProc}_{\vartheta}$ admit feasible distributional approximations. To illustrate this general phenomenon, the following theorem gives a result for the supremum of $\big| \widehat{\TProc}_{\vartheta} \big|$. We define $\widehat{\GProc}_{\vartheta}$ as a process whose law, conditional on the data, is a centered Gaussian with unit variance and correlation function $\widehat{\rho}_{\vartheta}$.

\begin{theorem}[Kolmogorov-Smirnov distance: suprema]\label{thm:Kolmogorov-Smirnov Distance: Suprema}
    Suppose Assumptions \ref{as:1-DGP} and \ref{as:2-Kernel} hold. If $n\log(n)h^{1+d+2\pOrder}\to 0$ and if $nh^{1+d}/\log (n)\to\infty$, then
    \begin{align*}
    &\sup_{u \in \mathbb{R}}\left| \Prob\Big[ \sup_{y \in \Ysupp, \bx \in \Xsupp} \big| \widehat{\TProc}_{\vartheta}(y,\bx) \big| \leq u \Big] - \Prob\Big[ \sup_{y \in \Ysupp, \bx \in \Xsupp}\big|\widehat{\GProc}_{\vartheta}(y,\bx)\big|\leq u \Big| \bX,\bY  \Big] \right| \precsim_\Prob \mathtt{r_{KS}}\\
    &\qquad\qquad\text{where }\mathtt{r_{KS}} = \sqrt{n\log(n)h^{1+d+2\pOrder}} + \left(\frac{\log^{2+2d}(n)}{nh^{1+d}}\right)^{\frac{1}{2+2d}} + \left(\frac{\log^{5}(n)}{nh^{1+d}}\right)^{\frac{1}{4}}.
    \end{align*}
\end{theorem}

To compare the rate of distributional approximation with existing results, we follow the literature and ignore the first (smoothing bias) term. Then, the rate matches what \citet{Chernozhukov-Chetverikov-Kato_2014b_AoS} obtained when $d=2$ (see their Remark 3.1(ii)), but it is strictly faster when $d=1$.

\section{Applications}\label{section:Applications}

This section illustrates our theoretical and methodological results by means of three applications. Before turning to these applications, we discuss bandwidth selection, a necessary step for implementation. It is customary to select the bandwidth by minimizing an approximation to the IMSE of $\widehat{f}^{(\vartheta)}(y|\bx)$. Employing Lemma \ref{lem:Equivalent Kernel Representation} and assuming that $\pOrder-\vartheta$ is even (as outlined in the local polynomial regression literature \citep{Fan-Gijbels_1996_Book}), we propose to select the bandwidth by minimizing a feasible analogue of the integrated mean squared error (IMSE)

\begin{align*}
    h^\star_\pOrder = \argmin_{h > 0}\iint_{\Ysupp \times \Xsupp} \left( h^{2\pOrder-2\vartheta}B_{\vartheta}(y,\bx)^2 + \frac{1}{nh^{1+2\vartheta + d}}V_{\vartheta}(y,\bx)  \right) \text{d}y\text{d}\bx,
\end{align*}
where $B_\vartheta(y,\bx)$ and $V_\vartheta(y,\bx)$ are the constants in the leading bias and variance, respectively, defined as
\begin{align*}
B_{\vartheta}(y,\bx) &= f^{(\pOrder)}(y|\bx)\be_{1+\vartheta}^\Trans\bS_y^{-1}\bc_{y,\pOrder+1} + \sum_{|\bnu|=\pOrder-\vartheta}\frac{\partial^{\bnu}}{\partial \bx^{\bnu}} f^{(\vartheta)}(y|\bx)\be_{\mathbf{0}}^\Trans\bS_{\bx}^{-1}\bc_{\bx,\bnu},\\
V_{\vartheta}(y,\bx) &= f(y|\bx)\Big(\be_{1+\vartheta}^\Trans\bS_y^\inv \bT_{y}\bS_y^\inv \be_{1+\vartheta}\Big)\Big(\be_{\mathbf{0}}^\Trans \bS_{\bx}^\inv  \bT_{\bx}\bS_{\bx}^\inv \be_{\mathbf{0}}\Big),
\end{align*}
with
\begin{align*}
\bc_{y,\pOrder+1} &= \int_\Ysupp \frac{1}{(\pOrder+1)!}\Big(\frac{u-y}{h}\Big)^{\pOrder+1}\frac{1}{h}\bP\Big(\frac{u-y}{h}\Big) \text{d}F_y(u),\quad
\bc_{\bx,\bnu} = \int_\Xsupp \frac{1}{\bnu!}\Big(\frac{\bv-\bx}{h}\Big)^{\bnu}\frac{1}{h^d}\bQ\Big(\frac{\bv-\bx}{h}\Big)  \mathrm{d}F_\bx(\bv).
\end{align*}
Both $B_\vartheta(y|\bx)$ and $V_\vartheta(y|\bx)$ involve the conditional PDF and its derivatives, which can be estimated with our proposed method. Other unknown quantities in the IMSE expression have the sample analogues:
\begin{align*}
    &\widehat{\bc}_{y,\pOrder+1} = \frac{1}{nh}\sumIN\frac{1}{(\pOrder+1)!}\left(\frac{y_{i}-y}{h}\right)^{\pOrder+1}\bP \left(\frac{y_{i}-y}{h}\right)^\Trans,\quad \hat{\bc}_{\bx,\bnu} = \frac{1}{nh^d}\sumIN\frac{1}{\bnu!}\left(\frac{\bx_i-\bx}{h}\right)^{\bnu}\bQ  \left(\frac{\bx_i-\bx}{h}\right),\\
    &\widehat{\bT}_{y} = \frac{1}{n^2h^3}\sum_{i,j=1}^n   \big(\min(y_i,y_j)-y\big)\bP\Big(\frac{y_i-y}{h}\Big)\bP\Big(\frac{y_j-y}{h}\Big)^\Trans,\quad \widehat{\bT}_{\bx} = \frac{1}{nh^d}\sumIN \bQ  \left(\frac{\bx_i-\bx}{h}\right)\bQ  \left(\frac{\bx_i-\bx}{h}\right)^\Trans.
\end{align*}

The bandwidth that minimizes the approximate IMSE, $h^\star_\pOrder$, is proportional to $ n^{-\frac{1}{1+d+2\pOrder}}$. Although this bandwidth delivers estimates that are approximately IMSE-optimal, a non-vanishing bias will be present in their asymptotic distribution, complicating statistical inference. To address this well-known problem, our construction of confidence bands and test statistics for parametric or shape restrictions employs robust bias correction \citep{Calonico-Cattaneo-Farrell_2018_JASA,Calonico-Cattaneo-Farrell_2022_Bernoulli}: one first constructs an IMSE-optimal point estimator, and then bias corrects the estimator and adjust the covariance function estimator accordingly to obtain a valid distributional approximation. More precisely, given an IMSE-optimal point estimator $\widehat{f}^{(\vartheta)}(y|\bx)$, robust bias correction relies on a test statistic of the form
\[\frac{\widehat{f}^{(\vartheta)}(y|\bx) - \widehat{\mathsf{Bias}}\big[\widehat{f}^{(\vartheta)}(y|\bx)\big]}
       {\sqrt{\widehat{\mathsf{Var}}\Big[\widehat{f}^{(\vartheta)}(y|\bx) - \widehat{\mathsf{Bias}}\big[\widehat{f}^{(\vartheta)}(y|\bx)\big]\Big]}},
\]
where $\widehat{\mathsf{Bias}}[\widehat{f}^{(\vartheta)}(y|\bx)]$ denotes a bias correction estimate of the IMSE-optimal point estimator $\widehat{f}^{(\vartheta)}(y|\bx)$, and $\widehat{\mathsf{Var}}\big[\widehat{f}^{(\vartheta)}(y|\bx) - \widehat{\mathsf{Bias}}[\widehat{f}^{(\vartheta)}(y|\bx)]\big]$ denotes an estimator of the variance of the bias-corrected estimate. The key idea underlying robust bias correction is to Studentize by the variance of the bias corrected estimate as opposed to by the variance of the original point estimator, an approach that leads to better distributional approximations \citep{Calonico-Cattaneo-Farrell_2018_JASA,Calonico-Cattaneo-Farrell_2022_Bernoulli}. Similarly, uniform robust bias correction constructs an estimator of the correlation function $\rho_{\vartheta}(y,\bx,y',\bx')$ taking into account the additional variability introduced by the bias correction.

A simple and intuitive way of operationalizing robust bias correction in local
polynomial settings is by increasing the polynomial order $\pOrder$  (recall that we set $\qOrder = \pOrder - \vartheta - 1$). That is, we first compute the bandwidth $h^\star_\pOrder$, and then form the final estimator with a local polynomial order of $\pOrder+1$. To make the procedure precise, we augment the notation so that it reflects the local polynomial order and the bandwidth used as needed. For example, the conditional density estimator using polynomial order $\pOrder$ and employing the bandwidth $h$ is written as $\widehat{f}^{(\vartheta)}_\pOrder(y|\bx;h)$.

\subsection{Confidence bands}

Confidence bands can be constructed using the process $(\widehat{\TProc}_{\vartheta,\pOrder+1}^{\mathtt{CB}}(y,\bx) : y\in\mathcal{Y},\bx\in\mathcal{X})$, where
\begin{align*}
    \widehat{\TProc}_{\vartheta,\pOrder+1}^{\mathtt{CB}}(y,\bx) &= \frac{\widehat{f}^{(\vartheta)}_{\pOrder+1}(y|\bx;h^\star_\pOrder) - f^{(\vartheta)}(y|\bx)}{\sqrt{\widehat{\mathsf{V}}_{\vartheta,\pOrder+1}(y,\bx;h^\star_\pOrder)}},
\end{align*}
By Theorem \ref{thm:Kolmogorov-Smirnov Distance: Suprema}, the distribution of $\sup_{y\in\mathcal{Y},\bx\in\mathcal{X}}|\widehat{\TProc}_{\vartheta,\pOrder+1}^{\mathtt{CB}}(y,\bx)|$ is approximated by the conditional (on the data) distribution of $\sup_{y\in\mathcal{Y},\bx\in\mathcal{X}}|\widehat{\GProc}_{\vartheta,\pOrder+1}(y,\bx)|$, with $\widehat{\GProc}_{\vartheta,\pOrder+1}$ being
a centered Gaussian process whose law, conditionally on the data, is Gaussian with unit variance and correlation $\widehat{\rho}_{\vartheta,\pOrder+1}(\cdot ;h^\star_\pOrder)$. Accordingly, let
\[ \mathrm{CB}_{\vartheta,\pOrder+1}(1-\alpha) = \Big[ \; \widehat{f}_{\pOrder+1}^{(\vartheta)}(y|\bx;h^\star_{\pOrder}) \ \pm\ \cval_{\vartheta,\pOrder+1}^{\mathtt{CB}}(\alpha) \sqrt{\widehat{\mathsf{V}}_{\vartheta,\pOrder+1}(y,\bx;h^\star_{\pOrder})} \; : \; y\in\mathcal{Y},\bx\in\mathcal{X} \Big],\]
where
\[\cval_{\vartheta,\pOrder+1}^{\mathtt{CB}}(\alpha)=\inf\Big\{u\in\mathbb{R}_+:\Prob\Big[ \sup_{y\in\mathcal{Y},\bx\in\mathcal{X}} \big|\widehat{\GProc}_{\vartheta,\pOrder+1}(y,\bx)\big| \leq u \;\Big| \;\bX,\bY \Big]\geq 1-\alpha \Big\}.\]

As the notation suggests, $\mathrm{CB}_{\vartheta,\pOrder+1}(1-\alpha)$ is a
$100(1-\alpha)\%$ confidence band. To be specific, we have the following theorem.
\begin{theorem}[Confidence bands]\label{thm:feasible Coverage Error Band}
Suppose Assumptions \ref{as:1-DGP} and \ref{as:2-Kernel} hold, $f^{(\pOrder+1)}(y|\bx)$ exists and is continuous, and $\partial^{\bnu} f^{(\vartheta)}(y|\bx)/\partial \bx^{\bnu} $ exists and is continuous for all $|\bnu|=\pOrder+1-\vartheta$. Then
\begin{align*}
 \left| \Prob\left[ f^{(\vartheta)} \in \mathrm{CB}_{\vartheta,\pOrder+1}(1-\alpha) \right] - (1-\alpha) \right| \precsim \log^{\frac{5}{4}}(n)\mathtt{r_{CB}},
\end{align*}
where $\mathtt{r_{CB}} = n^{-\frac{1}{1+d+2\pOrder}} +n^{-\frac{2\pOrder-2\vartheta +1}{4(1+d+2\pOrder)}} + n^{-\frac{\pOrder}{(1+d+2\pOrder)(1+d)}}$.
\end{theorem}

The confidence band $\mathrm{CB}_{\vartheta,\pOrder+1}(1-\alpha)$ is easy to construct because, by discretizing the index set of the Gaussian process, the critical value $\cval_{\vartheta,\pOrder+1}(1-\alpha)$ can be computed by simulation from a conditionally (on the data) multivariate Gaussian distribution. We illustrate the performance of our proposed confidence bands using simulated and real data in Section \ref{section: Simulations}.

Theorem \ref{thm:feasible Coverage Error Band} provides a formal, theoretical justification for employing strong approximation methods to construct confidence bands instead of relying on extreme value theory for approximating the distribution of the suprema of the process $\widehat{\TProc}_{\vartheta,\pOrder+1}^{\mathtt{CB}}$. More specifically, the coverage error rate $\mathtt{r_{CB}}$ is polynomial in $n$ for the former inference approach, while the latter inference approach would have a logarithmic in $n$ convergence rate \citep[see, e.g.,][and references therein]{Hall_1979_JAP,Hall_1993_JRSSB}. The same remark applies to the upcoming Theorems \ref{thm:Specification Testing} and \ref{thm:Shape Restrictions}, which characterize the error in rejection probability of two different classes of hypothesis testing procedures.

\subsection{Parametric specification testing}

Suppose the researcher postulates that the conditional density (derivative) belongs to the parametric class $\{ f^{(\vartheta)}(y|\bx;\bgamma): \bgamma\in \mathsf{\Gamma}_{\vartheta} \}$, where $\mathsf{\Gamma}_{\vartheta}$ is some parameter space. Abstracting away from the specifics of the estimation technique, we assume that the researcher also picks some estimator $\widehat{\bgamma}$ (e.g., maximum likelihood or minimum distance), which is assumed to converge in probability to some  $\bar{\bgamma}\in\mathsf{\Gamma}_{\vartheta}$. A natural statistic for the problem of testing
\begin{align*}
\mathsf{H}_0^{\mathtt{PS}}&:\ f^{(\vartheta)}(y|\bx;\bar{\bgamma})=f^{(\vartheta)}(y|\bx) \qquad \text{for all } (y,\bx)\in\Ysupp\times\Xsupp
\end{align*}
is
\begin{align*}
\sup_{y \in \Ysupp, \bx \in \Xsupp}\big|\widehat{\TProc}_{\vartheta,\pOrder+1}^{\mathtt{PS}} (\eval) \big|, \qquad
\widehat{\TProc}_{\vartheta,\pOrder+1}^{\mathtt{PS}}(\eval) = \frac{\widehat{f}^{(\vartheta)}_{\pOrder+1}(y|\bx;h^\star_{\pOrder}) - f^{(\vartheta)}(y|\bx;\widehat{\bgamma}) }{\sqrt{\widehat{\mathsf{V}}_{\vartheta,\pOrder+1}(y,\bx;h^\star_{\pOrder})}}.
\end{align*}
Assuming the estimation error of $\widehat{\bgamma}$ is asymptotically negligible, a valid $100\alpha\%$ critical value is given by $\cval_{\vartheta,\pOrder+1}^{\mathtt{CB}}(\alpha)$. To be specific, we have:

\begin{theorem}[Parametric specification testing]\label{thm:Specification Testing}
Suppose Assumptions \ref{as:1-DGP} and \ref{as:2-Kernel} hold, $f^{(\pOrder+1)}(y|\bx)$ exists and is continuous, and $\partial^{\bnu} f^{(\vartheta)}(y|\bx)/\partial \bx^{\bnu} $ exists and is continuous for all $|\bnu|=\pOrder+1-\vartheta$. If
\begin{align*}
n^{\frac{\pOrder-\vartheta}{1+d+2\pOrder}}\sup_{ y \in \Ysupp, \bx \in \Xsupp}\left| f^{(\vartheta)}(y|\bx;\widehat{\bgamma}) - f^{(\vartheta)}(y|\bx;\bar{\bgamma}) \right| \ \precsimConc\ \mathtt{r_{CB}},
\end{align*}
then, under $\mathsf{H}_0^{\mathtt{PS}}$,
\begin{align*}
\Big|\Prob\Big[ \sup_{y \in \Ysupp, \bx \in \Xsupp}|\widehat{\TProc}_{\vartheta,\pOrder+1}^{\mathtt{PS}}(\eval) | > \cval_{\vartheta,\pOrder+1}^{\mathtt{CB}}(\alpha) \Big] - \alpha\Big| \precsim
\log^{\frac{5}{4}}(n)\mathtt{r_{CB}},
\end{align*}
where $\mathtt{r_{CB}}$ is defined in Theorem~\ref{thm:feasible Coverage Error Band}.
\end{theorem}

\subsection{Testing shape restrictions}

As a third application, suppose the researcher wants to test shape restrictions on $f^{(\vartheta)}$. Letting $c_{\vartheta}$ be a pre-specified function, consider the problem of testing
\begin{align*}
    \mathsf{H}_0^{\mathtt{SR}} &:f^{(\vartheta)}(y|\bx) \leq c_{\vartheta}(y|\bx) \qquad \text{for all } (y,\bx)\in\Ysupp\times\Xsupp.
\end{align*}
For example, if $\vartheta=0$ and if $c_{\vartheta}(y|\bx)$ is some (positive) constant value $c$, the testing problem refers to whether the conditional density exceeds $c$ somewhere on its support. As another example, if $\vartheta=1$ and if $c_{\vartheta}(y|\bx)=0$, then the testing problem refers to whether the conditional density is non-increasing in $y$ for all values of $\bx$. More generally, the testing problem above can be used to test for monotonicity, convexity, and other shape features of the conditional density, possibly relative to the function $c_{\vartheta}(y|\bx)$.

A natural testing procedure rejects $\mathsf{H}_0^{\mathtt{SR}}$ whenever the test statistic
\[
\sup_{y \in \Ysupp, \bx \in \Xsupp} {\TProc}_{\vartheta,\pOrder+1}^{\mathtt{SR}}(\eval), \qquad
{\TProc}_{\vartheta,\pOrder+1}^{\mathtt{SR}}(\eval) = \frac{\widehat{f}^{(\vartheta)}_{\pOrder+1}(y|\bx;h^\star_{\pOrder}) - c_{\vartheta}(y|\bx)}{\sqrt{\widehat{\mathsf{V}}_{\vartheta,\pOrder+1}(y,\bx;h^\star_{\pOrder})}}
\]
exceeds a critical value of the form
\[\cval_{\vartheta,\pOrder+1}^{\mathtt{SR}}(\alpha)=\inf\Big\{u\in\mathbb{R}_+:\Prob\Big[ \sup_{y\in\mathcal{Y},\bx\in\mathcal{X}} \widehat{\GProc}_{\vartheta,\pOrder+1}(y,\bx) \leq u \;\Big| \;\bX,\bY \Big]\geq 1-\alpha \Big\}.\]

\begin{theorem}[Testing shape restriction]\label{thm:Shape Restrictions}
Suppose Assumptions \ref{as:1-DGP} and \ref{as:2-Kernel} hold, $f^{(\pOrder+1)}(y|\bx)$ exists and is continuous, and $\partial^{\bnu} f^{(\vartheta)}(y|\bx)/\partial \bx^{\bnu} $ exists and is continuous for all $|\bnu|=\pOrder+1-\vartheta$. Then, under $\mathsf{H}_0^{\mathtt{SR}}$,
\[
\Big|\Prob\Big[ \sup_{y \in \Ysupp, \bx \in \Xsupp}\widehat{\TProc}_{\vartheta,\pOrder+1}^{\mathtt{SR}}(\eval)  > \cval_{\vartheta,\pOrder+1}^{\mathtt{SR}}(\alpha) \Big] - \alpha\Big| \precsim
\log^{\frac{5}{4}}(n)\mathtt{r_{CB}},
\]
where $\mathtt{r_{CB}}$ is defined in Theorem~\ref{thm:feasible Coverage Error Band}.
\end{theorem}

\section{Imposing additional constraints for density estimation}\label{section:additional constraints}

Specific applications may require additional constraints on the estimates. For example, setting $\vartheta=0$ (PDF), it may be desirable to require that the estimator is non-negative and integrates to one. The nonnegativity constraint can be directly incorporated into the local polynomial regression \eqref{eq:local reg along y direction}:
\begin{align*}
\widehat{f}_{\mathtt{N}}(y|\bx) = \be_{1}^\Trans\widehat{\bbeta}_{\mathtt{N}}(y|\bx), \qquad
  \widehat{\bbeta}_{\mathtt{N}}(y|\bx) = \argmin_{\substack{\bu\in\mathbb{R}^{\pOrder+1}:\  \be_{1}^\Trans\bu\geq 0}} \sum_{i=1}^n \left(\widehat{F}(y_i|\bx) - \bp(y_i-y)^\Trans\bu \right)^2 K_h(y_i;y),
\end{align*}
where the subscript ``\texttt{N}'' stands for ``non-negative.'' While $\widehat{f}_{\mathtt{N}}(y|\bx)$ is non-negative by construction, it does not necessarily integrate to one. This follows from the fact that the estimator only exploits local features of the data and not global constraints. To address the second constraint, we propose and study the following enhanced estimator based on minimizing Kullback-Leibler divergence (the subscript ``\texttt{I}'' stands for ``integrating to one''):
\begin{align*}
    \widehat{f}_{\texttt{I}}(y|\bx) = \argmin_{g\in\mathcal{G}} \text{KL}\big(g \; \big\| \; \widehat{f}_{\mathtt{N}}(\cdot|\bx)\big), \qquad
    &\text{where }\text{KL}(g \; \big\| \; f) = \int_{\Ysupp} g(y) \log\Big(\frac{g(y)}{f(y)}\Big) \text{d}y,
\end{align*}
and $\mathcal{G} = \{g \geq 0: \int_{\Ysupp} g(y) \text{d}y = 1,\ g(y)=0\text{ for }y\not\in\mathcal{Y}\}$. It follows that our proposed conditional PDF estimator, $\widehat{f}_{\texttt{I}}(y|\bx)$, is non-negative and integrates to one. Furthermore, both $\widehat{f}_{\mathtt{N}}(y|\bx)$ and $\widehat{f}_{\mathtt{I}}(y|\bx)$ can be written in closed form (see Appendix \ref{appendix: derivation of constrained estimators}):
\begin{align}\label{eq:the constrained estimators}
\widehat{f}_{\mathtt{I}}(y|\bx) = \frac{\widehat{f}_{\mathtt{N}}(y|\bx)}{\int_{\Ysupp} \widehat{f}_{\mathtt{N}}(u|\bx) \text{d}u} \qquad\text{and}\qquad
  \widehat{f}_{\mathtt{N}}(y|\bx) = \max\big\{ \widehat{f}(y|\bx)\ ,\ 0 \big\}.
\end{align}
In practice, the support $\mathcal{Y}$ might be unknown, and in this case one can naturally replace it by the empirical support: $\widehat{\Ysupp} = [y_{(1)},y_{(n)}]$, defined by the smallest ($y_{(1)}$) and largest ($y_{(n)}$) order statistics of the observed $y_1,y_2,\dots,y_n$. Since $\widehat{\Ysupp}\subseteq \mathcal{Y}$, all the theoretical results discussed below remain valid on the empirical support $\widehat{\Ysupp}$.

We first establish stochastic linearization for both, $\widehat{f}_{\mathtt{N}}(y|\bx)$ and $\widehat{f}_{\mathtt{I}}(y|\bx)$.

\begin{lem}[Stochastic linearization]\label{lem:Equivalent Kernel Representation, NN and I estimators}
    Suppose Assumptions \ref{as:1-DGP} and \ref{as:2-Kernel} hold. If $nh^{1+d} / \log (n) \to \infty$ and $h\to 0$, then
   \begin{align*}
    &\sup_{y\in\Ysupp,\bx\in\Xsupp} \left| \widehat{f}_{\mathtt{N}}(y|\bx) - f(y|\bx) - \bar{f}^{(0)}(y|\bx) \right|\ \precsimConc\ \mathtt{r_{SL}},\\
    \text{and }\qquad&\sup_{y\in\Ysupp,\bx\in\Xsupp} \Big| \widehat{f}_{\mathtt{I}}(y|\bx) - f(y|\bx) - \Big(\bar{f}^{(0)}(y|\bx) - f(y|\bx)\int_\Ysupp \bar{f}^{(0)}(u|\bx) \diff u\Big)\Big|\ \precsimConc\ \mathtt{r_{SL}},
   \end{align*}
    where $\bar{f}^{(0)}(y|\bx)$ and $\mathtt{r_{SL}}$ are defined in Lemma \ref{lem:Equivalent Kernel Representation} by setting $\vartheta = 0$.
\end{lem}

The lemma provides a more refined stochastic linearization for $\widehat{f}_{\mathtt{I}}(y|\bx)$. We will show that the normalization in $\widehat{f}_{\mathtt{I}}(y|\bx)$ does not affect the uniform rate of convergence of the estimator. For distributional approximation, however, it is crucial to employ different Gaussian processes for the two estimators. In particular, we show that failing to capture the asymptotic contribution of the normalization in $\widehat{f}_{\mathtt{I}}(y|\bx)$ may lead to a slower rate for strong approximation.

\begin{theorem}[Probability concentration]\label{thm:Probability Concentration, NN and I estimators}
Suppose Assumptions \ref{as:1-DGP} and \ref{as:2-Kernel} hold. If $h\to 0$ and if $nh^{1+d}/\log (n)\to\infty$, then
\[
\sup_{y\in \Ysupp, \bx \in \Xsupp}\left| \widehat{f}_{\mathtt{N}}(y|\bx) - f(y|\bx) \right|\ \precsimConc\ \mathtt{r_{PC}},\quad
  \sup_{y\in \Ysupp, \bx \in \Xsupp}\left| \widehat{f}_{\mathtt{I}}(y|\bx) - f(y|\bx) \right|\ \precsimConc\ \mathtt{r_{PC}},
\]
where $\mathtt{r_{PC}}$ is defined in Theorem \ref{thm:Probability Concentration} (with $\vartheta=0$).
\end{theorem}

Finally, to state a strong approximation result, we define the following standardized processes
\[
    \widehat{\SProc}_{\mathtt{N}}(y,\bx) = \frac{\widehat{f}_{\mathtt{N}}(y|\bx) - f(y|\bx)}{\sqrt{\mathsf{V}_{0}(y,\bx)}},\quad \widehat{\SProc}_{\mathtt{I}}(y,\bx) = \frac{\widehat{f}_{\mathtt{I}}(y|\bx) - f(y|\bx)}{\sqrt{\mathsf{V}_{0}(y,\bx)}}.
\]
\begin{theorem}[Strong approximation]\label{thm:Strong Approximation S-hat, NN and I estimators}
    Suppose Assumptions \ref{as:1-DGP} and \ref{as:2-Kernel} hold. If ${nh^{1+d+2\pOrder}}\to 0$ and if $nh^{1+d}/\log (n)\to\infty$, then there exist three stochastic processes, $\widehat{\SProc}^{\prime}_{\mathtt{N}}$, $\widehat{\SProc}^{\prime}_{\mathtt{I}}$, and $\GProc$, in a possibly enlarged probability space, such that:
    \begin{enumerate}[label=(\roman*),noitemsep]
        \item $\widehat{\SProc}_{\mathtt{N}}$ and $\widehat{\SProc}^{\prime}_{\mathtt{N}}$ have the same distribution; $\widehat{\SProc}_{\mathtt{I}}$ and $\widehat{\SProc}^{\prime}_{\mathtt{I}}$ have the same distribution
        \item $\GProc$ is a centered Gaussian process with unit variance and correlation $\rho_{0}$;
        \item the following holds:
			\begin{align*}
             &\sup_{y \in \Ysupp, \bx \in \Xsupp}\left|\widehat{\SProc}^{\prime}_{\mathtt{N}}(y,\bx)-\GProc(y,\bx)\right|\ \precsimConc\  \mathtt{r_{SA}},\\
             \text{and}\quad &\sup_{y \in \Ysupp, \bx \in \Xsupp}\left|\widehat{\SProc}^{\prime}_{\mathtt{I}}(y,\bx)- \Big(\GProc(y,\bx) - f(y|\bx)\int_{\Ysupp} \sqrt{\frac{\mathsf{V}_{0}(u,\bx)}{\mathsf{V}_{0}(y,\bx)}}\GProc(u,\bx)\diff u\Big)\right| \precsimConc \mathtt{r_{SA}},
              \end{align*}
         where $\mathtt{r_{SA}}$ is defined in Theorem \ref{thm:Strong Approximation S-hat}.
    \end{enumerate}
\end{theorem}
The different Gaussian processes needed for distributional approximation to
$\widehat{\SProc}_{\mathtt{N}}$ and $\widehat{\SProc}_{\mathtt{I}}$ in
Theorem~\ref{thm:Strong Approximation S-hat, NN and I estimators} is due to the
normalization in $\widehat{\SProc}_{\mathtt{I}}$. Of course, it is possible to
couple $\widehat{\SProc}^{\prime}_{\mathtt{I}}$ with $\GProc$ directly, but a
slower rate may arise, particularly $\sup_{y \in \Ysupp, \bx \in \Xsupp}|\widehat{\SProc}^{\prime}_{\mathtt{I}}- \GProc(y,\bx)|\ \precsimConc \mathtt{r_{SA}} + \sqrt{\log(n)h}$.

Constructing analogues of Lemma \ref{lem:Covariance Estimation} and Theorem \ref{thm:Kolmogorov-Smirnov Distance: Suprema} from Section \ref{subsection:Suprema Approximation} for the constrained estimators $\widehat{f}_{\mathtt{N}}(y|\bx)$ and $\widehat{f}_{\mathtt{I}}(y|\bx)$ now follows directly. Additionally, confidence bands and hypothesis testing procedures as in Section \ref{section:Applications} can also be easily developed when employing the constrained density estimators. We omit details to avoid repetition.

\section{Numerical evidence}\label{section: Simulations}

We illustrate the effectiveness of our proposed methods with two Monte Carlo
experiments, where we set $d=1$ and simulate $\bx$ and $y$ from a joint normal distribution with variance $2$ and covariance
$-0.1$, truncated on $[-1, 1]^2$.
We employ 1000 Monte Carlo repetitions, each with the sample size $n=5000$. Replication files, additional simulation results, and details of the
companion \texttt{R} package, \texttt{lpcde}, can be found at
\url{https://nppackages.github.io/lpcde/} and in our companion software article
\citep{Cattaneo-Chandak-Jansson-Ma_2022_lpcde}.

In the first simulation experiment, we estimate
the conditional PDF for 20 equally spaced points on $[-1, 1]$ for $y$.
Table \ref{table:coverage} presents the simulation results
at three different conditioning values: (a) interior ($\bx=0$),
(b) near-boundary ($\bx=0.8$), and (c) at-boundary ($\bx=1$).
See the discussion at the end of Section \ref{section:Introduction} for a classification of  interior and (near) boundary evaluation points.

Table~\ref{table:coverage} reports average estimated bandwidth in column
``$\widehat{h}$'',
and average bias and standard error in the
``bias'' and ``se'' columns, respectively.
We consider bands formed by pointwise confidence intervals (columns
``pointwise''), which are not uniformly valid and hence should exhibit
considerable under coverage, as well as the uniform confidence bands discussed
in Section \ref{section:Applications} (columns ``uniform'').
We report their empirical uniform coverage probabilities (column ``Coverage'')
and the average width (column ``Width'').
For the non-bias corrected rows (``NBC''), the polynomial orders for bandwidth
selection, point estimation and statistical inference are $\pOrder = 2$ and
$\qOrder=1$, while those for robust bias-corrected statistical inference rows
(``RBC'') are $\pOrder=3$ and $\qOrder=2$.
\begin{table}[!ht]
  \caption{Empirical uniform coverage probabilities. }
\label{table:coverage}
\begin{tabular}{lr|ccc|cc|cc}
\hline & & & & & \multicolumn{2}{c|} { Coverage } & \multicolumn{2}{c} { Width } \\
& & $\widehat{h}$ & bias & se & pointwise & uniform & pointwise & uniform \\
\hline
  \multirow{2}{*}{$\bx=0$} & NBC &0.32 & 0.09 & 0.03 & 62.6 & 74.8 & 0.01 & 0.02\\
              & RBC &0.32  & 0.09 & 0.09& 83.4 & 93.9 & 0.05 & 0.05\\
  \hline
  \multirow{2}{*}{$\bx=0.8$} & NBC &0.30  & 0.10 & 0.04 & 72.8 & 89.4 & 0.02 & 0.03 \\
              & RBC &0.30 & 0.10 & 0.18 & 86.9 & 94.3 & 0.13 & 0.19 \\
  \hline
  \multirow{2}{*}{$\bx=1.0$} & NBC &  0.32 & 0.10 & 0.06 & 74.9 & 91.3 & 0.02 & 0.05 \\
              & RBC & 0.32 & 0.10 & 0.20 & 88.1 & 93.2 & 0.11 & 0.23 \\
   \hline
\end{tabular}%
\end{table}

\begin{table}[ht]
\centering
  \caption{Comparison between \texttt{FYT} and our method for conditional
    PDF estimation.
}
  \label{table:fyt_comparison}
  \begin{tabular}{r|cccc|cccccc}
    \hline
    & \multicolumn{4}{c|}{\texttt{FYT}} & \multicolumn{6}{c}{This paper: $\widehat{f}(y|\bx)$}\\
    \hline
    $\times\ h_{\text{MSE}}$ & $h$ & bias & se & rmse & $h$ & bias & se & rmse  & NBC CI & RBC CI\\
    \hline
    ($y=0\phantom{.0}, \bx=0$)\qquad 0.8 & 0.32 &0.08 & 0.13 & 0.15 & 0.23 & 0.07 & 0.07 & 0.11 & 0.87 & 0.99 \\
    0.9 & 0.36 & 0.09 & 0.12 & 0.15 & 0.26 & 0.07 & 0.05 & 0.09& 0.74 & 0.98 \\
    1 & 0.39 & 0.09 & 0.12 & 0.15 & 0.29 &0.07 & 0.04 & 0.08& 0.56 & 0.96 \\
    1.1& 0.43 & 0.09 & 0.12 & 0.15 & 0.32 & 0.06 & 0.03 & 0.07& 0.39 & 0.89 \\
    1.2& 0.46 & 0.10 & 0.11 & 0.15 & 0.35 & 0.07 & 0.02 & 0.07& 0.25 & 0.81 \\
    \hline
    ($y=0.8, \bx=0$)\qquad 0.8 & 0.29 & 0.14 & 0.10 & 0.18 & 0.26 & 0.03 & 0.02 & 0.04 & 0.90 & 0.99 \\
    0.9 & 0.33 & 0.14 & 0.10 & 0.17 & 0.30 & 0.03 & 0.01 & 0.03 & 0.83 & 0.98 \\
    1 & 0.36 & 0.14 & 0.09 & 0.17 & 0.33 & 0.03 & 0.01 & 0.03 & 0.75 & 0.93 \\
    1.1 & 0.39 & 0.14 & 0.09 & 0.17 & 0.36 & 0.03 & 0.01 & 0.04 & 0.70 & 0.87 \\
    1.2 & 0.43 & 0.14 & 0.08 & 0.16 & 0.40 & 0.03 & 0.01 & 0.04 & 0.64 & 0.80 \\
    \hline
    ($y=1\phantom{.0}, \bx=0$)\qquad 0.8 & 0.27 & 0.18 & 0.07 & 0.20 & 0.40 & 0.04 & 0.04 & 0.06 & 0.93 & 1.00 \\
    0.9 & 0.30 & 0.18 & 0.07 & 0.19 & 0.45 & 0.04 & 0.03 & 0.05 & 0.73 & 0.99 \\
    1 & 0.33 & 0.20 & 0.06 & 0.20& 0.50 & 0.04 & 0.02 & 0.04 & 0.51 & 0.96 \\
    1.1 & 0.36 & 0.21 & 0.06 & 0.21& 0.55 & 0.04 & 0.01 & 0.04 & 0.36 & 0.89 \\
    1.2 & 0.39 & 0.23 & 0.05 & 0.24& 0.60 & 0.04 & 0.01 & 0.04 & 0.22 & 0.80 \\
    \hline
  \end{tabular}
\end{table}

The simulation results in Table \ref{table:coverage} support our main
theoretical findings.
First, robust bias correction leads to uniformly better performance
of the inference procedures, both pointwise and uniformly over $\mathcal{Y}$.
Second, our uniform distributional approximation leads to feasible confidence
bands with good finite sample performance, when coupled with robust bias
correction methods.

For example, for $\bx=0$, the averaged (across simulations) estimated
approximate IMSE-optimal bandwidth choice is $\widehat{h}=0.32$, with
$\pOrder=2$ and $\qOrder=\pOrder-1$. Bands constructed with pointwise confidence
intervals have empirical uniform coverage of $62.6\%$ without bias correction,
and $83.4\%$ with robust bias correction, both are substantially below the
$95\%$ nominal level because they are not uniformly valid over the range of $y$.
The feasible confidence bands are designed to address that issue: our proposed
confidence bands have empirical coverage of $93.9\%$ when robust bias correction
is employed. It also highlights the importance of addressing the
misspecification (smoothing) bias for statistical inference. Without bias
correction, the uniform confidence bands only cover the true conditional PDF with
probability $74.8\%$.

The second simulation study compares
our estimator (\texttt{lpcde})
to the estimator proposed by
\citet{Fan-Yao-Tong_1996_Biometrika} (\texttt{FYT},
see Appendix \ref{Appendix: closed-form expression and comparison with FYT}
for details).
Table \ref{table:fyt_comparison} presents the simulation results for the
conditional PDF at three distinct evaluation points.
For a fair comparison, we first compute the MSE optimal bandwidth ($h_{\text{MSE}}$) for the two estimators at each evaluation point. We then investigate the performance of the two estimators over a grid of bandwidths, ranging from $0.8\times h_{\text{MSE}}$ (under smoothing) to $1.2 \times h_{\text{MSE}}$ (over-smoothing).

For each of the two estimators we report the average bandwidth,
bias, standard error, and root mean squared error.
Additionally, for our estimator we report the pointwise empirical coverage probabilities, both with and without bias correction.
Since \texttt{FYT} do not provide theory for statistical inference,
we do not report confidence interval information for the estimator.
Results in Table~\ref{table:fyt_comparison} suggest that our local polynomial conditional PDF estimator perform well across all three evaluation points, and the confidence intervals constructed thereof exhibits satisfactory empirical coverage property. In particular, at the boundary evaluation point ($y=1, \bx=0$), our estimator has
accurate coverage while \texttt{FYT} suffers from boundary bias.

\begin{figure}[h]
  \begin{subfigure}[b]{0.5\linewidth}
    \centering
    \includegraphics[width=0.95\linewidth]{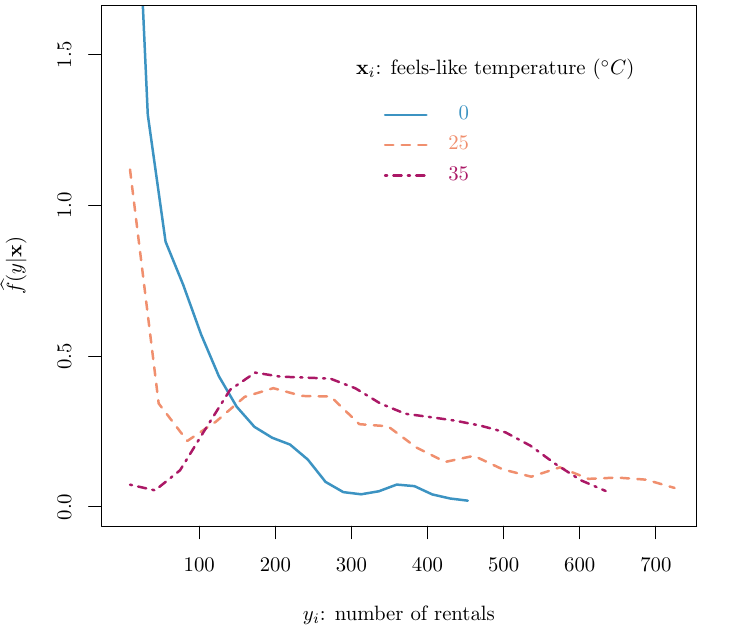}
    \caption{Estimated conditional PDFs.}
    \label{fig:pdf_std}
  \end{subfigure}
  \begin{subfigure}[b]{0.5\linewidth}
    \centering
    \includegraphics[width=0.95\linewidth]{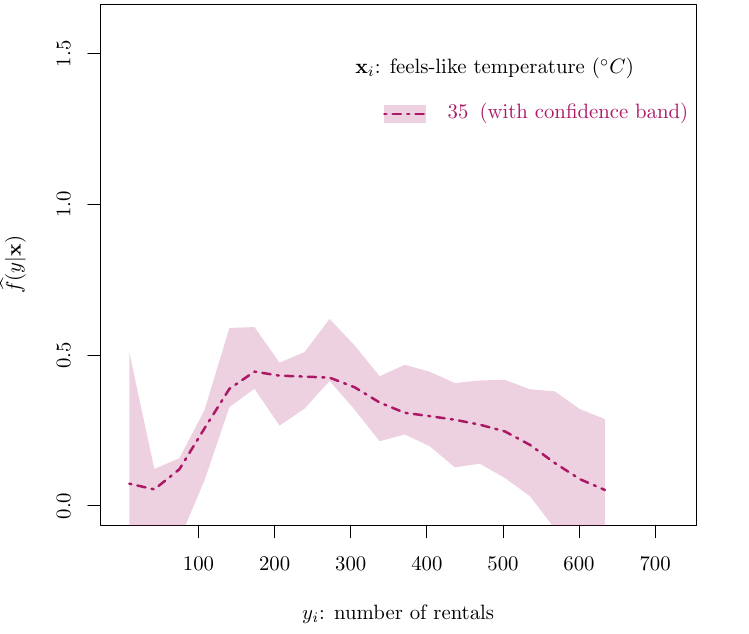}
    \caption{Illustration of confidence band.}
    \label{fig:pdf_rbc}
  \end{subfigure}
  \caption{Estimated relationship (conditional PDF) between bike rental counts and temperature.}
  \label{fig:point_est}
\end{figure}

Finally, we illustrate the performance of our estimator in Figure~\ref{fig:point_est} with real data. The data we employ is from Capital Bikeshare (available at \url{https://archive.ics.uci.edu/dataset/275/bike+sharing+dataset}). The outcome variable $y_i$ is the total number of bike rentals, and the covariate $\bx_i$ is the ``feels-like'' temperature in Celsius. Panel (a) shows the estimated conditional PDFs for three temperature levels, $\bx_i = 0$, $25$, and $35$ $^\circ C$. From the conditional density plots, more bike rental activities  happen in warmer days (i.e., the conditional distribution moves toward right). It is worth mentioning that the outcome variable has a lower boundary at $0$, and using a standard kernel density estimator for conditional PDF estimation will lead to a severe under-estimation bias for $f(y|\bx)$ whenever the evaluation point $y$ is close to zero. To avoid overcrowding the figure, we illustrate the confidence band with robust bias correction in panel (b).

\section{Conclusion}\label{section:conclusion}

We introduced a new boundary adaptive estimator of the conditional density and derivatives thereof. This estimator is conceptually distinct from prior proposals in the literature, as it relies on two (nested) local polynomial estimators. Our proposed estimation approach has several appealing features, most notably automatic boundary adaptivity. We provided an array of uniform estimation and distributional results, including a valid uniform equivalent kernel representation and uniform distributional approximations. Our methods are applicable in data science settings either where the conditional density or its derivatives are the main object of interest, or where they are preliminary estimands entering a multi-step statistical procedure.

\begin{appendix}
\section*{Appendix}\label{appn}

\subsection{Derivation of \texorpdfstring{\eqref{equation:closed-form expression} } { }and an alternative expression}\label{Appendix: closed-form expression and comparison with FYT}

To start, the conditional CDF estimation step is a weighted least squares problem, and has the solution
\[
\widehat{F}(y_j|\bx) = \be_\mathbf{0}^\Trans\Big(\sum_{i=1}^{n}\bq(\bx_i-\bx)\bq(\bx_i-\bx)^\Trans L_b(\bx_i;\bx)\Big)^{-1}\Big(\sum_{i=1}^{n}\bq(\bx_i-\bx) L_b(\bx_i;\bx) \Indicator(y_i\leq y_j)\Big).
\]
The second local polynomial regression takes $\widehat{F}(y_j|\bx)$ as the ``dependent variable,'' and therefore the final estimator takes the form
\[
\widehat{f}^{(\vartheta)}(y|\bx) = \be_{1+\vartheta}^\Trans\Big(\sum_{j=1}^{n}\bp(y_j-y)\bp(y_j-y)^\Trans K_h(y_j;y)\Big)^{-1}\Big(\sum_{j=1}^{n}\bp(y_j-y) K_h(y_j;y) \widehat{F}(y_j|\bx)\Big).
\]
The final expression in \eqref{equation:closed-form expression} then follows from re-normalizing $\bx_i-\bx$ to $(\bx_i-\bx)/b$ and $y_j-y$ to $(y_j-y)/h$, leading to the multiplicative factor $h^{-1-\vartheta}$. By changing the order of summation in $\widehat{\bR}_{y,\bx}$, we can also write $\widehat{f}^{(\vartheta)}(y|\bx)$ as
\[\widehat{f}^{(\vartheta)}(y|\bx) = \be_\mathbf{0}^\Trans\Big(\sum_{i=1}^{n}\bq(\bx_i-\bx)\bq(\bx_i-\bx)^\Trans L_b(\bx_i;\bx)\Big)^{-1}\Big(\sum_{i=1}^{n}\bq(\bx_i-\bx) L_b(\bx_i;\bx) \widehat{K}_h(y_i,y)\Big),\]
where
\[\widehat{K}_h(y_i,y) = \be_{1+\vartheta}^\Trans\Big(\sum_{j=1}^{n}\bp(y_j-y)\bp(y_j-y)^\Trans K_h(y_j;y)\Big)^{-1}\Big(\sum_{j=1}^{n}\bp(y_j-y) K_h(y_j;y)\Indicator(y_i\leq y_j) \Big).
\]
The above alternative expression shows that our proposed estimator can be understood as first forming $\widehat{K}_h(y_i,y)$, which is a data-driven kernel re-weighting of $y_i$ and then conducting local polynomial regression on $\bx_i$. To compare, the density estimator ($\vartheta=0$) introduced by \citet{Fan-Yao-Tong_1996_Biometrika} takes the form
\[
\widehat{f}_{\mathtt{FYT}}(y|\bx) = \be_\mathbf{0}^\Trans\Big(\sum_{i=1}^{n}\bq(\bx_i-\bx)\bq(\bx_i-\bx)^\Trans L_b(\bx_i;\bx)\Big)^{-1}\Big(\sum_{i=1}^{n}\bq(\bx_i-\bx) L_b(\bx_i;\bx) {K}_h(y_i,y)\Big),
\]
where ${K}_h(y_i,y) = K((y_i-y)/h)/h$ for some (second-order) kernel function $K$. The estimator, $\widehat{f}_{\mathtt{FYT}}(y|\bx)$, is consistent at the boundary of $\mathcal{X}$ (due to the local polynomial regression step on $\bx_i$), but is generally inconsistent at the boundary of $\mathcal{Y}$. Unlike their proposal, our estimator remains consistent at the boundaries of both $\mathcal{X}$ and $\mathcal{Y}$.

\subsection{A local smoothing based estimator}\label{appendix:local smoothing based estimator}

In this appendix we introduce a local smoothing based estimator for the conditional PDF and its derivatives. Recall from Section \ref{section:Introduction} that $\widehat{F}(y|\bx)$ is the estimated conditional CDF formed by a $\qOrder$-th order local polynomial regression. Now let $G$ be some nonnegative measure such that the Radon-Nikodym derivative with respect to the Lebesgue measure is continuous. Then instead of employing a local polynomial regression as in \eqref{eq:local reg along y direction}, we form a conditional PDF (and derivatives) estimator by local smoothing:
\begin{align*}
    \check{f}^{(\vartheta)}(y|\bx) = \be_{1+\vartheta}^\Trans\check{\bbeta}(y|\bx), \qquad
      \check{\bbeta}(y|\bx) = \argmin_{\bv\in\mathbb{R}^{\pOrder+1}} \int_{\mathcal{Y}} \left(\widehat{F}(u|\bx) - \bp(u-y)^\Trans\bv \right)^2 K_h(u;y)\diff G(u),
\end{align*}
which has the closed-form expression: $\check{f}^{(\vartheta)}(y|\bx) = \be_{1+\vartheta}^\Trans\bS_y^{-1} \bar{\bR}_{y,\bx} \widehat{\bS}_\bx^{-1}\be_\mathbf{0}$.
Here we define
\begin{align*}
    \bS_y &= \int_\Ysupp \bp\Big(\frac{u-y}{h}\Big)\frac{1}{h}\bP\Big(\frac{u-y}{h}\Big)^\Trans \text{d}G(u),\\
    \bar{\bR}_{y,\bx} &= \frac{1}{nh^{1+\vartheta}} \sum_{i=1}^{n} \left(\int_\mathcal{Y} \Indicator(y_i\leq u)
    \frac{1}{h}\bP \Big(\frac{u-y}{h}\Big) \diff G(u)\right) \frac{1}{b^d}\bQ\Big(\frac{\bx_i-\bx}{b}\Big)^\Trans.
\end{align*}
Compared to $\widehat{f}^{(\vartheta)}(y|\bx)$, the above local smoothing based estimator requires knowledge of the support $\mathcal{Y}$. On the other hand, $\check{f}^{(\vartheta)}(y|\bx)$ has the advantage that it is immune to low density regions of $y_i$; that is, the new estimator remains valid even when the density of $y_i$ is close to zero. Intuitively, this is because $\check{f}^{(\vartheta)}(y|\bx)$ employs a nonrandom local smoothing in the second step, while  $\widehat{f}^{(\vartheta)}(y|\bx)$ is based on two local polynomial regressions.

Due to space limitations, we investigate the theoretical properties of this estimator in the supplementary material \citep{Cattaneo-Chandak-Jansson-Ma_2023_OnlineSupp}.

\subsection{Proof of Lemma \ref{lem:Equivalent Kernel Representation}}\label{Appendix:proof of stochastic linear rep}

Define
{\small\begin{align*}
u_{i,j} &= \Bigg(\Big(\Indicator(y_i\leq y_j) - F(y_j|\bx_i) \Big)  \bP \Big(\frac{y_j-y}{h}\Big)  -\int_{\mathcal{Y}} \Big[\Indicator(y_i\leq u) - F(u|\bx_i) \Big]  \bP \Big(\frac{u-y}{h}\Big) \diff F_y(u)\Bigg)\bQ  \left(\frac{\bx_i-\bx}{h}\right)^\Trans.
\end{align*}}
We write
{\small\begin{align*}
\tag{I}&\widehat{f}^{(\vartheta)}(y|\bx) =  \frac{1}{nh^{2+d+\vartheta}}\be_{1+\vartheta}^\Trans\widehat{\bS}_y^{-1}\bigg(\sumIN  \int_{\mathcal{Y}} \Big(\Indicator(y_i\leq u) - F(u|\bx_i) \Big)  \bP \Big(\frac{u-y}{h}\Big) \diff F_y(u)\bQ  \left(\frac{\bx_i-\bx}{h}\right)^\Trans\bigg) \widehat{\bS}_\bx^{-1} \be_\mathbf{0}\\
\tag{II}&+ \frac{1}{n^2h^{2+d+\vartheta}} \be_{1+\vartheta}^\Trans\widehat{\bS}_y^{-1} \bigg(\sum_{j=1}^{n}\sum_{i=1}^{n}F(y_j|\bx_i)\bP\Big(\frac{y_j-y}{h}\Big)\bQ\Big(\frac{\bx_i-\bx}{h}\Big)^\Trans  \bigg) \widehat{\bS}_\bx^{-1} \be_\mathbf{0} \\
\tag{III + IV}&+ \frac{1}{n^2h^{2+d+\vartheta}}\be_{1+\vartheta}^\Trans\widehat{\bS}_y^{-1}\bigg(\sum_{i=1}^n u_{i,i}\bigg)\widehat{\bS}_\bx^{-1} \be_\mathbf{0}\quad +\quad  \frac{1}{n^2h^{2+d+\vartheta}}\be_{1+\vartheta}^\Trans\widehat{\bS}_y^{-1}\bigg(\sum_{{i,j=1,i\neq j}}^n u_{i,j}\bigg)\widehat{\bS}_\bx^{-1} \be_\mathbf{0}.
\end{align*}}
We first provide probability concentration results for the matrices $\widehat{\bS}_\bx$ and $\widehat{\bS}_y$. We will then show that term (II) encompasses the target parameter $f^{(\vartheta)}(y|\bx)$ and the smoothing bias. Next, we establish probabilistic orders for (III) and (IV). We analyze term (I) as the last step, which will close the proof.  \medskip

\noindent \textbf{\textit{Convergence of $\widehat{\bS}_\bx$ and $\widehat{\bS}_y$}}. To start, note that $\mathcal{X}$ is compact, then for any $\eta_n>0$, one can find
$\{ \bx_{\ell}:\ 1\leq \ell \leq M_n \}$,
such that
$\mathcal{X}\subseteq \cup_{1\leq \ell\leq M_n}B_{\ell}$,
where $B_{\ell}:=B(\bx_{\ell},\eta_n)$
is the Euclidean ball centered at $\bx_{\ell}$
with radius $\eta_n$.
Define $\mathtt{r} = \sqrt{\log (n)/(nh^d)}$.
Then,
\begin{align*}
\sup_{\bx \in \mathcal{X}} \big|\widehat{\bS}_\bx-{\bS}_\bx\big|
&\leq \max_{1\leq \ell\leq M_n} \big|\widehat{\bS}_{\bx_\ell}-{\bS}_{\bx_\ell}\big|
+ \sup_{1\leq \ell\leq M_n}\sup_{\bx \in B_{\ell}}\big|\widehat{\bS}_\bx - \widehat{\bS}_{\bx_\ell}\big|
+ \sup_{1\leq \ell\leq M_n}\sup_{\bx \in B_{\ell}}\big|{\bS}_\bx - {\bS}_{\bx_\ell}\big|.
\end{align*}
Consider the last term on the RHS. It is straightforward to show that $\bS_\bx$ is continuous with Lipschitz constant of order $h^{-1}$, which implies that $\sup_{1\leq \ell\leq M_n}\sup_{\bx \in B_{\ell}}|{\bS}_\bx - {\bS}_{\bx_\ell}| \precsim {\eta_n}/{h}$.
Similar technique applies to the second term on the RHS: the matrix $\widehat{\bS}_\bx$
is the average of continuous functions with Lipschitz constant of order $h^{-1-d}$, which means $\sup_{1\leq \ell\leq M_n}\sup_{\bx \in B_{\ell}}|\widehat{\bS}_\bx - \widehat{\bS}_{\bx_\ell}|  \precsim {\eta_n}/{h^{1+d}} $.

Now consider the first term. By employing the union bound, we have that, for any constant $\mathfrak{c}_1>0$,
\begin{align*}
&\ \Prob\Big[\max_{1\leq \ell\leq M_n} \big|\widehat{\bS}_{\bx_\ell}-{\bS}_{\bx_\ell}\big|>\mathfrak{c}_1 \mathtt{r}\Big] \leq M_n\max_{1\leq \ell\leq M_n}\Prob\Big[\big|\widehat{\bS}_{\bx_\ell}-{\bS}_{\bx_\ell}\big|>\mathfrak{c}_1 \mathtt{r}\Big].
\end{align*}
To proceed, we recall the formula of $\widehat{\bS}_{\bx}$, and it follows that the summands satisfy
\begin{align*}
\Var\Big[\frac{1}{h^d}\bq\left( \frac{\bx_i-\bx}{h} \right)\bQ  \left( \frac{\bx_i-\bx}{h} \right)^\Trans \Big] \leq  C'h^{-d},\qquad
\Big|\frac{1}{h^d}\bq\left( \frac{\bx_i-\bx}{h} \right)\bQ  \left( \frac{\bx_i-\bx}{h} \right)^\Trans \Big| \leq C' h^{-d},
\end{align*}
where $C'$ is a constant that does not depend on $n$, $h$ or the evaluation point $\bx$.
Applying Bernstein's inequality,
\begin{align*}
M_n\max_{1\leq \ell\leq M_n}\Prob\Big[\big|\widehat{\bS}_{\bx_\ell}-{\bS}_{\bx_\ell}\big|>\mathfrak{c}_1 \mathtt{r}\Big] &\leq  2\exp\Big\{ -\frac{1}{2}\frac{\mathfrak{c}_1^2 \log (n)}{C' + \frac{1}{3}\mathfrak{c}_1C'\mathtt{r}} + \log (M_n)\Big\}.
\end{align*}
To complete the proof, we note that $M_n$ is at most polynomial in $n$ as long
as $\eta_n$ is also polynomial in $n$.
Therefore, one can choose $\eta_n$ sufficiently small so that ${\eta_n}/{h^{1+d}}$ become negligible,
and hence for some constants $\mathfrak{c}_1$, $\mathfrak{c}_2$, and $\mathfrak{c}_3$,
\begin{align*}
\Prob\Big[ \sup_{\bx \in \mathcal{X}} \big|\widehat{\bS}_\bx-{\bS}_\bx\big| > \mathfrak{c}_1\mathtt{r} \Big] \leq \mathfrak{c}_2n^{-\mathfrak{c}_3},
\end{align*}
and $\mathfrak{c}_3$ can be made arbitrarily large with appropriate choices of
$\mathfrak{c}_1$. In other words, we have shown that $\sup_{\bx \in \mathcal{X}} |\widehat{\bS}_\bx-{\bS}_\bx| \precsimConc \sqrt{\log (n)/(nh^d)}$. Analogously, we can show the probability concentration result
$\sup_{y \in \mathcal{Y}} |\widehat{\bS}_y-{\bS}_y| \precsimConc \sqrt{{\log (n)}/{(nh)}}$.  \medskip

\noindent \textbf{\textit{Term (II), and the smoothing bias calculation}}. We start with a Taylor expansion of the conditional CDF up to some order $s$:
\begin{align*}
F(y_j|\bx_i) &= \sum_{\ell + |\bmm| \leq s}  \frac{\partial^{\ell}}{\partial y^{\ell}}\frac{\partial^{\bmm}}{\partial \bx^{\bmm}}F(y|\bx)\frac{1}{\ell!\bmm!} (y_j-y)^{\ell} (\bx_i-\bx)^{\bmm} +  o\bigg( \sum_{ \ell + |\bmm|= s} |y_j-y|^{\ell} |\bx_i-\bx|^{\bmm} \bigg).
\end{align*}
Then,
{\small\begin{align*}
&\frac{1}{n^2h^{2+d+\vartheta}}\sum_{i,j=1}^n\be_{1+\vartheta}^\Trans\widehat{\bS}_y^\inv F(y_j|\bx_i)\bP \left(\frac{y_j-y}{h}\right)\bQ \left(\frac{\bx_i-\bx}{h}\right)^\Trans\widehat{\bS}_{\bx}^\inv \be_{\mathbf{0}}\\
=&\frac{1}{n^2h^{2+d+\vartheta}}\sum_{i,j=1}^n\be_{1+\vartheta}^\Trans\widehat{\bS}_y^\inv \sum_{\ell + |\bmm| \leq s}  \frac{\partial^{\ell}}{\partial y^{\ell}}\frac{\partial^{\bmm}}{\partial \bx^{\bmm}}F(y|\bx)\frac{1}{\ell!\bmm!} (y_j-y)^{\ell} (\bx_i-\bx)^{\bmm} \bP \left(\frac{y_j-y}{h}\right)\bQ \left(\frac{\bx_i-\bx}{h}\right)^\Trans\widehat{\bS}_{\bx}^\inv \be_{\mathbf{0}}\\
&\quad +o\bigg(\frac{1}{n^2h^{2+d+\vartheta}}\be_{1+\vartheta}^\Trans\widehat{\bS}_y^\inv \sum_{i,j=1}^n  \sum_{ \ell + |\bmm|= s} |y_j-y|^{\ell} |\bx_i-\bx|^{\bmm} \left|\bP \left(\frac{y_j-y}{h}\right)\right|\left|\bQ \left(\frac{\bx_i-\bx}{h}\right)\right|\widehat{\bS}_{\bx}^\inv \be_{\mathbf{0}}\bigg)\\
=& f^{(\vartheta)}(y|\bx) + O_\Prob(h^{\qOrder+1} + h^{\pOrder-\vartheta}).
\end{align*}
}
To understand the stochastic order, we notice that the first nonzero term in the summation corresponds to $\ell = 1+\vartheta$ and $\bmm = \mathbf{0}$, which gives rise to the target parameter $f^{(\vartheta)}(y|\bx)$. The next nonzero terms in the summation will be the leading smoothing bias, and correspond to $\ell = 1+\vartheta$ and $|\bmm| = \qOrder+1$, or $\ell =  \pOrder + 1$ and $\bmm = \mathbf{0}$. The leading bias terms will involve random vectors and matrices that are sample averages, whose probabilistic orders can be established using the earlier method of combining discretization, union bound, and Bernstein's inequality. \medskip

\noindent \textbf{\textit{Term (III), the leave-in bias}}. This term arises
because the same observation is used twice: $y_i$ is used to construct the
conditional CDF estimator $\widehat{F}(y|\bx)$, and later as an evaluation point
in the second step local polynomial regression. Term (III) takes the form of a
sample average, and using the earlier method of combining discretization, union bound, and Bernstein's inequality, it is straightforward to show that it has the order
\begin{align*}
\sup_{ y \in \Ysupp, \bx \in \Xsupp} |\text{(III)}| \precsimConc \frac{1}{nh^{1+\vartheta}}\bigg( 1 + \sqrt{\frac{\log (n)}{nh^{1+d}}} \bigg).
\end{align*}

\noindent \textbf{\textit{Term (IV)}}. Term (IV) is a degenerate U-statistic. Take $C$ and $C'$ to be some large constant, and we set
\begin{align*}
A = C',\qquad B^2 = C'nh,\qquad D^2=C'n^2h^{d+1},\qquad t = C(\log (n))\sqrt{n^2h^{d+1}}.
\end{align*}
We apply Lemmas \ref{sa-lem:ustatistic concentration inequality} and \ref{lem:decouping for Ustat}, which give (the value of $C'$ may change for each line)
\begin{align*}
&\Prob\Big[\sup_{y \in \mathcal{Y}, \bx \in \mathcal{X}}\Big| \sum_{i,j=1,i\neq j}^n u_{i,j} \Big| > t \Big] \leq C'\exp\left\{ -\frac{1}{C'}\min\left[ \frac{t}{\sqrt{n^2h^{d+1}}},\ \frac{t^{2/3}}{(nh)^{1/3}},\ t^{1/2} \right] + \log (n) \right\}\\
&= C'\exp\Big\{ -\frac{\sqrt{C}}{C'}\min\Big[ \log (n),\ \left(\log^2 (n)nh^d\right)^{\frac{1}{3}},\ \left(\log^2 (n)n^2h^{1+d}\right)^{\frac{1}{4}} \Big] + \log (n) \Big\}.
\end{align*}
As a result,
\begin{align*}
\sup_{y \in \mathcal{Y}, \bx \in \mathcal{X}} \big|\text{(IV)}\big| \precsimConc \frac{\log (n)}{\sqrt{n^2h^{3+d+2\vartheta}}}.
\end{align*}

\noindent \textbf{\textit{Term (I)}}. To close the proof, we write $\text{(I)} - \bar{f}^{(\vartheta)}(y|\bx) = \text{(I.1)} + \text{(I.2)}$, where
{\small
\begin{align*}
\text{(I.1)} &= \frac{1}{nh^{2+d+\vartheta}}\be_{1+\vartheta}^\Trans(\widehat{\bS}_y^{-1} - {\bS}_y^{-1})\bigg(\sumIN  \int_{\mathcal{Y}} \Big(\Indicator(y_i\leq u) - F(u|\bx_i) \Big)  \bP \Big(\frac{u-y}{h}\Big) \diff F_y(u)\bQ  \left(\frac{\bx_i-\bx}{h}\right)^\Trans\bigg) \widehat{\bS}_\bx^{-1} \be_\mathbf{0},\\
\text{(I.2)}&= \frac{1}{nh^{2+d+\vartheta}}\be_{1+\vartheta}^\Trans{\bS}_y^{-1}\bigg(\sumIN  \int_{\mathcal{Y}} \Big(\Indicator(y_i\leq u) - F(u|\bx_i) \Big)  \bP \Big(\frac{u-y}{h}\Big) \diff F_y(u)\bQ  \left(\frac{\bx_i-\bx}{h}\right)^\Trans\bigg) (\widehat{\bS}_\bx^{-1}-{\bS}_\bx^{-1}) \be_\mathbf{0}.
\end{align*}}

To analyze term (I.1), we have shown that $\sup_{y \in \mathcal{Y}} |\widehat{\bS}_y-{\bS}_y| \precsimConc \sqrt{{\log (n)}/{(nh)}}$ and $\sup_{\bx \in \mathcal{X}}|\widehat{\bS}_\bx| \precsimConc 1 + \sqrt{{\log (n)}/{(nh^d)}}$. Notice that both ${\bS}_y$ and ${\bS}_\bx$ are invertible, which means the same rates apply after inverting the matrices. The middle matrix in (I.1) is a sample average that is mean zero and has variance of order $nh^{2+d}$. We can therefore apply the earlier technique of discretization, union bound, and Bernstein's inequality to show that the middle matrix has the order $\sqrt{\log (n)nh^{2+d}}$. Therefore,
\[
\text{(I.1)} \precsimConc \frac{1}{nh^{2+d+\vartheta}}\sqrt{\frac{\log (n)}{nh}}\sqrt{\log (n)nh^{2+d}} \Big( 1 + \sqrt{\frac{\log (n)}{nh^d}} \Big) \precsim \frac{\log (n)}{\sqrt{n^2h^{3+d+2\vartheta}}}.
\]

To analyze term (I.2), we us the fact that $\sup_{\bx \in \mathcal{X}} |\widehat{\bS}_\bx-{\bS}_\bx| \precsimConc \sqrt{{\log (n)}/{(nh^d)}}$ and the rest of the term is mean zero conditional on $\bx_i$. It remains to compute the variance.
{
\begin{align*}
&\Var\Big[ \be_{1+\vartheta}^\Trans{\bS}_y^{-1}\int_{\mathcal{Y}} \Big(\Indicator(y_i\leq u) - F(u|\bx_i) \Big)  \bP \Big(\frac{u-y}{h}\Big) \diff F_y(u)\bQ  \left(\frac{\bx_i-\bx}{h}\right)^\Trans \Big]\\
&= h^2\Expectation\Big[ \iint_{(\mathcal{Y}-y)/h} \Big(F(y+h(u_1\wedge u_2)|\bx_i) - F(y+hu_1|\bx_i)F(y+hu_2|\bx_i) \Big)    f_y(y+hu_1)f_y(y+hu_2) \\
&\qquad\qquad\qquad\qquad\be_{1+\vartheta}^\Trans{\bS}_y^{-1}\bP (u_1)\bP (u_2)^\Trans {\bS}_y^{-1}\be_{1+\vartheta}\diff u_1\diff u_2\bQ  \left(\frac{\bx_i-\bx}{h}\right)\bQ  \left(\frac{\bx_i-\bx}{h}\right)^\Trans \Big],
\end{align*}
}
where $f_y$ represents the marginal PDF of $y_i$. By a standard Taylor expansion (in $h$) exercise, one can show that the leading term is zero, which means the variance has the order $h^{3+d}$. We can therefore apply the earlier technique (discretization, union bound, and Bernstein's inequality) to show that
\[
\sup_{ y \in \Ysupp, \bx \in \Xsupp} \Big| \be_{1+\vartheta}^\Trans{\bS}_y^{-1}\Big(\sumIN  \int_{\mathcal{Y}} \Big(\Indicator(y_i\leq u) - F(u|\bx_i) \Big)  \bP \Big(\frac{u-y}{h}\Big) \diff F_y(u)\bQ  \left(\frac{\bx_i-\bx}{h}\right)^\Trans\Big) \Big| \precsimConc \sqrt{\log (n)nh^{3+d}}.
\]
As a result,
\[
\text{(I.2)} \precsimConc \frac{1}{nh^{3+d+\vartheta}}\sqrt{\log (n)nh^{3+d}}\sqrt{\frac{\log (n)}{(nh^d)}} = \frac{\log (n)}{\sqrt{n^2h^{1+2d+2\vartheta}}}.
\]

\subsection{Properties of the equivalent kernel}\label{Appendix:properties of K-circle}

In this appendix we prove some useful properties of the equivalent kernel function $\K_{\vartheta,h}^\circ$, which will be employed to establish the strong approximation result in Theorem \ref{thm:Strong Approximation S-hat}.

\begin{lem}[Leading variance]\label{lem:variance of f-bar}
    Suppose Assumptions \ref{as:1-DGP} and \ref{as:2-Kernel} hold. If $h\to 0$ and if $nh^{1+d}/\log (n)\to\infty$, then \eqref{eq:leading variance} holds.
\end{lem}

\begin{proof}[Proof of Lemma \ref{lem:variance of f-bar}]
To save notation, let $\bc_1 = \bS_y^\inv\be_{1+\vartheta}$ and $\bc_2 = \bS_{\bx}^\inv \be_{\mathbf{0}}$. Then
{\begin{align*}
&\Var\left[   \int_{\mathcal{Y}}\big(\Indicator(y_i\leq u) - F(u|\bx_i)\big)\bc_1^\Trans \frac{1}{h}\bP \left(\frac{u-y}{h}\right) f_y(u)\diff u  \frac{1}{h^d}\bQ  \left(\frac{\bx_i-\bx}{h}\right)^\Trans \bc_2 \right]\\
&=  \Expectation\Big[\iint_{\frac{\mathcal{Y}-y}{h}}\Big( F(y+h(u_1\wedge u_2)|\bx_i) - F(y+h u_1|\bx_i)F(y+h u_2|\bx_i) \Big) f_y(y+hu_1)f_y(y+hu_2) \\
\tag{I}&\qquad\qquad\qquad\qquad\qquad\bc_1^\Trans\bP \left(u_1\right)\bc_1^\Trans\bP \left(u_2\right) \diff u_1\diff u_2\Big(\bc_2^\Trans\hd \bQ  \left(\frac{\bx_i-\bx}{h}\right) \Big)^2\Big].
\end{align*}}
We make a further expansion:
\begin{align*}
&\ F(y+h(u_1\wedge u_2)|\bx_i) - F(y+h u_1|\bx_i)F(y+h u_2|\bx_i)\\
&= F(y|\bx_i)(1-F(y|\bx_i))
 + h(u_1\wedge u_2)f(y|\bx_i) - h(u_1+u_2)f(y|\bx_i)F(y|\bx_i)
 + O(h^{2}).
\end{align*}
Note that the remainder term, $O(h^2)$, holds uniformly for $y\in\mathcal{Y}$ and $\bx_i\in\mathcal{X}$ since the conditional distribution function is assumed to have bounded second derivatives. In addition, it is straightforward to verify that with the above Taylor expansion, the first term in (I) is zero, meaning that the leading variance term is
{\begin{align*}
\text{(I)}
&=
 h \left(\be_{1+\vartheta}^\Trans\bS_y^\inv \bT_y\bS_y^\inv \be_{1+\vartheta}\right)
\Expectation\Big[ f(y|\bx_i) \Big(\bc_2^\Trans\hd \bQ  \left(\frac{\bx_i-\bx}{h}\right) \Big)^2\Big]
+ O\Big(\frac{1}{h^{d-2}}\Big).
\end{align*}}
To conclude the proof, we compute the expectation,
\begin{align*}
\text{(I)} &= \frac{1}{h^{d-1}}f(y|\bx) \left(\be_{1+\vartheta}^\Trans\bS_{y}^{-1}\bT_y \bS_{y}^{-1}\be_{1+\vartheta}\right)\left(\be_{\mathbf{0}}^\Trans \bS_{\bx}^{-1} \bT_{\bx} \bS_{\bx}^{-1} \be_{\mathbf{0}}\right) + O\left(\frac{1}{h^{d-2}}\right).
\end{align*}
Therefore, \eqref{eq:leading variance} holds.
\end{proof}

\begin{lem}[Properties of $\K_{\vartheta,h}^\circ$]\label{sa-lem:equivalent kernel properties}
Let Assumptions \ref{as:1-DGP} and \ref{as:2-Kernel} hold. Then \\
(i) $\K_{\vartheta,h}^\circ\left( a,\bb; y,\bx \right)$ is bounded: $\sup_{a,\bb,y,\bx}|\K_{\vartheta,h}^\circ( a,\bb; y,\bx )| \precsim h^{-1-d-\vartheta}$.\\
(ii) $\K_{\vartheta,h}^\circ\left( a,\bb; y,\bx \right)$ is Lipschitz continuous:
\begin{align*}
&\sup_{|a- a'|+|\bb- \bb'|>0,y,\bx}\frac{ \Big|\K_{\vartheta,h}^\circ\left( a,\bb; y,\bx \right) - \K_{\vartheta,h}^\circ\left( a',\bb'; y,\bx \right)\Big|}{|a-a'|+|\bb-\bb'|}  = O \left(h^{-2-d-\vartheta}\right),\\
&\sup_{a,\bb, |y- y'|+ |\bx - \bx'|>0}\frac{ \Big|\K_{\vartheta,h}^\circ\left( a,\bb; y,\bx \right) - \K_{\vartheta,h}^\circ\left( a,\bb; y',\bx' \right)\Big|}{|y-y'| + |\bx-\bx'|}  = O \left(h^{-2-d-\vartheta}\right).
\end{align*}
\end{lem}

\begin{proof}[Proof of Lemma \ref{sa-lem:equivalent kernel properties}]
\textit{Part (i)}. We first rewrite the kernel using change-of-variable. Then, $h^{1+d+\vartheta}\K_{\vartheta,h}^\circ$ takes the form
\begin{align*}
\be_{1+\vartheta}^\Trans\bS_y^\inv \Big[  \int_{\frac{\mathcal{Y} -y}{h}}\Big(\Indicator(a\leq y+hv) - F(y+hv|\bb)\Big)\bP \left(v\right)f_y(y+hv)\diff v\Big] \bQ  \left(\frac{\bb-\bx}{h}\right)^\Trans  {\bS}_{\bx}^\inv \be_{\mathbf{0}}.
\end{align*}
It should be clear that the above is bounded.

\textit{Part (ii)}. From the expression in part (i), it is clear that $h^{1+d+\vartheta}\K_{\vartheta,h}^\circ$ is Lipschitz continuous in $\bb$ with a Lipschitz constant of order $h^{-1}$. Next consider the directions $a$. We have
\begin{align*}
&\sup_{\bb,y,\bx}h^{1+d+\vartheta}|\K_{\vartheta,h}^\circ\left( a,\bb; y,\bx \right) - \K_{\vartheta,h}^\circ\left( a',\bb; y,\bx \right)|\\
\precsim& \sup_{y}\Big|  \int_{\frac{\mathcal{Y} -y}{h}}\Big(\Indicator(a\leq y+hv)-\Indicator(a'\leq y+hv)\Big)\bP \left(v\right)f_y(y+hv)\diff v\Big| \\
\precsim& \sup_{y}\Big|  \int_{\frac{\mathcal{Y} -y}{h}\cap [-1,1] \cap \left[\frac{a-y}{h},\frac{a'-y}{h}\right]}\bP \left(v\right)f_y(y+hv)\diff v\Big|.
\end{align*}
Therefore, the kernel is also Lipschitz-$h^{-1}$ continuous with respect to $a$.

To conclude the proof, it is straightforward to show that $\bS_\bx$ and $\bS_y$ are Lipschitz continuous with respect to $\bx$ and $y$, with the Lipschitz constant of order $1/h$. The same holds for their inverses.
\end{proof}

\begin{lem}[Covering number]\label{lem:covering number for K}
Define
$\Kclass = \{ h^{1+d+\vartheta}\K_{\vartheta,h}^\circ\left( \cdot,\cdot;y,\bx \right):\ y\in\mathcal{Y},\ \bx\in\mathcal{X} \}$. Let Assumptions \ref{as:1-DGP} and \ref{as:2-Kernel} hold. Then
\begin{align*}
\sup_{P}N\Big(\epsilon,\ \Kclass,\ L^1(P)\Big) \leq \mathfrak{c}\frac{1}{\epsilon^{d+2}}+1,
\end{align*}
where the supremum is taken over all probability measures on $[0,1]^{d+1}$, and the constant $\mathfrak{c}$ does not depend on the bandwidth $h$.
\end{lem}

\begin{proof}[Proof of Lemma \ref{lem:covering number for K}]
To show this result, it suffices to consider the uncentered kernel function,
\begin{align*}
h^{1+d+\vartheta}\K_{\vartheta,h}( a,\bb; y,\bx )
          &= \be_{1+\vartheta}^\Trans\bS_y^\inv \int_{\mathcal{Y}} \Indicator(a\leq u)  \frac{1}{h}\bP\Big(\frac{u-y}{h}\Big)\diff F_y(u) \bQ\left(\frac{\bb-\bx}{h}\right)^\Trans \bS_{\bx}^\inv \be_{\mathbf{0}}\\
          &=\be_{1+\vartheta}^\Trans\bS_y^\inv \Big[  \int_{\frac{\mathcal{Y} -y}{h}}\Indicator(a\leq y+hv)\bP \left(v\right)f_y(y+hv)\diff v\Big] \bQ  \left(\frac{\bb-\bx}{h}\right)^\Trans  {\bS}_{\bx}^\inv \be_{\mathbf{0}}.
\end{align*}
We will first show that it has compact support. Consider two cases. If $(a-y)/h> 1$, then the integrand $\Indicator(a\leq y+hv)\bP \left(v\right)$ will be zero because $\bP (v)$ is zero for $v\geq 1$. Therefore, the kernel defined above will be zero as well. For the case that $(a-y)/h\leq -1$, we can simply drop the indicator, as again $\bP (v)$ will be zero for $v\leq -1$. Then the kernel becomes
\begin{align*}
h^{1+d+\vartheta}\K_{\vartheta,h}\left( a,\bb; y,\bx \right) &= \be_{1+\vartheta}^\Trans\bS_y^\inv \Big[  \int_{\frac{\mathcal{Y} -y}{h}}\bP \left(v\right)f_y(y+hv)\diff v\Big] \bQ  \left(\frac{\bb-\bx}{h}\right)^\Trans  {\bS}_{\bx}^\inv \be_{\mathbf{0}},\qquad a\leq -1.
\end{align*}
Note that the matrix, $\bS_y$, can be written as $\bS_y = \int_{\frac{\mathcal{Y} -y}{h}} \bP (v)\bp (v)^\Trans f_y(y+hv)\diff v$, which means its first column is $\int_{\frac{\mathcal{Y} -y}{h}}\bP \left(v\right)f_y(y+hv)\diff v$, showing that the expression above is zero. As for the second argument, $\bb$, we note that $\bQ  ((\bb-\bx)/h)$ is zero if $\bb$ lies outside of an $h$-cube around $\bx$.

With the above result, we can simply apply Lemmas \ref{sa-lem:equivalent kernel properties} and \ref{sa-lem:covering number} to conclude the covering number result for the class $\{ h^{1+d+\vartheta}\K_{\vartheta,h}\left( \cdot,\cdot;y,\bx \right):\ y\in\mathcal{Y},\ \bx\in\mathcal{X} \}$ (note that the boundedness and Lipschitz continuity results in Lemma \ref{sa-lem:equivalent kernel properties} also apply to $\K_{\vartheta,h}$). The same covering number then holds for $\Kclass$, as the two classes differ only by a centering.
\end{proof}

\subsection{Proof of Theorem \ref{thm:Probability Concentration}}\label{Appendix:Proof uniform convergence rate}

Given Lemma \ref{lem:Equivalent Kernel Representation}, we will only need to provide a probability concentration for $\bar{f}^{(\vartheta)}(y|\bx)$. We have established in Lemma \ref{lem:variance of f-bar} that
\begin{align*}
\Var[\K_{\vartheta,h}^\circ( a,\bb; y,\bx )] \leq C'\frac{1}{h^{1+d+2\vartheta}},\quad |\K_{\vartheta,h}^\circ( a,\bb; y,\bx )|\leq C'\frac{1}{h^{1+d+\vartheta}}.
\end{align*}
Then we apply the technique used in the proof of Lemma \ref{lem:Equivalent Kernel Representation} (discretization, union bound, and Bernstein's inequality), which leads to $\sup_{ y \in \Ysupp, \bx \in \Xsupp}|\bar{f}^{(\vartheta)}(y|\bx)| \precsimConc \sqrt{{\log (n)}/{(nh^{1+d+2\vartheta})}}$. To conclude the proof, we notice that the second component in $\mathtt{r_{SL}}$ satisfies
\begin{align*}
\frac{\log (n)}{\sqrt{n^2h^{1+2\vartheta+d + (2 \vee d)}}} &= \sqrt{\frac{\log (n)}{nh^{1+d+2\vartheta}}}\sqrt{\frac{\log (n)}{nh^{ 2 \vee d}}} = o\left(\sqrt{\frac{\log (n)}{nh^{1+d+2\vartheta}}}\right).
\end{align*}

\subsection{Proof of Theorem \ref{thm:Strong Approximation S-hat}}\label{appendix:proof of strong approximation}

It suffices to consider the process $\tilde{\SProc}_{\vartheta}(\eval) = \sumIN h^{1+d+\vartheta}\K_{\vartheta,h}^\circ\left( y_i,\bx_i;y,\bx \right)/\sqrt{n}$, which is the empirical process indexed by the function class $\Kclass$ (defined in Lemma \ref{lem:covering number for K} above). From Lemma \ref{sa-lem:equivalent kernel properties}, the functions in the above class are uniformly bounded. Lemma \ref{lem:covering number for K} shows that the function class above is of VC type, and the covering number does not depend on the bandwidth. The measurability condition required in Lemma~\ref{sa-lem:Rio coupling} also holds, as our function class is indexed by $(\eval) \in[0,1]^{d+1}$, and the functions in $\Kclass$ are continuous in $y$ and $\bx$.

Now the only missing ingredient is the total variation of the functions in $\Kclass$. First, note that the function $h^{1+d+\vartheta}\K_{\vartheta,h}^\circ\left( \cdot,\cdot;y,\bx \right)$ is Lipschitz continuous with respect to the arguments, and the Lipschitz constant is of order $h^{-1}$. Therefore, its total variation is bounded by
\begin{align*}
\mathrm{TV}_{(\eval)} = \mathrm{TV}\Big(  h^{1+d+\vartheta}\K_{\vartheta,h}^\circ\left( \cdot,\cdot;y,\bx \right)\Big) \precsim \frac{1}{h}\mathrm{vol}\Big(\mathrm{supp}\Big(\K_{\vartheta,h}( \cdot,\cdot;y,\bx )\Big)\Big),
\end{align*}
where $\mathrm{vol}\left(\mathrm{supp}\left(\cdot\right)\right)$ denotes the Euclidean volume of the support, and $\K_{\vartheta,h}$ is defined in the proof of Lemma \ref{lem:covering number for K}. We also showed in the proof of Lemma \ref{lem:covering number for K} that $\K_{\vartheta,h}$ has compact support, leading to $\mathrm{TV}_{\Kclass} = \sup_{y \in \mathcal{Y}, \bx \in \mathcal{X}}\mathrm{TV}_{(\eval)} \precsim h^d$.

Putting all pieces together, we conclude that there exists a centered Gaussian process, $\tilde{\mathbb{G}}_{\vartheta} $ which has the same covariance kernel as $\tilde{\SProc}_{\vartheta} $, such that
\begin{align*}
\Prob\Big[\sup_{y \in \mathcal{Y}, \bx \in \mathcal{X}} \left| \tilde{\SProc}^{\prime}_{\vartheta}(\eval)- \tilde{\mathbb{G}}_{\vartheta}(\eval) \right| \geq \mathfrak{c}_1 \Big(\sqrt{\frac{h^d\log n}{n^{\frac{1}{d+1}}}} + \sqrt{\frac{\log^3 n}{n}}\Big)\Big] \leq \mathfrak{c}_2n^{-\mathfrak{c}_3},
\end{align*}
where $\tilde{\SProc}^{\prime}_{\vartheta}(\eval)$ is a copy of $\tilde{\SProc}_{\vartheta}(\eval)$.

\subsection{Proof of Theorem \ref{thm:Kolmogorov-Smirnov Distance: Suprema}}\label{appendix:KS distance}

First consider $\widehat{\TProc}_{\vartheta}(\eval)$. The difference between $\widehat{\TProc}_{\vartheta}(\eval)$ and $\widehat{\SProc}_{\vartheta}(\eval)$ is
\begin{align*}
\widehat{\TProc}_{\vartheta}(\eval) - \widehat{\SProc}_{\vartheta}(\eval) &= \bigg( \sqrt{\frac{\mathsf{V}_{\vartheta}(\eval)}{\widehat{\mathsf{V}}_{\vartheta}(\eval)}} - 1 \bigg)\widehat{\SProc}_{\vartheta}(\eval).
\end{align*}
With Theorem \ref{thm:Probability Concentration}, Lemma \ref{lem:Covariance Estimation} and the variance bound in \eqref{eq:leading variance} (also see Lemma \ref{lem:variance of f-bar} in Appendix \ref{Appendix:properties of K-circle}), we have
\begin{align*}
\sup_{y \in \mathcal{Y}, \bx \in \mathcal{X}}\left| \widehat{\TProc}_{\vartheta}(\eval) - \widehat{\SProc}_{\vartheta}(\eval) \right| \precsimConc \mathtt{r_{VE}}\Big(h^{\pOrder-\vartheta} + \sqrt{\frac{\log(n)}{nh^{1+d+2\vartheta}}}\Big)\sqrt{nh^{1+d+2\vartheta}}
\precsim \sqrt{\log (n)}\mathtt{r_{VE}}.
\end{align*}

Next, we establish a Gaussian comparison result. Consider an $\epsilon$ discretization of $\mathcal{Y}\times \mathcal{X}$, which is denoted by $\mathcal{A}_{\epsilon} = \{(y_\ell,\bx_{\ell}^\Trans):\ 1\leq \ell \leq L\}$. Then one can define two Gaussian vectors, $\bz,\widehat{\bz}\in\mathbb{R}^L$, such that
\begin{align*}
\Cov[z_{\ell}, z_{\ell'}] = \rho_\vartheta(y_{\ell},\bx_{\ell},y_{\ell'},\bx_{\ell'}),\quad \Cov[\widehat{z}_{\ell}, \widehat{z}_{\ell'}|\bX,\bY] = \widehat{\rho}_\vartheta(y_{\ell},\bx_{\ell},y_{\ell'},\bx_{\ell'}).
\end{align*}
Then we apply the Gaussian comparison result in Lemma \ref{sa-lem:gaussian comparison} and the correlation estimation error rate in Lemma \ref{lem:Covariance Estimation}, which lead to
\begin{align*}
\sup_{u \in \mathbb{R}}\Big| \Prob\Big[  \sup_{1 \leq \ell \leq L}|\widehat{\GProc}_{\vartheta}(y_{\ell},\bx_{\ell})|\leq u \Big| \bY,\bX \Big] - \Prob\Big[ \sup_{1 \leq \ell \leq L}|\GProc_{\vartheta}(y_{\ell},\bx_{\ell})|\leq u  \Big] \Big| \precsim_\Prob \sqrt{\mathtt{r_{VE}}}\log\Big(\frac{1}{\epsilon}\Big).
\end{align*}
Since $\epsilon$ only enters the above error bound logarithmically, one can choose $\epsilon = n^{-\mathfrak{c}}$ for some $\mathfrak{c}$ large enough, so that the error that arises from discretization becomes negligible. In other words, we have
\begin{align*}
\sup_{u \in \mathbb{R}}\Big| \Prob\Big[  \sup_{ y \in \Ysupp, \bx \in \Xsupp}|\widehat{\GProc}_{\vartheta}(y,\bx)|\leq u \Big| \bY,\bX \Big] - \Prob\Big[ \sup_{ y \in \Ysupp, \bx \in \Xsupp}|\GProc_{\vartheta}(y,\bx)|\leq u  \Big] \Big| \precsim_\Prob   \log(n)\sqrt{\mathtt{r_{VE}}}.
\end{align*}

Now consider $\widehat{\TProc}_{\vartheta}(\eval)$ again. Given the bound on the difference, $\sup_{y \in \mathcal{Y}, \bx \in \mathcal{X}}| \widehat{\TProc}_{\vartheta}(\eval) - \widehat{\SProc}_{\vartheta}(\eval) |$, and the strong approximation in Theorem \ref{thm:Strong Approximation S-hat}, we clearly have
\begin{align*}
\Prob\bigg[\sup_{y \in \mathcal{Y}, \bx \in \mathcal{X}}| \GProc_{\vartheta}(\eval)| \leq u &- \mathfrak{c}_1(\sqrt{\log (n)}\mathtt{r_{VE}}+\mathtt{r_{SA}})\bigg] - \mathfrak{c}_2n^{-\mathfrak{c}_3}
\leq \Prob\bigg[ \sup_{y \in \mathcal{Y}, \bx \in \mathcal{X}}|\widehat{\TProc}_{\vartheta}(\eval) | \leq u \bigg] \\
&\leq \Prob\bigg[\sup_{y \in \mathcal{Y}, \bx \in \mathcal{X}}| \GProc_{\vartheta}(\eval)| \leq u + \mathfrak{c}_1(\sqrt{\log (n)}\mathtt{r_{VE}}+\mathtt{r_{SA}})\bigg] + \mathfrak{c}_2n^{-\mathfrak{c}_3}.
\end{align*}
Finally, we apply the Gaussian comparison result, which implies that
\begin{align*}
&\ \sup_{u \in \mathbb{R}}\Big| \Prob\Big[ \sup_{y \in \mathcal{Y}, \bx \in \mathcal{X}}|\widehat{\TProc}_{\vartheta}(\eval)|\leq u  \Big] -
\Prob\Big[  \sup_{y \in \mathcal{Y}, \bx \in \mathcal{X}}|\widehat{\GProc}_{\vartheta}(\eval)|\leq u \Big| \bY,\bX \Big]  \Big|\\
\precsim_\Prob&\  \mathfrak{c}_2n^{-\mathfrak{c}_3} + \log(n)\sqrt{\mathtt{r_{VE}}} + \sup_{u\in\mathbb{R}}\Prob\Big[ \sup_{y \in \mathcal{Y}, \bx \in \mathcal{X}}| \GProc_{\vartheta}(\eval)| \in [u, u+\mathfrak{c}_1(\sqrt{\log (n)}\mathtt{r_{VE}}+\mathtt{r_{SA}})] \Big].
\end{align*}
Finally, due to Lemma \ref{sa-lem:anti-concentration}, we have
\begin{align*}
\sup_{u\in\mathbb{R}}\Prob\Big[ \sup_{y \in \mathcal{Y}, \bx \in \mathcal{X}}| \GProc_{\vartheta}(\eval)| \in [u, u+\mathfrak{c}_1(\sqrt{\log (n)}\mathtt{r_{VE}}+\mathtt{r_{SA}})] \Big] \precsim \sqrt{\log (n)}(\sqrt{\log (n)}\mathtt{r_{VE}}+\mathtt{r_{SA}}).
\end{align*}

\subsection{Derivation of \texorpdfstring{\eqref{eq:the constrained estimators}} . }\label{appendix: derivation of constrained estimators}

First consider $\widehat{f}_{\mathtt{N}}(y|\bx)$. If the unconstrained estimator, $\widehat{f}(y|\bx)$, is already nonnegative, then the constraint in the least squares problem is not binding, which means in this case $\widehat{f}_{\mathtt{N}}(y|\bx) = \widehat{f}(y|\bx)$. Now assume $\widehat{f}(y|\bx) < 0$. Since the least squares objective function is strictly convex, the solution will be on the boundary of the set $\{\bu\in\mathbb{R}^{\pOrder+1}:\be_1^\Trans \bu \geq 0\}$, leading to $\widehat{f}_{\mathtt{N}}(y|\bx) = 0$. Therefore, we have the expression $\widehat{f}_{\mathtt{N}}(y|\bx) = \max\{0, \widehat{f}(y|\bx)\}$ in \eqref{eq:the constrained estimators}.

The expression of $\widehat{f}_{\mathtt{I}}(y|\bx)$ in \eqref{eq:the constrained estimators} follows from Jensen's inequality, which is binding if and only if $g(y)/\widehat{f}_{\mathtt{N}}(y|\bx)$ is constant (in $y$).

\subsection{Proof of Lemma \ref{lem:Equivalent Kernel Representation, NN and I estimators}, Theorems \ref{thm:Probability Concentration, NN and I estimators} and \ref{thm:Strong Approximation S-hat, NN and I estimators}}

We write $\widehat{f}_{\mathtt{N}}(y|\bx) = \widehat{f}(y|\bx) - \Indicator(\widehat{f}(y|\bx) < 0)\cdot\widehat{f}(y|\bx)$.
We first study the indicator function. Take $\mathtt{r}$ to be any sequence shrinking to 0, and $\mathfrak{c}_1$ some positive constant. Then
\begin{align*}
&\Prob\Big[\sup_{ y \in \mathcal{Y}, \bx \in \mathcal{X}}\Indicator\left(\widehat{f}(y|\bx) < 0\right) > \mathtt{r}\mathfrak{c}_1  \Big]\leq \Prob\Big[ \sup_{y\in\mathcal{Y},\bx\in\mathcal{X}}\Big|\widehat{f}(y|\bx) - f(y|\bx)\Big| > \inf_{y\in\mathcal{Y},\bx\in\mathcal{X}}f(y|\bx) \Big].
\end{align*}
Then by the probability concentration in Theorem \ref{thm:Probability Concentration}, it should be obvious that the the above probability vanishes faster than any polynomials of $n$ (recall that we assume the conditional density is uniformly bounded away from zero); that is,
$\sup_{ y \in \mathcal{Y}, \bx \in \mathcal{X}}\Indicator(\widehat{f}(y|\bx) < 0)\precsimConc \mathtt{r}$
for any vanishing sequence $\mathtt{r}$. This shows that $\sup_{ y \in \mathcal{Y}, \bx \in \mathcal{X}}|\widehat{f}_{\mathtt{N}}(y|\bx) - \widehat{f}(y|\bx)|\precsimConc \mathtt{r}$. By letting $\mathtt{r}$ shrinking to 0 fast enough, we have the stochastic linearization for $\widehat{f}_{\mathtt{N}}(y|\bx)$.

Next, For $\widehat{f}_{\mathtt{I}}(y|\bx)$, we employ the following decomposition:
\begin{align*}
\widehat{f}_{\mathtt{I}}(y|\bx)
= \widehat{f}_{\mathtt{N}}(y|\bx) - \frac{\widehat{f}_{\mathtt{N}}(y|\bx)}{\int_{\mathcal{Y}} \widehat{f}_{\mathtt{N}}(u|\bx) \diff u}\int_{\mathcal{Y}} \big(\widehat{f}_{\mathtt{N}}(u|\bx) - {f}(u|\bx)\big)\diff u.
\end{align*}
Then we can write
\begin{align*}
&\sup_{ y \in \mathcal{Y}, \bx \in \mathcal{X}}\Big|\widehat{f}_{\mathtt{I}}(y|\bx) -  f(y|\bx) - \Big(\bar{f}^{(0)}(y|\bx) - f(y|\bx)\int_{\mathcal{Y}} \bar{f}^{(0)}(u|\bx)\diff u\Big)\Big|\\
&\leq  \sup_{ y \in \mathcal{Y}, \bx \in \mathcal{X}}\Big|\widehat{f}_{\mathtt{N}}(y|\bx) -  f(y|\bx) - \bar{f}^{(0)}(y|\bx) \Big|  + \sup_{ y \in \mathcal{Y}, \bx \in \mathcal{X}}\Big|  \frac{\widehat{f}_{\mathtt{N}}(y|\bx)}{\int_{\mathcal{Y}} \widehat{f}_{\mathtt{N}}(u|\bx) \diff u} - f(y|\bx)\Big|\cdot\Big|\int_{\mathcal{Y}} \bar{f}^{(0)}(u|\bx)\diff u\Big|\\
&\qquad\qquad+ \sup_{ y \in \mathcal{Y}, \bx \in \mathcal{X}}\Big|  \frac{\widehat{f}_{\mathtt{N}}(y|\bx)}{\int_{\mathcal{Y}} \widehat{f}_{\mathtt{N}}(u|\bx) \diff u}\Big|\cdot \Big|\int_{\mathcal{Y}} \big(\widehat{f}_{\mathtt{N}}(u|\bx) - {f}(u|\bx)-\bar{f}^{(0)}(u|\bx)\big)\diff u \Big|\\
&\precsimConc  \mathtt{r_{SL}} + \Big(h^{\pOrder} + \sqrt{\frac{\log(n)}{nh^{1+d}}}\Big)\Big( \sqrt{\frac{\log(n)}{nh^{d}}}\Big) \precsim \mathtt{r_{SL}}.
\end{align*}
In the above, we have used the result that $\sup_{ \bx \in \Xsupp}\Var[ \int_{\mathcal{Y}} \bar{f}^{(0)}(u|\bx)\diff u ] \precsim (nh^d)^{-1}$, which shows that $\int_{\mathcal{Y}} \bar{f}^{(0)}(u|\bx)\diff u$ has a smaller asymptotic order compared to $\bar{f}^{(0)}(y|\bx)$.

To prove Theorem \ref{thm:Probability Concentration, NN and I estimators}, we combine the results in Lemma \ref{lem:Equivalent Kernel Representation, NN and I estimators} and Theorem \ref{thm:Probability Concentration}. The strong approximation for $\widehat{\SProc}_{\mathtt{N}}$ in Theorem \ref{thm:Strong Approximation S-hat, NN and I estimators} follows from Lemma \ref{lem:Equivalent Kernel Representation, NN and I estimators} and Theorem \ref{thm:Strong Approximation S-hat}. The strong approximation for $\widehat{\SProc}_{\mathtt{I}}$ also follows from Lemma \ref{lem:Equivalent Kernel Representation, NN and I estimators}, as the stochastic linearization of $\widehat{f}_{\mathtt{I}}$ is a linear functional of $\bar{f}^{(0)}$.

\subsection{A result on covering number}\label{Appendix: covering number}

In this appendix, we prove a general result on the uniform covering number for function classes consisting of kernels. Importantly, we allow the kernels in the function class to take different shapes and to depend on a range of bandwidths.

\begin{lem}[Covering number]\label{sa-lem:covering number}
Let $h>0$, and $\mathfrak{c}>0$ be a (large) generic constant which does not depend on $h$. Define the class of functions
$$\mathcal{G}=\left\{g_{\bz}\left(\frac{\cdot-\bz}{ah}\right):\ \bz\in [0,1]^{d},\ 1\leq a\leq \mathfrak{c}\right\}.$$
Assume (i) boundedness: $\sup_{\bz,\bz'}|g_{\bz}(\bz')| \leq \mathfrak{c}$. (ii) $g_{\bz}(\cdot)$ is supported in $ [-1, 1]^{d}$ for all $\bz$. (iii) Lipschitz continuity: $\sup_{\bz}|g_{\bz}(\bz')-g_{\bz}(\bz'')| \leq \mathfrak{c}|\bz'-\bz''|$ and $\sup_{\bz}|g_{\bz'}(\bz)-g_{\bz''}(\bz)| \leq \mathfrak{c}h^{-1} |\bz'-\bz''|$.
Then, for any probability measure $P$, the $L^1(P)$-covering number of the class $\mathcal{G}$ satisfies
\begin{align*}
N\big((2\mathfrak{c}+1)^{d+1}\epsilon,\ \mathcal{G},\ L^1(P)\big)\leq \mathfrak{c}'\frac{1}{\epsilon^{d+2}}+1,
\end{align*}
where $\mathfrak{c}'$ is some constant that depends only on $\mathfrak{c}$ and $d$.
\end{lem}

This rate, $\epsilon^{-d-2}$, is clearly suboptimal for very small $\epsilon$. The reason is that when we fix $h$ and consider how the covering number changes as $\epsilon\downarrow 0$, the optimal rate is $\epsilon^{-d-1}$, as in this case the class of functions is fixed (c.f. Theorem 2.7.11 in \cite{van-der-Vaart-Wellner_1996_Book}). Such suboptimality is introduced because we prefer a covering number that depends only on $\epsilon$ (but not $h$). The result we derived performs better for moderate and large $\epsilon$ (relative to $h$).

Now consider how the above (a sharper result for moderate and large $\epsilon$) manifests itself in our proof below. Take a fixed $\epsilon$. As the bandwidth shrinks to 0, we will be employing finer partitions of $[0,1]^d$. However, not all of the sets in the partition matter for bounding the covering number,
because there are at most $\epsilon^\inv $ sets carrying a probability mass larger than $\epsilon$.
Given that the functions we consider have compact support, most of them become irrelevant in our calculation of the covering number. Indeed, a function only makes a nontrivial contribution if its support intersects with some set in the (very fine) partition whose $P$-measure exceeds $\epsilon$.
Therefore, instead of considering all $h^{-d}$ partitions, we only need to focus on $\epsilon^\inv $ of them, which is why an extra $\epsilon^\inv $ term is introduced.

Finally, from the definition of $\mathcal{G}$, it is clear that the covering number obtained above allows for a range of bandwidths (captured by $ah$ with $1\leq a\leq \mathfrak{c}$). If we instead consider the restricted function class, $\{g_{\bz}((\cdot-\bz)/h):\ \bz\in [0,1]^{d}\}$, then a sharper bound will apply: $\mathfrak{c}'{\epsilon^{-d-1}}+1$.

\begin{proof}[Proof of Lemma \ref{sa-lem:covering number}]
This proof strategy is motivated by Lemma 4.1 in \cite{Rio_1994_PTRF}. Take $\ell=\lfloor 1/h\rfloor$, and partition each coordinate $[0,1]$ into $\ell$ intervals of equal length. This will lead to a partition $\mathcal{A}=\{A_j:\ 1\leq j\leq \ell^d\}$ of $[0,1]^d$. Next, consider sets whose $P$-measure exceeds $\epsilon$,
\[ \mathcal{A}_{P,\epsilon} = \{A\in \mathcal{A}:\ P[A]> \epsilon\}, \]
and their $\mathfrak{c}h$-enlargements
\[ \mathcal{A}_{P,\epsilon}^{\mathfrak{c}h} = \{A + [-\mathfrak{c}h,\mathfrak{c}h]^d:\ A\in \mathcal{A}_{P,\epsilon}\}. \]
\textbf{Case 1:} $\bz$ does not belong to any set in $\mathcal{A}_{P,\epsilon}^{\mathfrak{c}h}$. This implies that the support of the function $g_{\bz}\left( \frac{\cdot - \bz}{ah} \right)$ will not intersect with any set in $\mathcal{A}_{P,\epsilon}$. We also notice that
\begin{align*}
\int \Big|g_{\bz}\Big( \frac{\cdot - \bz}{ah} \Big)\Big|\diff P \leq \mathfrak{c} P\big[ ah\cdot\mathrm{supp}(g_{\bz}(\cdot)) + \bz \big] \leq \mathfrak{c} P\big[ \mathfrak{c}h\cdot\mathrm{supp}(g_{\bz}(\cdot)) + \bz \big].
\end{align*}
Define the complement of $\mathcal{A}_{P,\epsilon}$ as $\mathcal{A}_{P,\epsilon}^{\perp} = \{A\in \mathcal{A}:\ P[A]\leq \epsilon\}$, then the set $\mathfrak{c}h\cdot\mathrm{supp}(g_{\bz}(\cdot)) + \bz$ will be completely covered by sets in $\mathcal{A}_{P,\epsilon}^{\perp}$. To determine the maximum number of intersections between $\mathfrak{c}h\cdot\mathrm{supp}(g_{\bz}(\cdot)) + \bz$ and sets in $\mathcal{A}_{P,\epsilon}^{\perp}$, it suffices to consider the Euclidean volume of the enlarged set $\mathfrak{c}h\cdot\mathrm{supp}(g_{\bz}(\cdot)) + \bz + [-\ell^\inv ,\ell^\inv ]^d$, which is $(2\mathfrak{c}h + \ell^\inv )^d$. The Euclidean volume of each set in $\mathcal{A}_{P,\epsilon}^{\perp}$ is $\ell^{-d}$. Therefore, the set $\mathfrak{c}h\cdot\mathrm{supp}(g_{\bz}(\cdot)) + \bz$ can intersect with at most
\[ \frac{(2\mathfrak{c}h + \ell^\inv )^d}{\ell^{-d}} = \left( 2\mathfrak{c}h\ell + 1 \right)^d \leq \left( 2\mathfrak{c} + 1 \right)^d\]
sets in $\mathcal{A}_{P,\epsilon}^{\perp}$. As a result, we conclude that $\int |g_{\bz}\left( \frac{\cdot - \bz}{ah} \right)|\diff P \leq \mathfrak{c} \left( 2\mathfrak{c} + 1 \right)^d\epsilon$.
This leads to our first result. Let $A_{P,\epsilon}^{\mathfrak{c}h} = \cup \mathcal{A}_{P,\epsilon}^{\mathfrak{c}h}$ be the union of sets in $\mathcal{A}_{P,\epsilon}^{\mathfrak{c}h}$, then

\[N\Big((2\mathfrak{c}+1)^{d+1}\epsilon,\ \mathcal{G}_1,\ L^1(P)\Big) = 1,\quad \text{where }\mathcal{G}_1 = \Big\{ g_{\bz}\left( \frac{\cdot - \bz}{ah} \right):\ \bz\not\in A_{P,\epsilon}^{\mathfrak{c}h},\ 1\leq a\leq \mathfrak{c}\Big\}.\]
As remark, we note that the function class $\mathcal{G}_1$ changes with respect to $h$, $\epsilon$, as well as the probability measure $P$.\smallskip

\noindent \textbf{Case 2:} $\bz$ belongs to some set in $\mathcal{A}_{P,\epsilon}^{\mathfrak{c}h}$. Each set in $\mathcal{A}_{P,\epsilon}^{\mathfrak{c}h}$ is a cube with edge length $\ell^\inv  + 2\mathfrak{c}h \leq 2(\mathfrak{c}+1)h$, because $h\ell \geq 0.5$. Then the covering number of $A_{P,\epsilon}^{\mathfrak{c}h}$ (under the Euclidean distance) is
\[
N\Big(h\epsilon,\ A_{P,\epsilon}^{\mathfrak{c}h},\ |\cdot|\Big) \leq \sum_{A\in \mathcal{A}_{P,\epsilon}^{\mathfrak{c}h}}N(h\epsilon,\ A,\ |\cdot|)  \leq \mathrm{card}(\mathcal{A}_{P,\epsilon}^{\mathfrak{c}h})\cdot \mathfrak{c}'\frac{1}{\epsilon^d} \leq \mathfrak{c}'\frac{1}{\epsilon^{d+1}}.
\]
Here, $\mathfrak{c}'$ is some fixed number that only depends on $\mathfrak{c}$ and $d$. Using the Lipschitz property, we have
\begin{align*}
\int \Big| g_{\bz}\left( \frac{\cdot - \bz}{ah} \right) - g_{\bz'}\left( \frac{\cdot - \bz'}{a'h} \right) \Big|\diff P &\leq
 2\mathfrak{c}h^\inv |\bz-\bz'| + \mathfrak{c}^2 |a-a'|.
\end{align*}
Now define $\mathcal{G}_2 = \mathcal{G}\backslash \mathcal{G}_1  = \{ g_{\bz}((\cdot - \bz)/(ah)):\ \bz\in A_{P,\epsilon}^{\mathfrak{c}h},\ 1\leq a\leq \mathfrak{c}\}$, then
\begin{align*}
&N\Big((2\mathfrak{c}+1)^{d+1}\epsilon,\ \mathcal{G}_2,\ L^1(P)\Big) \leq  N\Big(\frac{(2\mathfrak{c}+1)^{d+1}}{4\mathfrak{c}}h\epsilon,\ A_{P,\epsilon}^{\mathfrak{c}h},\ |\cdot|\Big)  N\Big(\frac{(2\mathfrak{c}+1)^{d+1}}{2\mathfrak{c}^2}\epsilon,\ [1,\mathfrak{c}],\ |\cdot|\Big)\leq \frac{\mathfrak{c}'}{\epsilon^{d+2}}.
\end{align*}
This closes the proof.
\end{proof}

\subsection{Technical lemmas}

\begin{lem}[Equation (3.5) in \cite{Gine-Latala-Zinn_2000_Ustat}]\label{sa-lem:ustatistic concentration inequality}
Let $\{{z}_i, 1\leq i\leq n\}$ be independent random variables, and $\{\tilde{z}_i, 1\leq i\leq n\}$ be an independent copy of $\{{z}_i, 1\leq i\leq n\}$. For a degenerate and decoupled second order U-statistic, $\sum_{i,j=1,i\neq j}^n u_{ij}(z_i,\tilde{z}_j)$, the following holds:
\begin{align*}
\Prob\Big[ \Big|\sum_{i,j,i\neq j}^n u_{ij}(z_i,\tilde{z}_j)\Big| > t \Big] \leq \mathfrak{c}\exp\Big\{ -\frac{1}{\mathfrak{c}}\min\Big[ \frac{t}{A},\ \left(\frac{t}{B}\right)^{\frac{2}{3}},\ \left( \frac{t}{C} \right)^{\frac{1}{2}} \Big] \Big\},
\end{align*}
where $\mathfrak{c}$ is some absolute constant, and $A$, $B$ and $C$ are any constants satisfying
\begin{align*}
A^2 &\geq \sum_{i,j=1,i\neq j}^n \Expectation [u_{ij}(z_i,\tilde{z}_j)^2], \quad B^2 \geq \max_{1\leq i,j\leq n}\Big[  \sup_{w}\Big|\sum_{i=1}^n \Expectation [u_{ij}(z_i,w)^2]\Big| ,\  \sup_{v}\Big|\sum_{j=1}^n \Expectation [u_{ij}(v,\tilde{z}_j)^2]\Big|  \Big],\\
C   &\geq \max_{1\leq i,j\leq n}\sup_{v,w}| u_{ij}(v,w) |.
\end{align*}
\end{lem}

To apply the above lemma, an additional decoupling step is usually needed. Fortunately, the decoupling step only introduces an extra constant, but will not affect the order of the tail probability bound. Formally,
\begin{lem}[\cite{de1995decoupling}]\label{lem:decouping for Ustat}
Consider the setting of Lemma~\ref{sa-lem:ustatistic concentration inequality}. Then
\begin{align*}
\Prob\Big[ \Big|\sum_{i,j,i\neq j}^n u_{ij}(z_i,{z}_j)\Big| > t \Big] \leq \mathfrak{c}\cdot \Prob\Big[ \mathfrak{c}\Big|\sum_{i,j,i\neq j}^n u_{ij}(z_i,\tilde{z}_j)\Big| > t \Big],
\end{align*}
where $\mathfrak{c}$ is an absolute constant.
\end{lem}
As a result, we will apply Lemma~\ref{sa-lem:ustatistic concentration inequality} without explicitly mentioning the decoupling step or the extra constant it introduces.

\begin{lem}[Theorem 1.1 in \cite{Rio_1994_PTRF}]\label{sa-lem:Rio coupling}
Let $\bz_1,\bz_2,\dots,\bz_n$ be iid random vectors with continuous and strictly positive density on $[0,1]^d$, and $d\geq 2$. Let $\mathcal{G}$ be a class of functions from $[0,1]^d$ to $[-1,1]$, satisfying $\sup_{P}N(\epsilon,\ \mathcal{G},\ L^1(P)) \leq \mathfrak{c}_1\epsilon^{-\mathfrak{c}_2}$, where the supremum is taken over all probability measures on $[0,1]^d$, and $\mathfrak{c}_1$ and $\mathfrak{c}_2$ are constants that can depend on $\mathcal{G}$.  In addition, assume the following measurability condition holds: there exists a Suslin space $\mathcal{S}$ and a mapping $\FProc:\mathcal{S}\to \mathcal{G}$, such that $(s,\bz)\mapsto \FProc(s,\bz)$ is measurable. Let
\begin{align*}
\mathrm{TV}_{\mathcal{G}} &= \sup_{g\in\mathcal{G}}\sup_{\phi\in \mathcal{C}^{\infty}_1([0,1]^d)}\int_{[0,1]^d} g(\bz)\mathrm{div}\phi(\bz)\diff \bz,
\end{align*}
where $\mathrm{div}$ is the divergence operator, and $\mathcal{C}^{\infty}_1([0,1]^d)$ is the collection of infinitely differentiable functions with values in $\mathbb{R}^d$, support included in $[0,1]^d$, and supremum norm bounded by 1. Then on a possibly enlarged probability space, there exists a centered Gaussian process, $\mathbb{G}$, indexed by $\mathcal{G}$, such that (i) $\Cov[\mathbb{G}(g),\mathbb{G}(g')] = \Cov[g(\bz_i),g'(\bz_i)]$, and (ii) for any $t\geq \mathfrak{c}_3\log n$,
\begin{align*}
\Prob\Big[ \sqrt{n}\sup_{g\in\mathcal{G}}\left| \BProc(g) - \mathbb{G}(g) \right| \geq \mathfrak{c}_3\sqrt{n^{\frac{d-1}{d}}\ t\ \mathrm{TV}_{\mathcal{G}}} + \mathfrak{c}_3t\sqrt{\log (n)} \Big] \leq e^{-t}.
\end{align*}
In the above, $\BProc = \sumIN (g(\bz_i) - \Expectation[g(\bz_i)])/\sqrt{n}$ is the empirical process indexed by $\mathcal{G}$,
and $\mathfrak{c}_3$ is some constant that only depends on $d$, $\mathfrak{c}_1$, and $\mathfrak{c}_2$.
\end{lem}

\begin{lem}[Corollary 5.1 in \citet{chernozhukov2022central}]\label{sa-lem:gaussian comparison}
Let $\bz_1,\bz_2\in\mathbb{R}^{\ell_n}$ be two mean-zero Gaussian random vectors with covariance matrices $\bOmega_{1}$ and $\bOmega_{2}$, respectively. Further assume that the diagonal elements in $\bOmega_{1}$ are all one. Then
\begin{align*}
\sup_{\substack{ A\ \text{rectangular}}} \left| \Prob\left[ \bz_1 \in A \right]
-
\Prob\left[ \bz_2 \in A \right] \right| \leq \mathfrak{c}\sqrt{\Vert \bOmega_{1} - \bOmega_{2} \Vert_{\infty}}\log (\ell_n),
\end{align*}
where $\Vert \cdot\Vert_\infty$ denotes the supremum norm, and $\mathfrak{c}$ is an absolute constant.
\end{lem}

\begin{lem}[Theorem 2.1 in \cite{Chernozhukov-Chetverikov-Kato_2014b_AoS}]\label{sa-lem:anti-concentration}
Let $\GProc$ be a centered and separable Gaussian process indexed by $g\in\mathcal{G}$ such that $\Var[\GProc(g)]=1$ for all $g\in\mathcal{G}$. Assume $\sup_{g\in\mathcal{G}}\GProc(g)<\infty$ almost surely. Define $C_\mathcal{G} = \Expectation[\sup_{g\in \mathcal{G}}\GProc(g)]$. Then for all $\epsilon>0$,
\begin{align*}
\sup_{u\in\mathbb{R}}\Prob\Big[\Big| \sup_{g\in \mathcal{G}}\GProc(g) - u \Big|\leq \epsilon\Big] \leq 4\epsilon(C_\mathcal{G}+1).
\end{align*}
\end{lem}

\end{appendix}

\begin{acks}[Acknowledgments]
The authors thank
the editor,
two anonymous reviewers,
Jianqing Fan,
Jason Klusowski,
Will Underwood, 
Jingshen Wang, 
and Rae Yu for their thoughtful discussions and valuable feedback.
\end{acks}

\begin{funding}
Cattaneo gratefully acknowledges financial support from the
National Science Foundation through grants SES-1947805 and DMS-2210561, and
from the National Institute of Health (R01 GM072611-16).

Jansson gratefully acknowledges financial support from the
National Science Foundation through grant SES-1947662 and
the research support of CREATES.
\end{funding}

\begin{supplement}
\stitle{Supplementary material to ``Boundary adaptive local polynomial conditional density estimators''}
\sdescription{The supplementary material \citep{Cattaneo-Chandak-Jansson-Ma_2023_OnlineSupp} contains general theoretical results encompassing those discussed in the main paper, includes proofs of those general results, and discusses additional methodological and technical results.}
\end{supplement}


\bibliographystyle{imsart-nameyear} 
\bibliography{CCJM_2023_Bernoulli--bib}       

\end{document}


\begin{frontmatter}
\title{Supplementary material to ``Boundary adaptive local polynomial conditional density estimators''}
\runtitle{Supplementary material}

\begin{aug}
\author[A]{\inits{MC}\fnms{Matias D.}~\snm{Cattaneo}\ead[label=e1]{cattaneo@princeton.edu}}
\author[A]{\inits{RC}\fnms{Rajita}~\snm{Chandak}\ead[label=e2]{rchandak@princeton.edu}}
\author[B]{\inits{MJ}\fnms{Michael}~\snm{Jansson}\ead[label=e3]{mjansson@econ.berkeley.edu}}
\author[D]{\inits{XM}\fnms{Xinwei}~\snm{Ma}\ead[label=e4]{x1ma@ucsd.edu}\orcid{0000-0001-8827-9146}}
\address[A]{Department of Operations Research and Financial Engineering, Princeton University, Princeton NJ, United States\printead[presep={,\ }]{e1,e2}}

\address[B]{Department of Economics, UC Berkeley, Berkeley CA, United States\printead[presep={,\ }]{e3}}


\address[D]{Department of Economics, UC San Diego, La Jolla CA, United States\printead[presep={,\ }]{e4}}
\end{aug}

\begin{abstract}
This Supplementary Material contains general theoretical results encompassing those discussed in the main paper, includes proofs of those general results, and discusses additional methodological and technical results. A companion \texttt{R} package is available at
\url{https://nppackages.github.io/lpcde/}.
\end{abstract}

\begin{keyword}
\kwd{Conditional density estimation}
\kwd{confidence bands}
\kwd{local polynomial methods}
\kwd{specification testing}
\kwd{strong approximation}
\kwd{uniform inference}
\end{keyword}

\end{frontmatter}

\tableofcontents

\section{Setup}\label{section:setup}
Let $\bx_i\in\mathbb{R}^d$ and $y_i\in\mathbb{R}$ be continuously distributed random variables supported on $\mathcal{X}=[0,1]^d$ and $\mathcal{Y}=[0,1]$, respectively.
We are interested in estimating the conditional distribution function and its derivatives:
\begin{align*}
\theta_\order(y|\bx)  &= \frac{\partial^\mu }{\partial y^\mu}\frac{\partial^{\bnu} }{\partial \bx^{\bnu}} F(y|\bx),
\end{align*}
where $\mu\in\mathbb{N}$, and $\bnu\in\mathbb{N}^d$ representing multi-indices. (In the main paper we only consider the estimation of conditional density and derivatives thereof, that is, we set $\bnu=0$ and $\mu = \vartheta + 1\geq 1$.)

To present our estimation strategy, we start from $\theta_{0,\bnu}$, the conditional distribution function and its derivatives with respect to the conditioning variable, and apply the local polynomial method:
\begin{align*}
&\widehat{\frac{\partial^{\bnu} }{\partial \bx^{\bnu}} F}(y|\bx)  = \be^{\Trans}_{\bnu} \hat{\bgamma}(y|\bx), &
  \hat{\bgamma}(y|\bx) = \argmin_{\bv\in\mathbb{R}^{\qOrder_d+1}} \sum_{i=1}^n \Big(\Indicator(y_i\leq y) - \bq(\bx_i-\bx)^\Trans\bv \Big)^2 L_h(\bx_i;\bx),
 \end{align*}
where and $\be_{\bnu}^\Trans$ is a basis vector extracting the corresponding estimate. We can write the solution in closed form as
\begin{align*}
\widehat{\frac{\partial^{\bnu} }{\partial \bx^{\bnu}} F}(y|\bx) = \be_{\bnu}^\Trans\hat{\bS}_{\bx}^\inv \Big( \frac{1}{n}\sumIN\Indicator(y_i\leq y)\bq\Big(\frac{\bx_i-\bx}{h}\Big) L_h(\bx_i;\bx) \Big),
\end{align*}
where
\begin{align*}
\hat{\bS}_\bx &= \frac{1}{n}\sum_{i=1}^{n}\bq\Big(\frac{\bx_i-\bx}{h}\Big)\bq\Big(\frac{\bx_i-\bx}{h}\Big)^\Trans L_h(\bx_i;\bx).
\end{align*}

To estimate $\theta_\order$, we further smooth via local polynomials along the $y$-direction:
\begin{align*}
&\hat{\theta}_{\mu,\bnu} (y|\bx) = \be_{\mu}^\Trans\hat{\bbeta}(y|\bx), &
\hat{\bbeta}(y|\bx) = \argmin_{\bu\in\mathbb{R}^{\pOrder+1}} \sumIN \left(\hat{F}(y_i|\bx) - \bp (y_i-y)^\Trans\bu \right)^2 K_h(y_i;y).
\end{align*}
We can write the solution in closed-form as
\begin{align*}
\hat{\theta}_\order (y|\bx) = \be_{\mu}^\Trans\hat{\bS}_y^{-1} \hat{\bR}_{y,\bx} \hat{\bS}_\bx^{-1}\be_{\bnu},
\end{align*}
where
\begin{align*}
\hat{\bS}_y       &= \frac{1}{n}\sum_{i=1}^{n}\bp \Big(\frac{y_i-y}{h}\Big)\bp \Big(\frac{y_i-y}{h}\Big)^\Trans K_h(y_i;y), \qquad \text{and} \\
\hat{\bR}_{y,\bx} &= \frac{1}{n^2h^{\mu+|\bnu|}}\sum_{j=1}^{n}\sum_{i=1}^{n}\Indicator(y_i\leq y_j)\bp \Big(\frac{y_j-y}{h}\Big)K_h(y_j;y)\bq\Big(\frac{\bx_i-\bx}{h}\Big)^\Trans L_h(\bx_i;\bx).
\end{align*}
While in the above we considered local polynomial regressions along both the $\bx$- and $y$-directions, it is also possible to employ a local smoothing technique. To be precise, let $G$ be some function such that the following Lebesgue-Stieltjes integration is well-defined, then an alternative estimator can be constructed as
\begin{align*}
&\check{\theta}_{\order} (y|\bx)= \be_{\mu}^\Trans \check{\bbeta}(y|\bx), &
  \check{\bbeta}(y|\bx) = \argmin_{\bv\in\mathbb{R}^{\pOrder+1}} \int \left(\hat{F}(u|\bx) - \bp (u-y)^\Trans\bv \right)^2 K_h(u;y) \diff G(u),
\end{align*}
which has the solution
\begin{align*}
\check{\theta}_\order (y|\bx)
&=\be_{\mu}^{\Trans}\bS_y^{-1} \bar{\bR}_{y,\bx} \hat{\bS}_\bx^{-1}\be_\mathbf{\bnu},
\end{align*}
where
\begin{align*}
\bS_y       &= \int_\mathcal{Y} \bp \Big(\frac{u-y}{h}\Big)\bp \Big(\frac{u-y}{h}\Big)^\Trans K_h(u;y) \diff G(u), \quad \text{and}\\
\bar{\bR}_{y,\bx} &= \frac{1}{nh^{\mu+|\bnu|}} \sum_{i=1}^{n} \left(\int_\mathcal{Y} \Indicator(y_i\leq u)
\bp \Big(\frac{u-y}{h}\Big)K_h(u;y) \diff G(u)\right) \bq\Big(\frac{\bx_i-\bx}{h}\Big)^\Trans L_h(\bx_i;\bx).
\end{align*}

\subsection{Notation}
Limits are taken with respect to the sample size tending to infinity and the bandwidth shrinking to zero (i.e., $n\to\infty$ and $h\to 0$). For two positive sequences, $a_n\precsim b_n$ implies that $\limsup_{n\to\infty}|a_n/b_n|<\infty$. Similarly, $a_n\precsim_\Prob b_n$ means $|a_n/b_n|$ is asymptotically bounded in probability. We also adopt the small-o and big-O notation: $a_n = \Op(b_n)$ is just $a_n\precsim_\Prob b_n$, and $a_n = \op(b_n)$ means $a_n/b_n$ converges to zero in probability. Constants that do not depend on the sample size or the bandwidth will be denoted by $\mathfrak{c}$, $\mathfrak{c}_1$, $\mathfrak{c}_2$, etc.

We introduce another notation, $\Oconc$, which not only provides an asymptotic order but also controls the tail probability. To be specific, $a_n = \Oconc(b_n)$ if for any $\mathfrak{c}_1 > 0$, there exists some $\mathfrak{c}_2$ such that
\[\limsup_{n\to\infty} n^{\mathfrak{c}_1}\Prob\left[ a_n \geq \mathfrak{c}_2b_n \right] < \infty.\]
Here the subscript, $\mathtt{TC}$, stands for ``tail control.'' Finally, let $\bX = (\bx_1^\Trans,\cdots,\bx_n^\Trans)^\Trans$ and $\bY = (y_1,\cdots,y_n)^\Trans$ be the data matrices.

\begin{itemize}[leftmargin=*]
\item $F(y|\bx)$ and $f(y|\bx)$: the conditional distribution and density functions of $y_i$ (at $y$) given $\bx_i=\bx$. The marginal distributions and densities are denoted by $F_y$, $F_{\bx}$, $f_y$, and $f_{\bx}$, respectively.

\item $y$ and $\bx$: the evaluation points.

\item $\mathcal{X}=[0,1]^d$ and $\mathcal{Y}=[0,1]$, the support of $\bx_i$ and $y_i$, respectively.

\item $h$: the bandwidth sequence.

\item $K$: the kernel function, and $L$ is the product kernel: $L(\bx) = K(x_1)K(x_2)\cdots K(x_d)$.

\item $\bp ,\ \bq $: polynomial expansions.

\item $\bP$ and $\bQ$: defined as $\bp(\cdot)K(\cdot)$ and $\bq(\cdot)L(\cdot)$, respectively.

\item $\be_\mu$ and $\be_{\bnu}$: standard basis vectors extracting the $\mu$-th and $\bnu$-th element in the expansion of
$\bp $ and $\bq $ for univariate and multivariate arguments, respectively.

\item $G(\cdot)$ the weighting function used in $\check{\theta}_{\order} $, with its Lebesgue density denoted by $g(\cdot)$.

\item  Some matrices
{\begin{alignat*}{2}
\bS_y &= \int_{\frac{\mathcal{Y}-y}{h}}\bp \left(u\right)\bP \left(u\right)^\Trans g(y+hu)\diff u,\quad
&&\hat{\bS}_y = \frac{1}{nh}\sum_{i=1}^{n}\bp \Big(\frac{y_i-y}{h}\Big)\bP \Big(\frac{y_i-y}{h}\Big)^\Trans, \\
\bc_{y,\ell} &= \int_{\frac{\mathcal{Y}-y}{h}}\frac{u^{\ell}}{\ell!} \bP \left(u\right) g(y+hu)\diff u,
&&\hat{\bc}_{y,\ell} = \frac{1}{nh}\sumIN\frac{1}{\ell!}\left(\frac{y_{i}-y}{h}\right)^{\ell}\bP \left(\frac{y_{i}-y}{h}\right)^\Trans,\\
\bS_{\bx} &= \int_{\frac{\mathcal{X}-\bx}{h}} \bq \left(\bv\right)\bQ  \left(\bv\right)^\Trans f_{\bx}(\bx+h\bv)  \diff \bv,\qquad
&&\hat{\bS}_\bx = \frac{1}{nh^d}\sum_{i=1}^{n}\bq\Big(\frac{\bx_i-\bx}{h}\Big)\bQ  \Big(\frac{\bx_i-\bx}{h}\Big)^\Trans, \\
\bc_{\bx,\bmm} &= \int_{\frac{\mathcal{X}-\bx}{h}} \frac{\bv^{\bmm}}{\bmm!}\bQ  \left(\bv\right) f_{\bx}(\bx+h\bv)  \diff \bv
,\qquad &&
\hat{\bc}_{\bx,\bmm} = \frac{1}{nh^d}\sumIN\frac{1}{\bmm!}\left(\frac{\bx_i-\bx}{h}\right)^{\bmm}\bQ  \left(\frac{\bx_i-\bx}{h}\right),\\
\bT_{\bx} &= \int_{\frac{\mathcal{X}-\bx}{h}}   \bQ  \left(\bv\right)\bQ  \left(\bv\right)^\Trans f_{\bx}(\bx+h\bv)\diff \bv,\qquad
&&\hat{\bT}_{\bx} = \frac{1}{nh^d}\sumIN \bQ  \left(\frac{\bx_i-\bx}{h}\right)\bQ  \left(\frac{\bx_i-\bx}{h}\right)^\Trans,
\end{alignat*}
\vskip-0.5cm
\begin{align*}
  \bT_{y} &= \iint_{\frac{\mathcal{Y}-y}{h}}   \min(u_1, u_2)\bP
            \left(u_1\right)\bP \left(u_2\right)^\Trans g(y+hu_1)g(y+hu_2)\diff
            u_1\diff u_2,
  \\
  \widehat{\bT}_{y} &= \frac{1}{n^2h^3}\sum_{i,j=1}^n   \big(\min(y_i,y_j)-y\big)\bP\Big(\frac{y_i-y}{h}\Big)\bP\Big(\frac{y_j-y}{h}\Big)^\Trans,
\end{align*}
\vskip-0.5cm
\begin{align*}
\hat{\bR}_{y,\bx} &= \frac{1}{n^2h^{1+d+\mu+|\bnu|}}\sum_{j=1}^{n}\sum_{i=1}^{n}\Indicator(y_i\leq y_j)\bP \Big(\frac{y_j-y}{h}\Big)\bQ  \Big(\frac{\bx_i-\bx}{h}\Big)^\Trans, \\
\bar{\bR}_{y,\bx} &= \frac{1}{nh^{1+d+\mu+|\bnu|}} \sum_{i=1}^{n} \left(\int_{\mathcal{Y}} \Indicator(y_i\leq u)
\bP \Big(\frac{u-y}{h}\Big) \diff G(u)\right) \bQ  \Big(\frac{\bx_i-\bx}{h}\Big)^\Trans.
\end{align*}}

\item Equivalent kernels:
{\small\begin{align*}
\check{\K}_{\mu,\bnu,h}^\circ\left( a,\bb; y,\bx \right) &= \frac{1}{h^{\mu+|\bnu|}}\be_{\mu}^\Trans\bS_y^\inv \left[  \int_{\mathcal{Y}}\Big(\Indicator(a\leq u) - \hat{F}(u|\bb)\Big)\frac{1}{h}\bP \left(\frac{u-y}{h}\right) \diff G(u)\right] \frac{1}{h^d}\bQ  \left(\frac{\bb-\bx}{h}\right)^\Trans \hat{\bS}_{\bx}^\inv \be_{\bnu},\\
\hat{\K}_{\mu,\bnu,h}^{\circ}\left( a,\bb; y,\bx \right) &= \frac{1}{h^{\mu+|\bnu|}}\be_{\mu}^\Trans \hat{\bS}_y^\inv \left[\frac{1}{n}\sum_{j=1}^n \Big(\Indicator(a\leq y_j) - \hat{F}(y_j|\bb) \Big) \frac{1}{h}\bP \Big(\frac{y_j-y}{h}\Big)\right] \frac{1}{h^d}\bQ  \left(\frac{\bb-\bx}{h}\right)^\Trans \hat{\bS}_{\bx}^\inv  \be_{\bnu},\\
\K_{\mu,\bnu,h}^\circ\left( a,\bb; y,\bx \right) & = \frac{1}{h^{\mu+|\bnu|}}\be_{\mu}^\Trans\bS_y^\inv \left[  \int_{\mathcal{Y}}\Big(\Indicator(a\leq u) - F(u|\bb)\Big)\frac{1}{h}\bP \left(\frac{u-y}{h}\right)  \diff G(u)\right] \frac{1}{h^d}\bQ  \left(\frac{\bb-\bx}{h}\right)^\Trans   {\bS}_{\bx}^\inv \be_{\bnu},\\
\K_{\mu,\bnu,h}\left( a,\bb; y,\bx \right) & = \frac{1}{h^{\mu+|\bnu|}}\be_{\mu}^\Trans\bS_y^\inv \left[  \int_{\mathcal{Y}}\Indicator(a\leq u)\frac{1}{h}\bP \left(\frac{u-y}{h}\right)  \diff G(u)\right] \frac{1}{h^d}\bQ  \left(\frac{\bb-\bx}{h}\right)^\Trans   {\bS}_{\bx}^\inv \be_{\bnu}.
\end{align*}}

\item Some rates:
\begin{align*}
\mathtt{r_B} &= h^{\qOrder+1-|\bnu|} + h^{\pOrder+1-\mu},\quad \mathtt{r_V} = \sqrt{\frac{1}{nh^{d + 2|\bnu|+2\mu-1}}},\\
\mathtt{r_{BE}} &=
\begin{cases}
\frac{1}{\sqrt{nh^d}} & \text{if $\mu=0$, and $\theta_{0, \mathbf{0}}\neq 0$ or $1$}\\
\frac{1}{\sqrt{nh^{d+1}}}  &\text{if $\mu > 0$, \ \ or $\theta_{0, \mathbf{0}}= 0$ or $1$}
\end{cases},\\
\mathtt{r_{VE}} &= h^{ \qOrder+\frac{1}{2}} + \sqrt{\frac{\log (n)}{nh^{d+1}}},\quad \mathtt{r_{SE}}=\sqrt{\log (n) }\mathtt{r_{VE}},\quad \mathtt{r_{SA}} = \left(\frac{\log^{d+1} (n)}{nh^{d+1}}\right)^{\frac{1}{2d+2}}.
\end{align*}
\end{itemize}

\subsection{Overview}\label{sa-subsec:Overview and Decomposition}

In this subsection we provide an overview of the main results. Underlying
assumptions and precise statements of the lemmas and theorems will be given in
later sections. First consider $\check{\theta}_{\mu,\bnu}(y|\bx)$, with a conditional expectation decomposition:
\begin{align*}
&\check{\theta}_{\order}(y|\bx)
=h^{-\mu -|\bnu|}\be_{\mu}^{\Trans}\bS_y^{-1} \left[\frac{1}{n} \sum_{i=1}^{n} \left(\int_{\mathcal{Y}} F(u|\bx_i)
\frac{1}{h}\bP \Big(\frac{u-y}{h}\Big) \diff G(u)\right) \frac{1}{h^d}\bQ  \Big(\frac{\bx_i-\bx}{h}\Big)^\Trans \right] \hat{\bS}_\bx^{-1}\be_\mathbf{\bnu}\\
&\quad+h^{-\mu -|\bnu|}\be_{\mu}^{\Trans}\bS_y^{-1} \left[\frac{1}{n} \sum_{i=1}^{n} \left(\int_{\mathcal{Y}} \Big(\Indicator(y_i \leq u) - F(u|\bx_i)\Big)
\frac{1}{h}\bP \Big(\frac{u-y}{h}\Big) \diff G(u)\right) \frac{1}{h^d}\bQ  \Big(\frac{\bx_i-\bx}{h}\Big)^\Trans \right] \hat{\bS}_\bx^{-1}\be_\mathbf{\bnu}.
\end{align*}
As we will show in Section~\ref{sa-sec:Pointwise Large-sample Properties}, the first term above consists of the centering of the estimator (i.e., the parameter of interest $\theta_{\mu,\bnu}(y|\bx)$) and the smoothing bias. The second term, on the other hand, gives the asymptotic representation of the estimator. To be precise, we have
\begin{align*}
\check{\theta}_{\mu,\bnu}(y|\bx) - \theta_{\mu,\bnu}(y|\bx)&=  \frac{1}{n}\sumIN \K_{\mu,\bnu,h}^\circ\left( y_i,\bx_i;y,\bx \right) \\
&\qquad+ \Op\left( h^{\qOrder+1-|\bnu|} +  h^{\pOrder+1-\mu} + \sqrt{\mathsf{V}_{\order}(\eval)}\frac{\log (n)}{\sqrt{nh^d}} \right).
\end{align*}
As a result, we can focus on establishing properties of the the first term, which provides an equivalent kernel expression. Denote its variance by $\mathsf{V}_{\order}(\eval)$. Then we show that the standardized process,
\begin{align*}
\bar{\SProc}_{\order}(\eval) &= \frac{1}{n\sqrt{\mathsf{V}_{\order}(\eval)}}\sumIN \K_{\mu,\bnu,h}^\circ\left( y_i,\bx_i;y,\bx \right),
\end{align*}
is approximately normally distributed both pointwise and uniformly for $y\in\mathcal{Y}$ and $\bx\in\mathcal{X}$. To be even more precise, we establish a strong approximation result, meaning that there exists a copy $\bar{\SProc}^{\prime}_{\order}(\eval)$, and a Gaussian process $\GProc_{\order}(\eval)$ with the same covariance structure, such that
\[
\sup_{ y \in \mathcal{Y}, \bx \in \mathcal{X}}\left|\bar{\SProc}^{\prime}_{\order}(\eval) - \GProc_{\order}(\eval)\right| = \Op\left(\frac{\log^{d+1} (n)}{nh^{d+1}}\right)^{\frac{1}{2d+2}}.
\]
Together with a feasible variance-covariance estimator, the strong approximation result not only allows us to construct confidence bands for the target parameter and test shape restrictions, but also provides an explicit characterization of the coverage error probability for those procedures.

Inside the remainder term, $h^{\qOrder+1-|\bnu|} +  h^{\pOrder+1-\mu}$ is the order of the leading smoothing bias, and $\log(n) \sqrt{\mathsf{V}_{\order}(\eval)/(nh^d)}$ arises from the linearization step which replaces the random matrix $\hat{\bS}_{\bx}$ by its large-sample analogue $\bS_{\bx}$. It is worth mentioning that the order of the remainder term is uniformly valid for $y\in\mathcal{Y}$ and $\bx\in\mathcal{X}$, which is why an extra logarithmic factor is present.

Now consider the other estimator, $\hat{\theta}_{\mu,\bnu}(y|\bx)$. While it is not possible to take a conditional expectation, we can still ``center'' the estimator with the conditional distribution function. That is,
\begin{align*}
&\hat{\theta}_{\order}(y|\bx)
=h^{-\mu -|\bnu|}\be_{\mu}^{\Trans}\hat{\bS}_y^{-1} \bigg[\frac{1}{n^2}\sum_{i=1}^{n}\sum_{j=1}^{n}F(y_j|\bx_i) \frac{1}{h}\bP \Big(\frac{y_j-y}{h}\Big)\frac{1}{h^d}\bQ  \Big(\frac{\bx_i-\bx}{h}\Big)^\Trans \bigg] \hat{\bS}_\bx^{-1}\be_\mathbf{\bnu}\\
&\quad+h^{-\mu -|\bnu|}\be_{\mu}^{\Trans}\hat{\bS}_y^{-1} \bigg[\frac{1}{n^2}\sum_{i=1}^{n}\sum_{j=1}^{n}\Big(\Indicator(y_i\leq y_j)-F(y_j|\bx_i)\Big)\frac{1}{h}\bP \Big(\frac{y_j-y}{h}\Big)\frac{1}{h^d}\bQ  \Big(\frac{\bx_i-\bx}{h}\Big)^\Trans \bigg] \hat{\bS}_\bx^{-1}\be_\mathbf{\bnu}.
\end{align*}
As before, the first term captures the target parameter and the smoothing bias. The analysis of the second term is more involved. Besides the asymptotic linear representation term, it also consists of a leave-in bias term (since the same observation is used twice) and a second order U-statistic. We show that the following expansion holds uniformly for $y\in\mathcal{Y}$ and $\bx\in\mathcal{X}$:
\begin{align*}
\hat{\theta}_{\order}(y|\bx) - \theta_{\mu,\bnu}(y|\bx)&=  \frac{1}{n}\sumIN \K_{\mu,\bnu,h}^\circ\left( y_i,\bx_i;y,\bx \right) \\
&\qquad+ \Op\left( h^{\qOrder+1-|\bnu|} +  h^{\pOrder+1-\mu} + \sqrt{\mathsf{V}_{\order}(\eval)}\frac{\log (n)}{\sqrt{nh^d}} + \frac{\log (n)}{\sqrt{n^2h^{d+2\mu+2|\bnu|+1}}} \right).
\end{align*}
Here, the contribution of the U-statistic is represented by the order $\log  (n)/\sqrt{n^2h^{d+2\mu+2|\bnu|+1}}$ in the remainder term. Interestingly, this term is negligible compared to the standard error, $\sqrt{\mathsf{V}_{\order}(\eval)}$, provided that $\log (n)/(nh^2)\to \infty$.

The above demonstrates that important large-sample properties of the local regression based estimator, $\hat{\theta}_{\order}(y|\bx)$ --- such as pointwise and uniform normal approximation --- stem from the equivalent kernel representation. Here we note that the representation holds by setting $G = F_y$. In other words, $\hat{\theta}_{\order}(y|\bx)$ is first-order asymptotically equivalent to $\check{\theta}_{\order}(y|\bx)$ with the (infeasible) local smoothing using the marginal distribution $F_y$.

\subsection{Assumptions}\label{sa-sec:assumptions}

We make the following assumptions on the joint distribution, the kernel function, and the weighting $ G$.

\begin{assumDGP}[Data generating process]\label{sa-assu:1 dgp} (i) $\{y_i, \bx_i\}_{1\leq i\leq n}$ is a random sample from the absolutely continuous joint distribution $F$ supported on $\mathcal{Y} \times \mathcal{X} = [0,1]^{1+d}$.
(ii) The joint density, $f$, is continuous and is bounded away from zero.
(iii) $\theta_{2,\mathbf{0}}$ exists and is continuous.
\end{assumDGP}

\begin{assumK}[Kernel]\label{sa-assu:2 kernel}\ \\
The kernel function $K$ is nonnegative, symmetric, supported on $[-1,1]$, Lipschitz continuous, and integrates to one.
\end{assumK}

\begin{assumW}[Weighting function]\label{sa-assu:3 density}\ \\
The weighting function $ G$ is continuously differentiable with a Lebesgue density denoted by $g$.
\end{assumW}

\section{Pointwise large-sample properties}\label{sa-sec:Pointwise Large-sample Properties}

We first present several uniform convergence results which will be used later to establish pointwise and uniform properties of our estimators.

\begin{lem}[Matrix convergence]\label{sa-lem:Taking hat off}
Let Assumptions \ref{sa-assu:1 dgp}, \ref{sa-assu:2 kernel}, and \ref{sa-assu:3 density} hold with $h\to 0$, $nh^d/\log (n)\to \infty$, and $G = F_y$. Then
\begin{alignat*}{2}
& \sup_{y \in \mathcal{Y}} \left|\hat{\bS}_y-\bS_y\right| =\Oconc \left(\sqrt{\frac{\log (n)}{nh}}\right),
\quad
&&\sup_{y \in \mathcal{Y}} \left|\hat{\bc}_{y,\ell}-\bc_{y,\ell}\right| = \Oconc\left(\sqrt{\frac{\log (n)}{nh}}\right),\\
&\sup_{\bx\in \mathcal{X}} \left|\hat{\bS}_\bx-\bS_\bx\right| = \Oconc \left(\sqrt{\frac{\log (n)}{nh^d}}\right),
\quad
&&\sup_{\bx\in \mathcal{X}} \left|\hat{\bc}_{\bx,\bmm}-\bc_{\bx,\bmm}\right|  = \Oconc\left(\sqrt{\frac{\log (n)}{nh^d}}\right) ,\\
&\sup_{\bx\in \mathcal{X}} \left|\hat{\bT}_\bx-\bT_\bx\right| = \Oconc \left(\sqrt{\frac{\log (n)}{nh^d}}\right). \qquad &&
\end{alignat*}
If in addition that $nh^{d+1}/\log (n) \to \infty$, then
\begin{alignat*}{2}
& \sup_{y \in \mathcal{Y}, \bx \in \mathcal{X}} \Big|\be_{\mu}^{\Trans}\bS_y^{-1}\left(\bar{\bR}_{\eval} - \Expectation \left[\bar{\bR}_{\eval}| \bX \right]\right)\Big|  = \Oconc \left(\mathtt{r}_1\right) ,\quad \text{where }\mathtt{r}_1 = \begin{cases}
\sqrt{\frac{\log (n)}{nh^{d+2\mu+2|\bnu|}}}     & \text{if $\mu=0$}\\
\sqrt{\frac{\log (n)}{nh^{d+2\mu+2|\bnu|-1}}} & \text{if $\mu>0$}
\end{cases}.
\end{alignat*}
\end{lem}

We now follow the decomposition in Section \ref{sa-subsec:Overview and Decomposition} and study the leading bias of our estimators.

\begin{lem}[Bias]\label{sa-lem:Bias}
Let Assumptions \ref{sa-assu:1 dgp}, \ref{sa-assu:2 kernel} and \ref{sa-assu:3 density} hold with $h \to 0$ and $nh^{d} / \log (n)\to\infty$. In addition, $\theta_{\mu',\bnu'}$ exists and is continuous for all $\mu' + |\bnu'|= \max\{ \qOrder+1+\mu,\ \pOrder+1+|\bnu| \}$. Then
\begin{align*}
&\be_{\mu}^{\Trans}\bS_y^{-1} \Big[ \frac{1}{nh^{\mu+|\bnu|}} \sum_{i=1}^{n}\Big(\int_\mathcal{Y} F(u|\bx_i)
\frac{1}{h}\bP \Big(\frac{u-y}{h}\Big) \diff G(u)\Big) \frac{1}{h^d}\bQ  \Big(\frac{\bx_i-\bx}{h}\Big)^\Trans \Big] \hat{\bS}_\bx^{-1}\be_\mathbf{\bnu}\\
&= \theta_{\order}(y|\bx) + \mathsf{B}_{\order}(y,\bx) + \op\left( h^{\qOrder+1-|\bnu|} +  h^{\pOrder+1-\mu} \right),
\end{align*}
where
\begin{align*}
\mathsf{B}_{\order} (\eval) &=
h^{\qOrder+1-|\bnu|}\underbrace{\sum_{|\bmm|=\qOrder+1}\theta_{\mu,\bmm}(y|\bx) \bc_{\bx,\bmm}^\Trans \bS_{\bx}^{-1}\be_{\bnu} }_{\textstyle B_{(i),\qOrder+1}(y,\bx) }
+ h^{\pOrder+1-\mu}\underbrace{\theta_{\pOrder+1,\bnu}(y|\bx)\bc_{y,\pOrder+1}^\Trans\bS_{y}^{-1}\be_{\mu}}_{\textstyle B_{(ii),\pOrder+1}(y,\bx)}.
\end{align*}
Similarly,
\begin{align*}
&\be_{\mu}^{\Trans}\hat{\bS}_y^{-1} \Big[\frac{1}{n^2h^{\mu +|\bnu|}}\sum_{i=1}^{n}\sum_{j=1}^{n}F(y_j|\bx_i) \frac{1}{h}\bP \Big(\frac{y_j-y}{h}\Big)\frac{1}{h^d}\bQ  \Big(\frac{\bx_i-\bx}{h}\Big)^\Trans \Big] \hat{\bS}_\bx^{-1}\be_\mathbf{\bnu} \\
=&\ \theta_{\mu,\bnu}(y|\bx) + \mathsf{B}_{\order} (\eval)
 + \op\left( h^{\qOrder+1-|\bnu|} +  h^{\pOrder+1-\mu} \right).
\end{align*}
\end{lem}

For future reference, we define the order of the leading bias as
\begin{align*}
\mathtt{r_B} &= h^{\qOrder+1-|\bnu|} + h^{\pOrder+1-\mu}.
\end{align*}

\begin{remark}[Higher-order bias]
Because the leading bias established in the lemma can be exactly zero, one may need to extract higher-order terms for bandwidth selection:
\begin{align*}
\mathsf{B}_{\order} (\eval) &=
h^{\qOrder+1-|\bnu|}B_{(i),\qOrder+1}(y,\bx)
+ h^{\pOrder+1-\mu}B_{(ii),\pOrder+1}(y,\bx)\\
&+ h^{\qOrder+2-|\bnu|}B_{(i),\qOrder+2}(y,\bx)+ h^{\pOrder+2-\mu}B_{(ii),\pOrder+2}(y,\bx)
+ h^{\pOrder+\qOrder+2-\mu-|\bnu|}B_{(iii),\pOrder+1, \qOrder+1}(y,\bx),
\end{align*}
where
\begin{align*}
B_{(i),\qOrder+2}(y,\bx) &= \sum_{|\bmm|=\qOrder+2}\theta_{\mu,\bmm}(y|\bx) \bc_{\bx,\bmm}^\Trans \bS_{\bx}^{-1}\be_{\bnu},\qquad
B_{(ii),\pOrder+2}(y,\bx) = \theta_{\pOrder+2,\bnu}(y|\bx)\bc_{y,\pOrder+2}^\Trans\bS_{y}^{-1}\be_{\mu},\\
B_{(iii),\pOrder+1, \qOrder+1}(y,\bx) & = \be_{\mu}^\Trans\bS_y^\inv \bc_{y,\pOrder+1} \bigg(\sum_{|\bmm|=\qOrder+1} \theta_{\pOrder+1,\bmm}(y|\bx) \bc_{\bx,\bmm}^\Trans\bigg)\bS_{\bx}^\inv  \be_{\bnu}.
\end{align*}
Note that the last term, $h^{\pOrder+\qOrder+2-\mu-|\bnu|}B_{(iii),\pOrder+1, \qOrder+1}(y,\bx)$, is present only if $\mu=\pOrder$ and $|\bnu|=\qOrder$.
\end{remark}

Next we study the leading variance of our estimator, defined as
\begin{align*}
\mathsf{V}_{\order} (\eval)&= \Var\Big[\frac{1}{n}\sumIN \K_{\mu,\bnu,h}^\circ\left( y_i,\bx_i;y,\bx \right)\Big].
\end{align*}

\begin{lem}[Variance]
\label{sa-lem:Variance}
Let Assumptions \ref{sa-assu:1 dgp}, \ref{sa-assu:2 kernel} and \ref{sa-assu:3 density} hold with $h \to 0$ and $nh^d / \log (n)\to\infty$. Then\\
(i) $\mu=0$ and $\theta_{0, \mathbf{0}}\neq 0$ or $1$:
\begin{align*}
\mathsf{V}_{0,\bnu}(\eval)
=& \frac{1}{nh^{d+2|\bnu|}} \theta_{0,\mathbf{0}}(y|\bx)(1-\theta_{0,\mathbf{0}}(y|\bx))\Big(\be_{\bnu}^\Trans \bS_{\bx}^{-1} \bT_{\bx} \bS_{\bx}^{-1} \be_{\bnu}\Big) + O\left( \frac{1}{nh^{d+2|\bnu|-1}} \right).
\end{align*}
(ii) $\mu=0$ and $\theta_{0, \mathbf{0}}= 0$ or $1$: $\mathsf{V}_{0,\bnu}(\eval)$ has the order $\frac{1}{nh^{d+2|\bnu|-1}}$.\\
(iii) $\mu > 0$:
\begin{align*}
\mathsf{V}_{\order} (\eval) &= \frac{1}{nh^{d + 2|\bnu|+2\mu-1}} \theta_{1,\mathbf{0}}(y|\bx)\Big(\be_{\mu}^\Trans\bS_y^\inv \bT_{y}\bS_y^\inv \be_{\mu}\Big)\Big(\be_{\bnu}^\Trans \bS_{\bx}^\inv  \bT_{\bx}\bS_{\bx}^\inv \be_{\bnu}\Big) + O\left(\frac{1}{nh^{d+2\mu+2|\bnu|-2}}\right).
\end{align*}
\end{lem}

For future reference, we will define
\begin{align*}
\mathtt{r_V} &= \sqrt{\frac{1}{nh^{d + 2|\bnu|+2\mu-1}}}.
\end{align*}

\begin{remark}[Vanishing boundary variance when $\mu=0$]
In case (ii), the true conditional distribution function is 0 or 1, which is why the leading variance shrinks faster. We do not provide a formula as the leading variance in this case takes a complicated form.
\end{remark}

Now, we propose two estimators for the variance that are valid for all three cases of Lemma \ref{sa-lem:Variance},
and hence will be useful for establishing a self-normalized distributional approximation later.
Define
\begin{align*}
\check{\mathsf{V}}_{\order}(\eval) &= \frac{1}{n^2}\sumIN \check{\K}_{\mu,\bnu,h}^\circ\left( y_i,\bx_i; y,\bx \right)^2,\quad
\hat{\mathsf{V}}_{\order} (\eval)= \frac{1}{n^2}\sumIN \hat{\K}_{\mu,\bnu,h}^\circ\left( y_i,\bx_i; y,\bx \right)^2.
\end{align*}

Note that $\hat{\mathsf{V}}_{\order}(\eval)$ is simply the plug-in variance estimator for $\hat{\theta}_{\order}(y|\bx)$
and $\check{\mathsf{V}}_{\order}(\eval)$ is the plug-in variance estimator for $\check{\theta}_{\order}(y|\bx)$.
The next lemma provides pointwise convergence results for the two variance estimators.
\begin{lem}[Variance estimation]\label{sa-lem:Variance estimation}
Let Assumptions \ref{sa-assu:1 dgp}, \ref{sa-assu:2 kernel} and \ref{sa-assu:3 density} hold with $h \to 0$ and $nh^{d+1} / \log (n)\to\infty$. In addition, $\theta_{0,\bnu}$ exists and is continuous for all $|\bnu|\leq \qOrder+1$.
Then\\
(i) $\mu=0$ and $\theta_{0, \mathbf{0}}\neq 0$ or $1$:
\begin{align*}
& \Big| \frac{\check{\mathsf{V}}_{0,\bnu}(\eval) - {\mathsf{V}}_{0,\bnu}(\eval)}{{\mathsf{V}}_{0,\bnu}(\eval)} \Big| = \Op\Big( h^{\qOrder + 1} + \sqrt{\frac{\log (n)}{nh^{d}}}\Big).
\end{align*}
(ii) $\mu > 0$, or $\theta_{0, \mathbf{0}}= 0$ or $1$:
\begin{align*}
& \Big| \frac{\check{\mathsf{V}}_{\order}(\eval) - {\mathsf{V}}_{\order}(\eval)}{{\mathsf{V}}_{\order}(\eval)} \Big| = \Op\Big( h^{\qOrder+\frac{1}{2}} + \sqrt{\frac{\log (n)}{nh^{d+1}}}\Big).
\end{align*}
Let $G = F_y$, then the same conclusions hold for $\hat{\mathsf{V}}_{\mu,\bnu}(\eval)$.
\end{lem}

Next, we study the large-sample distributional properties of the infeasible, standardized statistic
\begin{align*}
\bar{\SProc}_{\order}(\eval)
&=
\frac{1}{n\sqrt{\mathsf{V}_{\order}(\eval)} }\sumIN \K_{\mu,\bnu,h}^{\circ}\left( y_i,\bx_i; y,\bx \right).
\end{align*}
Note that this is equivalent to the scaled asymptotic linear representation of the estimator.

\begin{theorem}[Asymptotic normality]\label{sa-thm:Asymptotic normality}
Let Assumptions \ref{sa-assu:1 dgp}, \ref{sa-assu:2 kernel} and \ref{sa-assu:3 density} hold with $h \to 0$. Then
\begin{align*}
\sup_{u \in \mathbb{R}} \Big| \Prob \left[ \bar{\SProc}_{\order} (\eval) \leq u \right] - \Phi(u) \Big| = O\left(\mathtt{r_{BE}}\right),\quad \text{where }\mathtt{r_{BE}}=
\begin{cases}
\frac{1}{\sqrt{nh^d}} & \text{if $\mu=0$,\ \  and $\theta_{0, \mathbf{0}}\neq 0$ or $1$}\\
\frac{1}{\sqrt{nh^{d+1}}} & \text{if $\mu > 0$, or if $\theta_{0, \mathbf{0}}= 0$ or $1$}
\end{cases}.
\end{align*}
\end{theorem}

While the theorem focuses on asymptotic normality of the infeasible t-statistic, $\bar{\SProc}_{\order}^{\circ} (\eval)$, we show in the following remark that similar conclusions can be made for the t-statistics constructed with the estimators, $\hat{\theta}_{\order}(y|\bx)$ and $\check{\theta}_{\order}(y|\bx)$.

\begin{remark}[Asymptotic normality of standardized statistics]\label{remark:pointwise berry-esseen}
We first introduce the statistic
\begin{align*}
\check{\SProc}_{\order}^\circ(\eval) = \frac{\check{\theta}_{\order}(y|\bx)-\Expectation \left[ \check{\theta}_{\order}(y|\bx) | \bX \right]}{\sqrt{{\mathsf{V}}_{\order}(\eval)}},
\end{align*}
which is based on $\check{\theta}_{\order}(y|\bx)$. (In the main paper we directly center all statistics at the target parameter $\theta_{\mu,\bnu}$. For clarity, however, we will separate the discussion on distributional convergence from the smoothing bias in this supplementary material. This is reflected by the superscript ``circle.'') By combining the results of Lemmas~\ref{sa-lem:Taking hat off} and~\ref{sa-lem:Variance}, we have
\begin{align*}
\sup_{y \in \mathcal{Y}, \bx \in \mathcal{X}}\left| \check{\SProc}_{\order}^\circ(\eval) - \bar{\SProc}_{\order} (\eval) \right| = \Oconc\Big(\frac{\log (n)}{\sqrt{nh^d}}\Big).
\end{align*}
As a result,
\begin{align*}
\sup_{u \in \mathbb{R}} \left| \Prob \left[ \check{\SProc}_{\order}^\circ (\eval) \leq u \right] - \Phi(u) \right| = O\Big( \frac{\log (n)}{\sqrt{nh^d}} + \mathtt{r_{BE}}\Big).
\end{align*}

To present the pointwise distributional approximation result for the estimator $\hat{\theta}_{\order}(y|\bx)$, we define the following statistic
\begin{align*}
\hat{\SProc}_{\order}^\circ(\eval)
&=
\frac{1}{nh^{d+\mu+|\bnu|}\sqrt{\mathsf{V}_{\order}(\eval)} }\sumIN \be_{\mu}^\Trans \hat{\bS}_y^\inv \bigg[\frac{1}{n}\sum_{j=1}^n \Big[\Indicator(y_i\leq y_j) - F(y_j|\bx_i) \Big] \frac{1}{h}\bP \Big(\frac{y_j-y}{h}\Big)\bigg]\\ &\qquad\qquad\qquad\qquad\qquad\qquad\qquad\qquad\qquad\qquad\bQ  \left(\frac{\bx_i-\bx}{h}\right)^\Trans \hat{\bS}_{\bx}^\inv  \be_{\bnu}.
\end{align*}
It is worth mentioning that $\hat{\SProc}_{\order}^\circ(\eval)$ is not exactly
centered and therefore, it is not mean zero. Nevertheless, by the results of Lemmas~\ref{sa-lem:Taking hat off} and~\ref{sa-lem:Variance}, and the concentration inequality for second order U-statistics in Equation (3.5) of \cite{Gine-Latala-Zinn_2000_Ustat} (Lemmas 7 and 8 in the main paper), we have
\begin{align*}
\sup_{y \in \mathcal{Y}, \bx \in \mathcal{X}}\left| \hat{\SProc}_{\order}^\circ(\eval) - \check{\SProc}_{\order}^\circ (\eval) \right| =\Oconc \Big(\frac{\log (n)}{\sqrt{nh^2}}\Big).
\end{align*}
Then we can conclude that the coverage error satisfies
\begin{align*}
\sup_{u \in \mathbb{R}} \left| \Prob \left[ \hat{\SProc}_{\order}^\circ (\eval) \leq u \right] - \Phi(u) \right| = O\Big(\frac{\log (n)}{\sqrt{nh^{d\vee 2}}} + \mathtt{r_{BE}}\Big).
\end{align*}
\end{remark}

\section{Uniform large-sample properties}\label{sa-sec:Uniform Large-sample Properties}

To conduct statistical inference on the entire function $\theta_{\order}$, such as constructing confidence bands or testing shape restrictions, we need uniform distributional approximations to our estimators. In this section, we will consider large-sample properties of our estimator which hold uniformly on $\mathcal{Y} \times \mathcal{X} = [0,1]^{d+1}$. In the following remark, we demonstrate that the local sample size is uniformly large on the support $\mathcal{Y} \times \mathcal{X}$.

\begin{remark}[Local sample size]\label{SA-remark: local sample size}
Consider an evaluation point $(y,\bx)$ in $\mathcal{Y} \times \mathcal{X}$. We can define the local sample size by
\[
n_{y,\bx} = \sumIN \Indicator(|y_i-y|\leq \mathfrak{c}_1h)\Indicator(|\bx_i-\bx|\leq \mathfrak{c}_1h).
\]
We employed the Euclidean norm in the definition, which is innocuous for our purposes, as all norms are equivalent in finite dimensional spaces. For this reason, we also introduced the constant $\mathfrak{c}_1$. The purpose of this remark is to provide a uniform control on the local sample size. In particular, we have the following result: for some positive constant $\mathfrak{c}_2$ and any shrinking sequence $\mathfrak{r}$,
\begin{align*}
\Indicator\left(\inf_{y\in\Y,\bx\in\X}\left| n_{y,\bx} \right| < \mathfrak{c}_2 \frac{\log (n)}{nh^{d+1}}\right) = \Oconc(\mathfrak{r}).
\end{align*}
The result above builds on the lemma:
\begin{lem}[Probabilistic bound on the smallest multinomial cell]\label{SA-lem: smallest cell}
Let $\bz=(z_1,z_2,\dots,z_{J_n})^\Trans$ follow a multinomial distribution with parameters $n$ (number of trials), $J_n$ (number of cells), and $1/J_n$ (probability for each cell),  $\delta_n\in(0,1)$, and $\pi_n = n/(J_n\log (n))$. If $\delta_n^2\pi_n\to \infty$, then for any $\mathfrak{c}_1>0$,
\begin{align*}
\limsup_{n\to\infty}n^{\mathfrak{c}_1}\Prob\left[ \min_{1\leq j\leq J_n}z_j < (1-\delta_n) \frac{n}{J_n} \right] <\infty.
\end{align*}
\end{lem}
\end{remark}

We now establish the uniform convergence rate of our estimator.

\begin{lem}[Uniform rate of convergence]\label{sa-lem:Uniform rate of convergence}
Let Assumptions \ref{sa-assu:1 dgp}, \ref{sa-assu:2 kernel} and \ref{sa-assu:3 density} hold with $h \to 0$ and $nh^{d+1} / \log (n)\to\infty$. In addition, $\theta_{\mu',\bnu'}$ exists and is continuous for all $\mu' + |\bnu'|= \max\{ \qOrder+1+\mu,\ \pOrder+1+|\bnu| \}$. Then\\
(i) $\mu=0$:
\begin{align*}
&\sup_{y\in \mathcal{Y}, \bx \in \mathcal{X}}\left| \check{\theta}_{0,\bnu}(y|\bx) - \theta_{0,\bnu}(y|\bx) \right| = \Oconc\Big(h^{\qOrder+1-|\bnu|} +  h^{\pOrder+1} + \sqrt{\frac{\log (n)}{nh^{d+2|\bnu|}}}\Big);
\end{align*}
(ii) $\mu>0$:
\begin{align*}
&\sup_{y\in \mathcal{Y}, \bx \in \mathcal{X}}\left| \check{\theta}_{\order}(y|\bx) - \theta_{\order}(y|\bx) \right| = \Oconc\Big(h^{\qOrder+1-|\bnu|} +  h^{\pOrder+1-\mu} + \sqrt{\frac{\log (n)}{nh^{d+2\mu+2|\bnu|-1}}}\Big).
\end{align*}
The same conclusions hold for $\hat{\theta}_{\order}(y|\bx)$.
\end{lem}

In the next lemma, we characterize the uniform convergence rate of the variance estimators introduced in the previous section.

\begin{lem}[Uniform variance estimation]\label{sa-lem:Uniform variance estimation}
Let Assumptions \ref{sa-assu:1 dgp}, \ref{sa-assu:2 kernel} and \ref{sa-assu:3 density} hold with
$h \to 0$ and
$nh^{d+1} / \log (n)\to\infty$. In addition, $\theta_{0,\bnu}$ exists and is continuous for all $|\bnu|\leq \qOrder+1$. Then
\begin{align*}
\sup_{y \in \mathcal{Y}, \bx \in \mathcal{X}}\bigg| \frac{\check{\mathsf{V}}_{\order}(\eval)  - {\mathsf{V}}_{\order}(\eval) }{{\mathsf{V}}_{\order}(\eval) } \bigg| = \Oconc\left(\mathtt{r_{VE}}\right),\quad \text{where }\mathtt{r_{VE}}= h^{ \qOrder+\frac{1}{2}} + \sqrt{\frac{\log (n)}{nh^{d+1}}}.
\end{align*}
Let $G = F_y$, then the same conclusions hold for $\hat{\mathsf{V}}_{\mu,\bnu}(\eval)$.
\end{lem}

Now, we introduce the Studentized processes for each of the estimators,
$\hat{\theta}_{\order}$ and $\check{\theta}_{\order}$:
\begin{align*}
&\check{\TProc}_{\order}^\circ(\eval) = \sqrt{\frac{\mathsf{V}_{\order} (\eval)}{\check{\mathsf{V}}_{\order} (\eval)}}\check{\SProc}_{\order}^\circ(\eval),\qquad \hat{\TProc}_{\order}^\circ(\eval) = \sqrt{\frac{\mathsf{V}_{\order} (\eval)}{\hat{\mathsf{V}}_{\order} (\eval)}}\hat{\SProc}_{\order}^\circ(\eval).
\end{align*}
In the following lemma we study the error that arises from the Studentization of our estimators.

\begin{lem}[Studentization error]\label{sa-lem:uniform Studentization error}
Let Assumptions \ref{sa-assu:1 dgp}, \ref{sa-assu:2 kernel} and \ref{sa-assu:3 density} hold with
$h \to 0$,
$nh^{d+1} / \log (n)\to\infty$. In addition, $\theta_{0,\bnu}$ exists and is continuous for all $|\bnu|\leq \qOrder+1$. Then
\begin{align*}
\sup_{y \in \Y, \bx \in \X}\left| \check{\TProc}_{\order}^\circ(\eval) - \check{\SProc}_{\order}^\circ(\eval) \right| =\Oconc\left(\mathtt{r_{SE}}\right),\qquad \text{where } \mathtt{r_{SE}}=\sqrt{\log (n) }\mathtt{r_{VE}}.
\end{align*}
The same holds for $\hat{\TProc}_{\order}^\circ(\eval) - \hat{\SProc}_{\order}^\circ(\eval)$.
\end{lem}

Our next goal is to establish a uniform normal approximation to the process $\bar{\SProc}_{\order}(y,\bx)$. See Appendix A.4 of the main paper for important properties of the equivalent kernel.

\begin{theorem}[Strong approximation]\label{sa-lem:strong approximation}
Let Assumptions \ref{sa-assu:1 dgp}, \ref{sa-assu:2 kernel} and \ref{sa-assu:3 density} hold with $h \to 0$ and $nh^{d+1} / \log (n)\to\infty$. Also assume $\mu \geq 1$. Define
\begin{align*}
\mathtt{r_{SA}} = \left(\frac{\log^{d+1} n}{nh^{d+1}}\right)^{\frac{1}{2d+2}}.
\end{align*}
Then there exist two centered processes, $\bar{\SProc}^{\prime}_{\order}(\eval) $ and  $\GProc_{\order}(\eval) $, such that (i) $\bar{\SProc}_{\order}(\eval) $ and $\bar{\SProc}^{\prime}_{\order}(\eval) $ have the same distribution, (ii) $\GProc_{\order}(\eval)$ is a Gaussian process and has the same covariance kernel as $\bar{\SProc}_{\order}(\eval)$, and (iii)
\begin{align*}
\sup_{y \in \Y, \bx \in \X}\left|\bar{\SProc}^{\prime}_{\order}(\eval)-\GProc_{\order}(\eval)\right|  = \Oconc\left(\mathtt{r_{SA}}\right).
\end{align*}
\end{theorem}

The Gaussian approximation in the above lemma is not feasible, as its covariance kernel depends on unknowns.
To be more precise, the covariance kernel takes the form
\begin{align*}
&\ \rho_{\order}(y,\bx,y',\bx') = \Cov\left[ \bar{\SProc}_{\order}(\eval), \bar{\SProc}_{\order}(y',\bx') \right] = \frac{\mathsf{C}_{\order}(y,\bx,y',\bx')}{ \sqrt{\mathsf{V}_{\order}(\eval)\mathsf{V}_{\order}(y',\bx')}},
\end{align*}
where
\begin{align*}
&\ \mathsf{C}_{\order}(y,\bx,y',\bx')
=\ \frac{1}{n}\Cov\Big[\K_{\mu,\bnu,h}^\circ\left( y_i,\bx_i;y,\bx \right),\K_{\mu,\bnu,h}^\circ\left( y_i,\bx_i;y',\bx' \right)  \Big].
\end{align*}
We consider two estimators of the covariance kernel
\begin{align*}
\check{\rho}_{\order}(y,\bx,y',\bx') &= \frac{\check{\mathsf{C}}_{\order}(y,\bx,y',\bx')}{ \sqrt{\check{\mathsf{V}}_{\order}(\eval)\check{\mathsf{V}}_{\order}(y',\bx')}},\quad
\hat{\rho}_{\order}(y,\bx,y',\bx') = \frac{\hat{\mathsf{C}}_{\order}(y,\bx,y',\bx')}{ \sqrt{\hat{\mathsf{V}}_{\order}(\eval)\hat{\mathsf{V}}_{\order}(y',\bx')}},
\end{align*}
and
\begin{align*}
\check{\mathsf{C}}_{\mu,\bnu}(y,\bx,y',\bx') &= \frac{1}{n^2}\sumIN \check{\K}_{\mu,\bnu,h}^\circ\left( y_i,\bx_i; y,\bx \right)\check{\K}_{\mu,\bnu,h}^\circ\left( y_i,\bx_i; y',\bx' \right)\\
\hat{\mathsf{C}}_{\mu,\bnu}(y,\bx,y',\bx') &= \frac{1}{n^2}\sumIN \hat{\K}_{\mu,\bnu,h}^\circ\left( y_i,\bx_i; y,\bx \right)\hat{\K}_{\mu,\bnu,h}^\circ\left( y_i,\bx_i; y',\bx' \right).
\end{align*}

\begin{lem}[Uniform consistency of the correlation estimator]\label{sa-lem:uniform covariance estimation}
Let Assumptions \ref{sa-assu:1 dgp}, \ref{sa-assu:2 kernel} and \ref{sa-assu:3 density} hold with $h \to 0$ and $nh^{d+1} / \log (n)\to\infty$. In addition, $\theta_{0,\bnu}$ exists and is continuous for all $|\bnu|\leq \qOrder+1$. Then
\begin{align*}
\sup_{y,y' \in \Y, \bx,\bx' \in \X}\left| \check{\rho}_{\order}(y,\bx,y',\bx') - \rho_{\order}(y,\bx,y',\bx') \right|  =\Oconc\left(\mathtt{r_{VE}}\right),
\end{align*}
where $\mathtt{r_{VE}}$ is defined in Lemma \ref{sa-lem:Uniform variance estimation}. Let $G = F_y$, then the same conclusion holds for $\hat{\rho}_{\order}(y,\bx,y',\bx')$.
\end{lem}

\begin{lem}[Gaussian comparison]\label{sa-lem:feasible Gaussian approx}
Let Assumptions \ref{sa-assu:1 dgp}, \ref{sa-assu:2 kernel} and \ref{sa-assu:3 density} hold with $h \to 0$, $nh^{d+1} / \log (n)\to\infty$. In addition, $\theta_{0,\bnu}$ exists and is continuous for all $|\bnu|\leq \qOrder+1$. Then conditional on the data there exists a centered Gaussian process, $\check{\GProc}_{\order}(\eval) $ with unit variance and correlation function $\check{\rho}_{\order}$, and another centered Gaussian process, $\hat{\GProc}_{\order}(\eval) $ with unit variance and correlation kernel $\hat{\rho}_{\order}$, such that
\begin{align*}
&\sup_{u \in \mathbb{R}}\Big| \Prob\Big[  \sup_{y \in \Y, \bx \in \X}|\check{\GProc}_{\order}(\eval) |\leq u \Big| \bY,\bX \Big] - \Prob\Big[ \sup_{y \in \Y, \bx \in \X}|\GProc_{\order}(\eval) |\leq u  \Big] \Big| = \Op \left(\log(n)\sqrt{\mathtt{r_{VE}}}\right),\\
&\sup_{u \in \mathbb{R}}\Big| \Prob\Big[  \sup_{y \in \Y, \bx \in \X}|\hat{\GProc}_{\order}(\eval) |\leq u \Big| \bY,\bX \Big] - \Prob\Big[ \sup_{y \in \Y, \bx \in \X}|\GProc_{\order}(\eval) |\leq u  \Big] \Big| = \Op \left(\log(n)\sqrt{\mathtt{r_{VE}}}\right).
\end{align*}
\end{lem}

\begin{theorem}[Feasible normal approximation]\label{sa-thm:feasible Gaussian approx}
Let Assumptions \ref{sa-assu:1 dgp}, \ref{sa-assu:2 kernel} and \ref{sa-assu:3 density} hold with $h \to 0$ and $nh^{d+1} / \log (n)\to\infty$. In addition, $\theta_{0,\bnu}$ exists and is continuous for all $|\bnu|\leq \qOrder+1$.
Also assume $\mu\geq 1$. Then
\begin{align*}
&\sup_{u \in \mathbb{R}}\Big| \Prob\Big[ \sup_{y \in \Y, \bx \in \X}|\check{\TProc}_{\order}^\circ (\eval) |\leq u  \Big] -
\Prob\Big[ \sup_{y \in \Y, \bx \in \X}|\check{\GProc}_{\order}(\eval)|\leq u \Big| \bY,\bX \Big]  \Big|\\
&\qquad\qquad = \Op \left(\sqrt{\log(n)}\mathtt{r_{SA}} + \log(n)\sqrt{\mathtt{r_{VE}}}\right),\\
&\sup_{u \in \mathbb{R}}\Big| \Prob\Big[ \sup_{y \in \Y, \bx \in \X}|\hat{\TProc}_{\order}^\circ (\eval) |\leq u  \Big] -
\Prob\Big[ \sup_{y \in \Y, \bx \in \X}|\hat{\GProc}_{\order}(\eval)|\leq u \Big| \bY,\bX \Big]  \Big|\\
&\qquad\qquad   = \Op \left(\sqrt{\log(n)}\mathtt{r_{SA}} + \log(n)\sqrt{\mathtt{r_{VE}}}\right).
\end{align*}
\end{theorem}

\section{Applications}

\subsection{Confidence bands}

A natural corollary of Theorem~\ref{sa-thm:feasible Gaussian approx} is that one can employ critical values computed from $\check{\GProc}_{\order}(\eval) $ and $\hat{\GProc}_{\order}(\eval)$ to construct confidence bands. To be very precise, define
\begin{align*}
\check{\mathtt{cv}}_{\order}(\alpha) &= \inf\Big\{ u:\ \Prob\Big[  \sup_{y \in \Y, \bx \in \X}|\check{\GProc}_{\order}(\eval)|\leq u \Big| \bY,\bX \Big] \geq 1-\alpha  \Big\},\\
\hat{\mathtt{cv}}_{\order}(\alpha) &= \inf\Big\{ u:\ \Prob\Big[ \sup_{y \in \Y, \bx \in \X}|\hat{\GProc}_{\order}(\eval)|\leq u \Big| \bY,\bX \Big] \geq 1-\alpha  \Big\}.
\end{align*}
Then level $(1-\alpha)$ confidence bands can be constructed as
\begin{align*}
\check{\mathcal{C}}_{\order}(1-\alpha) &= \Big\{ \check{\theta}_{\order}(y|\bx) \pm  \check{\mathtt{cv}}_{\order}(\alpha)\sqrt{\check{\mathsf{V}}_{\order}(\eval)}:\ (y,\bx)\in\mathcal{Y}\times \mathcal{X} \Big\},\\
\hat{\mathcal{C}}_{\order}(1-\alpha) &= \Big\{ \hat{\theta}_{\order}(y|\bx) \pm  \hat{\mathtt{cv}}_{\order}(\alpha)\sqrt{\hat{\mathsf{V}}_{\order}(\eval)}:\ (y,\bx)\in\mathcal{Y}\times \mathcal{X} \Big\},
\end{align*}
whose coverage error is given in the following theorem.

\begin{theorem}[Confidence band]\label{sa-thm:confidence bands}
Consider the setting of Theorem~\ref{sa-thm:feasible Gaussian approx}. In addition, $\theta_{\mu',\bnu'}$ exists and is continuous for all $\mu' + |\bnu'|= \max\{ \qOrder+1+\mu,\ \pOrder+1+|\bnu| \}$. Then
\begin{align*}
\Prob\left[ \theta_{\mu,\bnu}(y|\bx) \in \check{\mathcal{C}}_{\order}(1-\alpha),\ \forall (y,\bx)\in\mathcal{Y}\times \mathcal{X} \right] \geq 1-\alpha - O\Big(\sqrt{\log(n)}\Big(\mathtt{r_{SA}} + \frac{\mathtt{r_{B}}}{\mathtt{r_{V}}} \Big) + \log(n)\sqrt{\mathtt{r_{VE}}}\Big),\\
\Prob\left[ \theta_{\mu,\bnu}(y|\bx) \in \hat{\mathcal{C}}_{\order}(1-\alpha),\ \forall (y,\bx)\in\mathcal{Y}\times \mathcal{X} \right] \geq 1-\alpha - O\Big(\sqrt{\log(n)}\Big(\mathtt{r_{SA}} + \frac{\mathtt{r_{B}}}{\mathtt{r_{V}}} \Big) + \log(n)\sqrt{\mathtt{r_{VE}}}\Big).
\end{align*}
\end{theorem}

\subsection{Parametric specification testing}

In applications, it is not uncommon to estimate conditional densities or higher-order derivatives by specifying a parametric family of distributions. While such parametric restrictions may provide reasonable approximations, it is still worthwhile to conduct specification testing. To be specific, assume the researcher postulates the following class
\begin{align*}
\big\{ \theta_{\mu,\bnu}(y|\bx;\bgamma):\ \bgamma\in \mathsf{\Gamma}_{\mu,\bnu} \big\},
\end{align*}
where $\mathsf{\Gamma}_{\mu,\bnu}$ is some compact parameter space. We abstract away from the specifics of the estimation technique, and assume that the researcher also picks some estimator (maximum likelihood, minimum distance, etc.) $\hat{\bgamma}$. Under fairly mild conditions, the estimator will converge in probability to some (possibly pseudo-true) parameter $\bar{\bgamma}$ in the parameter space $\mathsf{\Gamma}_{\mu,\bnu}$. As before, we will denote the true parameter as $\theta_{\mu,\bnu}(y|\bx)$, and consider the following competing hypotheses:
\begin{align*}
\mathsf{H}_0:\ \theta_{\mu,\bnu}(y|\bx;\bar{\bgamma})=\theta_{\mu,\bnu}(y|\bx)\quad \text{vs}.\quad
\mathsf{H}_1:\ \theta_{\mu,\bnu}(y|\bx;\bar{\bgamma})\neq\theta_{\mu,\bnu}(y|\bx).
\end{align*}
The test statistics we employ takes the following form
\begin{align*}
\check{\TProc}_{\mathtt{PS}}(\eval) = \frac{\check{\theta}_{\order}(y|\bx) - \theta_{\order}(y|\bx;\hat{\bgamma}) }{\sqrt{\check{\mathsf{V}}_{\order}(\eval)}},\quad
\hat{\TProc}_{\mathtt{PS}}(\eval) = \frac{\hat{\theta}_{\order}(y|\bx) - \theta_{\order}(y|\bx;\hat{\bgamma}) }{\sqrt{\hat{\mathsf{V}}_{\order}(\eval)}}.
\end{align*}

\begin{theorem}[Parametric specification testing]\label{sa-thm:HT-specification}
Consider the setting of Theorem~\ref{sa-thm:feasible Gaussian approx}. In addition, $\theta_{\mu',\bnu'}$ exists and is continuous for all $\mu' + |\bnu'|= \max\{ \qOrder+1+\mu,\ \pOrder+1+|\bnu| \}$. Assume the parametric estimate satisfies
\begin{align*}
\sup_{ y \in \Y, \bx \in \X}\left| \theta_{\mu,\bnu}(y|\bx;\hat{\bgamma}) - \theta_{\mu,\bnu}(y|\bx;\bar{\bgamma}) \right| = \Oconc\left(\mathtt{r_{PS}}\right),
\end{align*}
for some $\mathtt{r_{PS}}$. Then under the null hypothesis,
\begin{align*}
\Prob\Big[ \sup_{y \in \Y, \bx \in \X}|\check{\TProc}_{\mathtt{PS}}(\eval) | > \check{\mathtt{cv}}_{\order}(\alpha) \Big] \leq \alpha  + O\Big(\sqrt{\log(n)}\Big(\mathtt{r_{SA}} + \frac{\mathtt{r_{B}}+\mathtt{r_{PS}}}{\mathtt{r_{V}}} \Big) + \log(n)\sqrt{\mathtt{r_{VE}}}\Big),\\
\Prob\Big[ \sup_{y \in \Y, \bx \in \X}|\hat{\TProc}_{\mathtt{PS}}(\eval) | > \hat{\mathtt{cv}}_{\order}(\alpha) \Big] \leq \alpha  + O\Big(\sqrt{\log(n)}\Big(\mathtt{r_{SA}} + \frac{\mathtt{r_{B}}+\mathtt{r_{PS}}}{\mathtt{r_{V}}} \Big) + \log(n)\sqrt{\mathtt{r_{VE}}}\Big).
\end{align*}
\end{theorem}

\subsection{Testing shape restrictions}

Now consider shape restrictions on the conditional density or its derivatives. Let $c(y,\bx)$ be a pre-specified function, and we study the following one-sided competing hypotheses.
\begin{align*}
    \mathsf{H}_0 :\theta_{\mu,\bnu}(y|\bx) \leq c(y,\bx)
    \quad \text{vs.} \quad
    \mathsf{H}_1 :\theta_{\mu,\bnu}(y|\bx) > c(y,\bx).
\end{align*}
The statistic we employ takes the form
\begin{align*}
\check{\TProc}_{\mathtt{SR}}(\eval) = \frac{\check{\theta}_{\order}(y|\bx) - c(y,\bx) }{\sqrt{\check{\mathsf{V}}_{\order}(\eval)}},\quad
\hat{\TProc}_{\mathtt{SR}}(\eval) = \frac{\hat{\theta}_{\order}(y|\bx) - c(y,\bx) }{\sqrt{\hat{\mathsf{V}}_{\order}(\eval)}}.
\end{align*}
and we will reject the null hypothesis if the test statistic exceeds a critical value.
\begin{theorem}[Shape restriction testing]\label{sa-thm:HT-restriction}
Consider the setting of Theorem~\ref{sa-thm:feasible Gaussian approx}. In addition, $\theta_{\mu',\bnu'}$ exists and is continuous for all $\mu' + |\bnu'|= \max\{ \qOrder+1+\mu,\ \pOrder+1+|\bnu| \}$. Then under the null hypothesis,
\begin{align*}
\Prob\Big[ \sup_{y \in \Y, \bx \in \X}\check{\TProc}_{\mathtt{SR}}(\eval) > \check{\mathtt{cv}}_{\order}(\alpha) \Big] \leq \alpha  + O\Big(\sqrt{\log(n)}\Big(\mathtt{r_{SA}} + \frac{\mathtt{r_{B}}}{\mathtt{r_{V}}} \Big) + \log(n)\sqrt{\mathtt{r_{VE}}}\Big),\\
\Prob\Big[ \sup_{y \in \Y, \bx \in \X}\hat{\TProc}_{\mathtt{SR}}(\eval) > \hat{\mathtt{cv}}_{\order}(\alpha) \Big] \leq \alpha  + O\Big(\sqrt{\log(n)}\Big(\mathtt{r_{SA}} + \frac{\mathtt{r_{B}} }{\mathtt{r_{V}}} \Big) + \log(n)\sqrt{\mathtt{r_{VE}}}\Big).
\end{align*}
\end{theorem}

\section{Bandwidth selection}
We assume throughout this section that $\mu >0$. Using the bias expression derived in Lemma~\ref{sa-lem:Bias}, and the leading variance is as characterized in Lemma~\ref{sa-lem:Variance},
we can derive precise expressions for bandwidth selection.
\subsection{Pointwise asymptotic MSE minimization}
\label{sa-sec:mse_bw}
Following from \cite{Fan-Gijbels_1996_Book},
the pointwise MSE-optimal bandwidth is defined as a minimizer of the following optimization problem
\begin{align*}
h^*_{\pOrder,\qOrder,\mu, \bnu}(\eval)=\argmin_{h>0}\left[ \mathsf{V}_\order (\eval) + \mathsf{B}^{2}_{\order}(\eval) \right]
\end{align*}
The solution to this equation
gives an MSE-optimal bandwidth that depends on (i) the order of the polynomials, (ii) the order of the derivative to be estimated, and (iii) the position of the evaluation point. \

\subsubsection*{Case 1:  $\qOrder-|\bnu| = \pOrder-\mu$, odd}

In this case, both the leading bias constants, $B_{(i),\qOrder+1}(y,\bx)$ and $B_{(ii),\pOrder+1}(y,\bx)$, are nonzero. Therefore, the MSE-optimal bandwidth is
\begin{align*}
h^*_{\pOrder,\qOrder,\mu, \bnu}(\eval) &=
\argmin_{h>0}\Big[ \frac{1}{nh^{d + 2|\bnu|+2\mu-1}}V_\order(\eval)
+ h^{\pOrder+\qOrder + 2-\mu-|\bnu|}
\left(B_{(i),\qOrder+1}(y,\bx) +B_{(ii),\pOrder+1}(y,\bx) \right)^{2}\Big]
\\ &
= \Big[ \frac{(d + 2|\bnu|+2\mu-1)V_\order(\eval)}
{(\pOrder+\qOrder + 2-\mu-|\bnu|)\left(B_{(i),\qOrder+1}(y,\bx) +B_{(ii),\pOrder+1}(y,\bx) \right)^{2}}\frac{1}{n} \Big]
^{\frac{1}{d + \pOrder+\qOrder+ |\bnu|+\mu+1 }}.
\end{align*}

\subsubsection*{Case 2:  $\qOrder-|\bnu| = \pOrder-\mu$, even; either $\bx$ or $y$ is at or near the boundary}

In this case, at least one of the leading bias constants, $B_{(i),\qOrder+1}(y,\bx)$ and $B_{(ii),\pOrder+1}(y,\bx)$, is nonzero. Therefore, the MSE-optimal bandwidth is the same as in Case 1:
\begin{align*}
h^*_{\pOrder,\qOrder,\mu, \bnu}(\eval) &=
\argmin_{h>0}\Big[ \frac{1}{nh^{d + 2|\bnu|+2\mu-1}}V_\order(\eval)
+ h^{\pOrder+\qOrder + 2-\mu-|\bnu|}
\left(B_{(i),\qOrder+1}(y,\bx) +B_{(ii),\pOrder+1}(y,\bx) \right)^{2}\Big]
\\ &
= \Big[ \frac{(d + 2|\bnu|+2\mu-1)V_\order(\eval)}
{(\pOrder+\qOrder + 2-\mu-|\bnu|)\left(B_{(i),\qOrder+1}(y,\bx) +B_{(ii),\pOrder+1}(y,\bx) \right)^{2}}\frac{1}{n} \Big]
^{\frac{1}{d + \pOrder+\qOrder+ |\bnu|+\mu+1 }}.
\end{align*}

\subsubsection*{Case 3:  $\qOrder-|\bnu| = \pOrder-\mu \neq 0$, even; both $\bx$ and $y$ are interior}

In this case, both leading bias constants are zero. Therefore, the MSE-optimal bandwidth will depend on higher-order bias terms:
\begin{align*}
h^*_{\pOrder,\qOrder,\mu, \bnu}(\eval) &=
\argmin_{h>0}\Big[ \frac{1}{nh^{d + 2|\bnu|+2\mu-1}}V_\order(\eval)
+ h^{\pOrder+\qOrder + 4-\mu-|\bnu|}
\left(B_{(i),\qOrder+2}(y,\bx) +B_{(ii),\pOrder+2}(y,\bx) \right)^{2}\Big]
\\ &
= \Big[ \frac{(d + 2|\bnu|+2\mu-1)V_\order(\eval)}
{(\pOrder+\qOrder + 4-\mu-|\bnu|)\left(B_{(i),\qOrder+2}(y,\bx) +B_{(ii),\pOrder+2}(y,\bx) \right)^{2}}\frac{1}{n} \Big]
^{\frac{1}{d + \pOrder+\qOrder+ |\bnu|+\mu+3 }}.
\end{align*}

\subsubsection*{Case 4:  $\qOrder-|\bnu| = \pOrder-\mu = 0$, even; both $\bx$ and $y$ are interior}

As in Case 3, both leading bias constants are zero. The difference, however, is that the leading bias will involve an extra term:
\begin{align*}
h^*_{\pOrder,\qOrder,\mu, \bnu}(\eval) &=
\argmin_{h>0}\Big[ \frac{1}{nh^{d + 2|\bnu|+2\mu-1}}V_\order(\eval)
\\
&\qquad\qquad\qquad+ h^{4}
\left(B_{(i),\qOrder+2}(y,\bx) +B_{(ii),\pOrder+2}(y,\bx) + B_{(iii),\pOrder+1, \qOrder+1}(y,\bx) \right)^{2}\Big]
\\ &
= \Big[ \frac{(d + 2|\bnu|+2\mu-1)V_\order(\eval)}
{4\left(B_{(i),\qOrder+2}(y,\bx) +B_{(ii),\pOrder+2}(y,\bx) + B_{(iii),\pOrder+1, \qOrder+1}(y,\bx) \right)^{2}}\frac{1}{n} \Big]
^{\frac{1}{d+ 2|\bnu|+2\mu+3 }}.
\end{align*}

\subsubsection*{Case 5:  $\qOrder-|\bnu| < \pOrder-\mu$, $\qOrder-|\bnu|$ odd}

In this case, the leading bias will involve only one term:
\begin{align*}
h^*_{\pOrder,\qOrder,\mu, \bnu}(\eval) &=
\argmin_{h>0}\Big[ \frac{1}{nh^{d + 2|\bnu|+2\mu-1}}V_\order(\eval)
+ h^{2\qOrder+2-2|\bnu|}B_{(i),\qOrder+1}(y,\bx)^{2}\Big]
\\ &
= \Big[ \frac{(d + 2|\bnu|+2\mu-1)V_\order(\eval)}
{(2\qOrder+2-2|\bnu|)B_{(i),\qOrder+1}(y,\bx)^{2}}\frac{1}{n} \Big]^{\frac{1}{d+ 2\qOrder + 2\mu+1 }}.
\end{align*}

\subsubsection*{Case 6:  $\qOrder-|\bnu| = \pOrder-\mu-1$, $\qOrder-|\bnu|$ even; $\bx$ is interior}

In this case, the leading bias will involve two terms:
\begin{align*}
h^*_{\pOrder,\qOrder,\mu, \bnu}(\eval) &=
\argmin_{h>0}\Big[ \frac{1}{nh^{d + 2|\bnu|+2\mu-1}}V_\order(\eval)
+ h^{\pOrder+\qOrder + 3-\mu-|\bnu|}
\left(B_{(i),\qOrder+2}(y,\bx) +B_{(ii),\pOrder+1}(y,\bx) \right)^{2}\Big]
\\ &
= \Big[ \frac{(d + 2|\bnu|+2\mu-1)V_\order(\eval)}
{(\pOrder+\qOrder + 3-\mu-|\bnu|)\left(B_{(i),\qOrder+2}(y,\bx) +B_{(ii),\pOrder+1}(y,\bx) \right)^{2}}\frac{1}{n} \Big]^{\frac{1}{d+ \pOrder + \qOrder + \mu + |\bnu|+2 }}.
\end{align*}

\subsubsection*{Case 7:  $\qOrder-|\bnu| < \pOrder-\mu-1$, $\qOrder-|\bnu|$ even; $\bx$ is interior}

In this case, the leading bias will involve only one term:
\begin{align*}
h^*_{\pOrder,\qOrder,\mu, \bnu}(\eval) &=
\argmin_{h>0}\Big[ \frac{1}{nh^{d + 2|\bnu|+2\mu-1}}V_\order(\eval)
+ h^{2\qOrder + 4 -2|\bnu|}
B_{(i),\qOrder+2}(y,\bx)^{2}\Big]
\\ &
= \Big[ \frac{(d + 2|\bnu|+2\mu-1)V_\order(\eval)}
{(2\qOrder + 4 -2|\bnu|)B_{(i),\qOrder+2}(y,\bx)^{2}}\frac{1}{n} \Big]^{\frac{1}{d+ 2\qOrder + 2\mu + 3 }}.
\end{align*}

\subsubsection*{Case 8:  $\qOrder-|\bnu| > \pOrder-\mu$, $\pOrder-\mu$ odd}

In this case, the leading bias will involve only one term:
\begin{align*}
h^*_{\pOrder,\qOrder,\mu, \bnu}(\eval) &=
\argmin_{h>0}\Big[ \frac{1}{nh^{d + 2|\bnu|+2\mu-1}}V_\order(\eval)
+ h^{2\pOrder+2-2\mu}B_{(ii),\pOrder+1}(y,\bx)^{2}\Big]
\\ &
= \Big[ \frac{(d + 2|\bnu|+2\mu-1)V_\order(\eval)}
{(2\pOrder+2-2\mu)B_{(ii),\pOrder+1}(y,\bx)^{2}}\frac{1}{n} \Big]^{\frac{1}{d+ 2\pOrder + 2|\bnu|+1 }}.
\end{align*}

\subsubsection*{Case 9:  $\qOrder-|\bnu|-1 = \pOrder-\mu$, $\qOrder-|\bnu|$ even; $y$ is interior}

In this case, the leading bias will involve two terms:
\begin{align*}
h^*_{\pOrder,\qOrder,\mu, \bnu}(\eval) &=
\argmin_{h>0}\Big[ \frac{1}{nh^{d + 2|\bnu|+2\mu-1}}V_\order(\eval)
+ h^{\pOrder+\qOrder + 3-\mu-|\bnu|}
\left(B_{(i),\qOrder+1}(y,\bx) +B_{(ii),\pOrder+2}(y,\bx) \right)^{2}\Big]
\\ &
= \Big[ \frac{(d + 2|\bnu|+2\mu-1)V_\order(\eval)}
{(\pOrder+\qOrder + 3-\mu-|\bnu|)\left(B_{(i),\qOrder+1}(y,\bx) +B_{(ii),\pOrder+2}(y,\bx) \right)^{2}}\frac{1}{n} \Big]^{\frac{1}{d+ \pOrder + \qOrder + \mu + |\bnu|+2 }}.
\end{align*}

\subsubsection*{Case 10:  $\qOrder-|\bnu|-1 > \pOrder-\mu$, $\pOrder-\mu$ even; $y$ is interior}

In this case, the leading bias will involve only one term:
\begin{align*}
h^*_{\pOrder,\qOrder,\mu, \bnu}(\eval) &=
\argmin_{h>0}\Big[ \frac{1}{nh^{d + 2|\bnu|+2\mu-1}}V_\order(\eval)
+ h^{2\pOrder + 4 -2\mu}
B_{(ii),\pOrder+2}(y,\bx)^{2}\Big]
\\ &
= \Big[ \frac{(d + 2|\bnu|+2\mu-1)V_\order(\eval)}
{(2\pOrder + 4 -2\mu)B_{(ii),\pOrder+2}(y,\bx)^{2}}\frac{1}{n} \Big]^{\frac{1}{d+ 2\pOrder + 2|\bnu| + 3 }}.
\end{align*}

\subsection{Rule-of-thumb bandwidth selection}
This section outlines the methodology
that the companion \texttt{R} package, \texttt{lpcde},
uses to construct the rule-of-thumb bandwidth selection.

The rule-of-thumb estimation uses the following assumptions in order
to compute the optimal bandwidth:
\begin{itemize}
    \item  the data is jointly normal,
    \item  $\bX$ and $\bY$ are independent,
    and,
    \item  $p-\mu = q-|\nu| = 1$.
\end{itemize}
Using these assumptions, each of the terms in the formula
given in Case 1 of Section~\ref{sa-sec:mse_bw} are computed
as follows:
\begin{enumerate}
\item
The densities and relevant derivatives
are evaluated based on the
joint normal distribution assumption.
\item
$\bS_y$, $\bT_y$, $\bT_{\bx}$ and $\bS_{\bx}$
matrices are computed by plugging in for
the range of the data,
the evaluation point,
the respective marginal densities,
and
the kernel used.
\item
Similarly, the $\bc_{y}$ and $\bc_{\bx}$
vectors are computed by using the
range of the data,
the evaluation point,
kernel function,
and
the respective marginal densities.
\item
Bias and variance estimates are constructed using
the relevant entries of the
vectors and matrices.
\end{enumerate}

\section{Alternative variance estimators}
\subsection{V-statistic variance estimator}
We propose here an alternative variance estimator that is quick to implement in practice.
We start by first observing that the estimator
$\hat{\theta}_\order(y|\bx)$
is a V-statistic.
That is,
\begin{align}
\hat{\theta}_\order(y|\bx)
 &= \frac{1}{n^2 h^{1+\mu+d+|\nu|}} \sum_{i,j}
  \Indicator(y_i\leq y_j)
  \be^{\Trans}_{\mu} \hat\bS_{y}^\inv \bP \left( \frac{y_j-y}{h} \right)
  \bQ  ^{\Trans}\left( \frac{\bx_i-\bx}{h} \right)\hat\bS_{\bx}^\inv \be_{\bnu}
  \nonumber\\&
\label{eq:vstat}
  =
  \frac{1}{n^2} \sumIN
  a(y_{i}, y) b(\bx_{i}, \bx)
  +
  \frac{1}{n^2} \sum_{i \neq j}
  \Indicator(y_i\leq y_j)
  a(y_{j}, y) b(\bx_{i}, \bx) ,
\end{align}
where,
\begin{align*}
a(y_{i}, y) &=  h^{-1-\mu}\be^{\Trans}_{\mu} \hat\bS_{y}^\inv \bP \left( \frac{y_i-y}{h} \right),
&
b(\bx_{i}, \bx) &= h^{-d-|\bnu|} \be^{\Trans}_{\bnu}\hat\bS_{\bx}^\inv \bQ  \left( \frac{\bx_i-\bx}{h} \right).
\end{align*}
Note that $a(\cdot)$ and $b(\cdot)$ are scalar functions that
are non-zero only for data points that are
within $h$ distance of the evaluation point.
The second term in \eqref{eq:vstat}
can now be symmetrized and treated as a U-statistic.
Applying the Hoeffding decomposition to the
symmetrized version of the second term and
plugging it back into Equation~\ref{eq:vstat},
we get
\begin{align}
\label{eq:vstat_est}
\hat{\theta}_\order(y|\bx)
&=
\frac{1}{n}\Expectation\left[ a(y_{i}, y)b(\bx_{i}, \bx) \right]
+ \frac{n-1}{n}\Expectation [u_{i,j}]
\nonumber\\&
  + \frac{1}{n^2} \sumIN \left(
    a(y_{i}, y) b(\bx_{i}, \bx)
  -  \Expectation\left[ a(y_{i}, y) b(\bx_{i}, \bx) \right] \right)
  + \frac{n-1}{n} L_\order(\eval)
   +\frac{n-1}{n} W_\order(\eval)
\end{align}
where
\begin{align*}
  u_{i,j} =\frac{1}{2} \left( \Indicator ( y_i \leq y_j) a(y_{j}, y)b(\bx_{i}, \bx)
  + \Indicator(y_j \leq y_i) a(y_{i}, y)b(\bx_{i}, \bx) \right),
\end{align*}
and
\begin{align*}
  L_\order(\eval)
  &= \frac{1}{n} \sumIN 2\left( \Expectation\left[u_{i,j}|y_i, \bx_i\right] -\Expectation \left[ u_{i,j} \right]\right),
  \\
  W_\order(\eval)
  &= {n \choose 2}^{-1} \sum_{i,j=1,i\neq j}^n \left(
      u_{i,j}
      - \Expectation\left[u_{i,j}|y_i, \bx_i\right]
      - \Expectation\left[u_{i,j}|y_j, \bx_j\right]
      + \Expectation \left[ u_{i,j} \right]\right).
\end{align*}
Dependence on polynomial orders is suppressed for notational simplicity.
Since each of the terms in \eqref{eq:vstat}
are orthogonal, the variance of the estimator can be expressed as
the sum of the variance of each of the terms on the right hand side.
Furthermore, we note that the first three terms and $W_\order(\eval)$
have higher-order variance.
Thus, we only need to look at the variance of
$L_\order(\eval)$.
\begin{align*}
\Var [L_\order(\eval)]  &=
\Var \Big[ \frac{2}{n}
\sumIN  \left( \Expectation \left[ u_{i,j} \mid y_i, \bx_i \right] - \Expectation \left[ u_{i,j} \right] \right) \Big]
= \frac{1}{n} \Var \Big[  2\Expectation \left[ u_{i,j} \mid y_i, \bx_i \right] - 2\Expectation \left[ u_{i,j} \right] \Big]
\end{align*}
where we know
\begin{align*}
2\Expectation \left[ u_{i,j} \mid y_i, \bx_i \right] &=
\int_{\mathcal{Y}} \Indicator(y_i \leq u) a(u, y) \diff F(u|\bx_i) b(\bx_{i}, \bx)
+ F(y_i|\bx_i) a(y_{i}, y) b(\bx_{i}, \bx).
\end{align*}
We can expand and simplify this to get
\begin{align*}
\Var \left[L_\order(\eval) \right] &=
\Expectation \left[ \left(
    \int_{\mathcal{Y}} \Indicator(y_i \leq u) a(u, y) \diff F(u|\bx_i) b(\bx_{i}, \bx)
+ F(y_i|\bx_i) a(y_{i}, y) b(\bx_{i}, \bx)
    \right)^2 \right]
     \\ &=
\Expectation \left[ \iint_{\mathcal{Y}} \Indicator(y_i \leq \min\{u, v\}) a(u, y) a(v, y)\diff F(u|\bx_i) \diff F(v|\bx_i) b^{2}(\bx_{i}, \bx)
  \right.\\&\left. \qquad
  + \int_{\mathcal{Y}} \Indicator(y_i \leq u) a(u, y) \diff F(u|\bx_i) F(y_i|\bx_i) a(y_{i}, y) b^{2}(\bx_{i}, \bx)
  \right.\\&\left.\qquad
   +  \left( F(y_i|\bx_i) a(y_{i}, y) b(\bx_{i}, \bx) \right)^2  \right].
\end{align*}
Note that this expression is identical to the
variance expression derived in the proof of
Lemma \ref{sa-lem:Variance}.
This leads to a natural alternative jackknife covariance estimator,
\begin{align*}
\hat{\mathsf{C}}_\order (\eval, y', \bx')
= \frac{1}{n} \sumIN
\hat{L}_{(i), \order}(\eval)
\hat{L}_{(i), \order}(y', \bx').
\end{align*}
where
\begin{align*}
\hat{L}_{(i), \order}(\eval)
=  \frac{2}{n-1} \sum_{j \neq i} \left(u_{i,j}
- \hat{\theta}_\order(y|\bx) \right).
\end{align*}

In particular,
note that if the two evaluation points are equivalent,
we return the variance estimator,
\begin{align*}
\hat{\mathsf{C}}_\order (\eval, \eval)
\equiv
\hat{\mathsf{V}}_\order(\eval)
= \frac{1}{n-1} \sumIN \hat{L}^{2}_{(i), \order}(\eval).
\end{align*}

\subsection{Asymptotic variance estimator}
Another alternative variance estimator
is a sample version of the asymptotic variance
derived in Lemma~\ref{sa-lem:Variance}.
That is,
each of the matrices
in the formula are replaced with
sample analogs.
That is,
\\
(i) $\mu=0$:
\begin{align*}
\hat{\mathsf{V}}_{0,\bnu}(\eval)
  =& \frac{1}{nh^{d+2|\bnu|}} \hat{\theta}_{0,\mathbf{0}}(y|\bx)(1-\hat{\theta}_{0,\mathbf{0}}(y|\bx))\Big(\be_{\bnu}^\Trans \hat{\bS}_{\bx}^{-1} \hat{\bT}_{\bx} \hat{\bS}_{\bx}^{-1} \be_{\bnu}\Big)
\end{align*}
(ii) $\mu > 0$:
\begin{align*}
\hat{\mathsf{V}}_{\order} (\eval) &= \frac{1}{nh^{d + 2|\bnu|+2\mu-1}} \hat{\theta}_{1,\mathbf{0}}(y|\bx)\Big(\be_{\mu}^\Trans\hat{\bS}_y^\inv \hat{\bT}_{y}\hat{\bS}_y^\inv \be_{\mu}\Big)\Big(\be_{\bnu}^\Trans \hat{\bS}_{\bx}^\inv  \hat{\bT}_{\bx}\hat{\bS}_{\bx}^\inv \be_{\bnu}\Big).
\end{align*}
The covariance can be estimated using similar idea.

\section{Proofs}
\label{sa-sec:proofs}

\subsection{Proof of Lemma~\ref{sa-lem:Taking hat off}}
\textbf{Part (i).} See the proof of Lemma 1 in the main paper (i.e., Appendix A.3).

\bigskip \noindent \textbf{Part (ii).} Next consider $\be_{\mu}^{\Trans}\bS_y^{-1}\big(\bar{\bR}_{\eval} - \Expectation[\bar{\bR}_{\eval} |\bX]\big)$, which takes the form
\begin{align*}
& \frac{1}{nh^{\mu+|\bnu|}} \sum_{i=1}^{n} \be_{\mu}^{\Trans}\bS_y^{-1}\int_{\frac{\Y-y}{h}} \Big(\Indicator(y_i\leq y+hu) - F(y+hu|\bx_i)\Big)
\bP (u) g(y+hu)\diff u\frac{1}{h^d}\bQ  \Big(\frac{\bx_i-\bx}{h}\Big)^\Trans.
\end{align*}
It is straightforward to see that
\begin{align*}
\Big|\be_{\mu}^{\Trans}\bS_y^{-1}\int_{\frac{\Y-y}{h}} \Big(\Indicator(y_i\leq y+hu) - F(y+hu|\bx_i)\Big)
\bP (u) g(y+hu)\diff u \frac{1}{h^d}\bQ  \Big(\frac{\bx_i-\bx}{h}\Big)^\Trans\Big| \leq C' h^{-d}
\end{align*}
for some $C'$ that holds uniformly for $y\in \mathcal{Y}$ and $\bx\in\mathcal{X}$. We also have the following bound on the variance
\begin{align*}
&\Var\Big[\be_{\mu}^{\Trans}\bS_y^{-1}\int_{\frac{\Y-y}{h}} \Big(\Indicator(y_i\leq y+hu) - F(y+hu|\bx_i)\Big)
\bP (u) g(y+hu)\diff u \frac{1}{h^d}\bQ  \Big(\frac{\bx_i-\bx}{h}\Big)^\Trans \Big]\\
&\leq
C'\begin{cases}
h^{-d}   & \text{if $\mu=0$} \\
h^{-d+1} & \text{if $\mu>0$}.
\end{cases}.
\end{align*}
Consider the first case above ($\mu=0$). By a discretization $\{ (y_\ell,\bx_\ell): 1\leq \ell\leq M_n \}$ of $\mathcal{Y}\times\mathcal{X}$, we have the probabilistic bound due to Bernstein's inequality
\begin{align*}
\Prob\Big[h^{\mu+|\bnu|}\max_{1\leq \ell\leq M_n}\left|\be_{\mu}^{\Trans}\bS_y^{-1}(\bar{\bR}_{\eval} - \Expectation[\bar{\bR}_{\eval} |\bX])\right|>\mathfrak{c}_1 \mathtt{r}\Big] &\leq
2\exp\Big\{ -\frac{1}{2}\frac{\mathfrak{c}_1^2 \log (n)}{C' + \frac{1}{3}\mathfrak{c}_1C'\mathtt{r}} + \log (M_n)\Big\},
\end{align*}
provided that we set $\mathtt{r} = \sqrt{\log (n)/(nh^d)}$.$M_n$ is at most polynomial in $n$, and the error from discretization can be ignored. This concludes the proof for the $\mu = 0$ case.

For $\mu > 0$, we set $\mathtt{r} = \sqrt{\log (n)/(nh^{d-1})}$, and the probabilistic bound takes the form
\begin{align*}
&\ \Prob\Big[h^{\mu+|\bnu|}\max_{1\leq \ell\leq M_n}\left|\be_{\mu}^{\Trans}\bS_y^{-1}(\bar{\bR}_{\eval} - \Expectation[\bar{\bR}_{\eval} |\bX])\right|>\mathfrak{c}_1 \mathtt{r}\Big]\\ &\leq 2\exp\Big\{ -\frac{1}{2}\frac{\mathfrak{c}_1^2n^2 \mathtt{r}^2}{nC'h^{-d+1} + \frac{1}{3}\mathfrak{c}_1C'h^{-d}n\mathtt{r}} + \log (M_n)\Big\}
= 2\exp\Big\{ -\frac{1}{2}\frac{\mathfrak{c}_1^2 \log (n) }{C' + \frac{\mathfrak{c}_1C'}{3}\sqrt{\frac{\log (n)}{nh^{d+1}}}} + \log (M_n)\Big\}.
\end{align*}
This concludes the proof for the second case, where $\mu>0$.

\subsection{Proof of Lemma~\ref{sa-lem:Bias}}

The conditional expectation of $\bar{\bR}_{\eval}$ in $\check{\theta}_{\order}$ is
\begin{align*}
&\ \Expectation\Big[\frac{1}{nh^{\mu+|\bnu|}}\sumIN \Big[\int_{\Y}\Indicator(y_i\leq u)\frac{1}{h}\bP \left(\frac{u-y}{h}\right)\diff G(u) \Big]\frac{1}{h^d}\bQ  \left(\frac{\bx_i-\bx}{h}\right)^\Trans\Big| \bX\Big]\\
&=  \frac{1}{nh^{\mu+|\bnu|}}\sumIN\Big[\int_{\Y}F(u|\bx_i)\frac{1}{h}\bP \left(\frac{u-y}{h}\right)\diff G(u)\Big]\frac{1}{h^d}\bQ  \left(\frac{\bx_i-\bx}{h}\right)^\Trans.
\end{align*}
To proceed, we employ a Taylor expansion of the conditional distribution function to order $s$:
\begin{align*}
F(u|\bx_i) &= \sum_{\ell + |\bmm| \leq s}  \theta_{\ell,\bmm}(y|\bx)\frac{1}{\ell!\bmm!} (u-y)^{\ell} (\bx_i-\bx)^{\bmm} +  o\Big( \sum_{ \ell + |\bmm|= s} |u-y|^{\ell} |\bx_i-\bx|^{\bmm} \Big).
\end{align*}
Then, the conditional expectation can be simplified as
\begin{align*}
&\ \frac{1}{n}\sumIN\Big[\int_{\Y}F(u|\bx_i)\frac{1}{h}\bP \left(\frac{u-y}{h}\right)\diff G(u)\Big]\frac{1}{h^d}\bQ  \left(\frac{\bx_i-\bx}{h}\right)^\Trans\\
&= \sum_{\ell + |\bmm| \leq s}\theta_{\ell,\bmm}(y|\bx)\Big[\int_{\Y}  \frac{1}{\ell!} (u-y)^{\ell} \frac{1}{h}\bP \left(\frac{u-y}{h}\right)\diff G(u)\Big]\Big[\frac{1}{n}\sumIN\frac{1}{h^d}\frac{1}{\bmm!}(\bx_i-\bx)^{\bmm}\bQ  \left(\frac{\bx_i-\bx}{h}\right)^\Trans\Big]\\
&\qquad+ o\Big(\sum_{\ell + |\bmm| = s}\Big[\int_{\Y}  \frac{1}{\ell!} |u-y|^{\ell} \frac{1}{h}\left|\bP \left(\frac{u-y}{h}\right)\right|\diff G(u)\Big]\Big[\frac{1}{n}\sumIN\frac{1}{h^d}\frac{1}{\bmm!}|\bx_i-\bx|^{\bmm}\Big|\bQ  \left(\frac{\bx_i-\bx}{h}\right)\Big|\Big]\Big)\\
&= \sum_{\ell + |\bmm| \leq s}h^{\ell + |\bmm|}\theta_{\ell,\bmm}(y|\bx)\bc_{y,\ell}\hat{\bc}_{\bx,\bmm}^\Trans + \op\left(  h^{s}\right).
\end{align*}
We note that
\begin{align*}
\bS_{y}^{-1} \bc_{y,\ell} = \be_{\ell} \text{ for all $0\leq \ell\leq \pOrder$},
\end{align*}
and
\begin{align*}
\hat{\bS}_{\bx}^{-1} \hat{\bc}_{\bx,\bmm} = \be_{\bmm} \text{ for all $0\leq |\bmm|\leq \qOrder$}.
\end{align*}
Therefore,
\begin{align*}
&\Expectation[\check{\theta}_{\order}(y|\bx)|\bX]\\
 &= \theta_{\order}(y|\bx)
+ h^{\qOrder+1-|\bnu|}\sum_{|\bmm|=\qOrder+1}\theta_{\mu,\bmm}(y|\bx)\hat{\bc}_{\bx,\bmm}^\Trans\hat{\bS}_{\bx}^{-1}\be_{\bnu} + h^{\pOrder+1-\mu}\theta_{\pOrder+1,\bnu}(y|\bx)\bc_{y,\pOrder+1}^\Trans\bS_{y}^{-1}\be_{\mu}\\
&\qquad + \op\left( h^{\qOrder+1-|\bnu|} +  h^{\pOrder+1-\mu} \right).
\end{align*}
By Lemma~\ref{sa-lem:Taking hat off}, the second term on the right-hand side satisfies
\begin{align*}
h^{\qOrder+1-|\bnu|}\sum_{|\bmm|=\qOrder+1}\theta_{\mu,\bmm}(y|\bx)\hat{\bc}_{\bx,\bmm}^\Trans\hat{\bS}_{\bx}^{-1}\be_{\bnu} &= h^{\qOrder+1-|\bnu|}\sum_{|\bmm|=\qOrder+1}\theta_{\mu,\bmm}(y|\bx){\bc}_{\bx,\bmm}^\Trans{\bS}_{\bx}^{-1}\be_{\bnu} + \Op\Big( h^{\qOrder+1-|\bnu|}\sqrt{\frac{\log (n)}{nh^d}} \Big),
\end{align*}
which means we can denote the leading bias as
\begin{align*}
\mathsf{B}_{\order}(\eval) &= h^{\qOrder+1-|\bnu|}\sum_{|\bmm|=\qOrder+1}\theta_{\mu,\bmm}(y|\bx){\bc}_{\bx,\bmm}^\Trans{\bS}_{\bx}^{-1}\be_{\bnu} + h^{\pOrder+1-\mu}\theta_{\pOrder+1,\bnu}(y|\bx)\bc_{y,\pOrder+1}^\Trans\bS_{y}^{-1}\be_{\mu}.
\end{align*}

For the second claim of this lemma, we again consider a Taylor expansion
\begin{align*}
F(y_j|\bx_i) &= \sum_{\ell + |\bmm| \leq s}  \theta_{\ell,\bmm}(y|\bx)\frac{1}{\ell!\bmm!} (y_j-y)^{\ell} (\bx_i-\bx)^{\bmm} +  o\Big( \sum_{ \ell + |\bmm|= s} |y_j-y|^{\ell} |\bx_i-\bx|^{\bmm} \Big).
\end{align*}
Then
\begin{align*}
&\frac{1}{n^2h^{d+1+\mu+|\bnu|}}\sum_{i,j=1}^n\be_{\mu}^\Trans\hat{\bS}_y^\inv\left[ F(y_j|\bx_i)\bP \left(\frac{y_j-y}{h}\right)\right]\bQ \left(\frac{\bx_i-\bx}{h}\right)^\Trans\hat{\bS}_{\bx}^\inv \be_{\bnu}\\
=&\frac{1}{n^2h^{d+1+\mu+|\bnu|}}\sum_{i,j=1}^n\be_{\mu}^\Trans\hat{\bS}_y^\inv\Big[ \sum_{\ell + |\bmm| \leq s}  \theta_{\ell,\bmm}(y|\bx)\frac{1}{\ell!\bmm!} (y_j-y)^{\ell} (\bx_i-\bx)^{\bmm} \bP \left(\frac{y_j-y}{h}\right)\Big]\bQ \left(\frac{\bx_i-\bx}{h}\right)^\Trans\hat{\bS}_{\bx}^\inv \be_{\bnu}\\
&\quad +o\Big(\frac{1}{n^2h^{d+1+\mu+|\bnu|}}\be_{\mu}^\Trans\hat{\bS}_y^\inv \sum_{i,j=1}^n\Big[  \sum_{ \ell + |\bmm|= s} |y_j-y|^{\ell} |\bx_i-\bx|^{\bmm} \Big|\bP \Big(\frac{y_j-y}{h}\Big)\Big|\Big|\bQ \Big(\frac{\bx_i-\bx}{h}\Big)\Big|\Big]\hat{\bS}_{\bx}^\inv \be_{\bnu}\Big)\\
=& \theta_{\mu,\bnu}(y|\bx) + h^{\qOrder+1-|\bnu|}\sum_{|\bmm|=\qOrder+1}\theta_{\mu,\bmm}(y|\bx)\bc_{\bx,\bmm}^\Trans\bS_{\bx}^{-1}\be_{\bnu} + h^{\pOrder+1-\mu}\theta_{\pOrder+1,\bnu}(y|\bx)\bc_{y,\pOrder+1}^\Trans\bS_{y}^{-1}\be_{\mu}\\
&\quad + \op(h^{\qOrder+1-|\bnu|} +h^{\pOrder+1-\mu} ).
\end{align*}

\subsection{Proof of Lemma~\ref{sa-lem:Variance}}

Let $\bc_1 = \bS_y^\inv\be_{\mu}$ and $\bc_2 = \bS_{\bx}^\inv \be_{\bnu}$.
\begin{align*}
&\Var\Big[   \int_{\Y}\Big[\Indicator(y_i\leq u) - F(u|\bx_i)\Big]\bc_1^\Trans \frac{1}{h}\bP \Big(\frac{u-y}{h}\Big) \diff G(u)  \frac{1}{h^d}\bQ  \left(\frac{\bx_i-\bx}{h}\right)^\Trans \bc_2 \Big]\\
&= \Expectation \Big[ \Var\Big[  \int_{\Y}\Big[\Indicator(y_i\leq u) - F(u|\bx_i)\Big]\bc_1^\Trans \frac{1}{h}\bP \left(\frac{u-y}{h}\right) \diff G(u)  \frac{1}{h^d}\bQ  \left(\frac{\bx_i-\bx}{h}\right)^\Trans \bc_2 \Big| \bX \Big] \Big]\\
&=  \Expectation\Big[\iint_{\frac{\Y-y}{h}}\Big( F(y+h(u_1\wedge u_2)|\bx_i) - F(y+h u_1|\bx_i)F(y+h u_2|\bx_i) \Big)\bc_1^\Trans\bP \left(u_1\right)\bc_1^\Trans\bP \left(u_2\right)   \\
\tag{*}&\qquad\qquad\qquad\qquad\qquad\qquad\qquad g(y+hu_1)g(y+hu_2)\diff u_1\diff u_2\left(\bc_2^\Trans\hd \bQ  \left(\frac{\bx_i-\bx}{h}\right) \right)^2\Big]
\end{align*}
We make a further expansion:
\begin{align*}
&\ F(y+h(u_1\wedge u_2)|\bx_i) - F(y+h u_1|\bx_i)F(y+h u_2|\bx_i)\\
=&\ F(y|\bx_i)(1-F(y|\bx_i))
 + h(u_1\wedge u_2)f(y|\bx_i) - h(u_1+u_2)f(y|\bx_i)F(y|\bx_i)
 + O(h^{2}).
\end{align*}
Note that the remainder term, $O(h^2)$, holds uniformly for $y\in\mathcal{Y}$ and $\bx_i\in\mathcal{X}$ since the conditional distribution function is assumed to have bounded second derivative. Therefore,
{\begin{align*}
\text{(*)}
=&\  \Big(\bc_1^\Trans\bc_{y,0}\bc_{y,0}^\Trans \bc_1\Big)\Expectation\Big[ F(y|\bx_i)(1-F(y|\bx_i))\Big(\bc_2^\Trans\hd \bQ  \Big(\frac{\bx_i-\bx}{h}\Big) \Big)^2\Big]
 \\
&+ h \left(\bc_1^\Trans\bT_y \bc_1\right)
\Expectation\Big[ f(y|\bx_i) \Big(\bc_2^\Trans\hd \bQ  \left(\frac{\bx_i-\bx}{h}\right) \Big)^2\Big]\\
& - h \left[\bc_1^\Trans\Big(\bc_{y,1}\bc_{y,0}^\Trans + \bc_{y,0}\bc_{y,1}^\Trans\Big) \bc_1\right]
\Expectation\Big[f(y|\bx_i)F(y|\bx_i)\Big[\bc_2^\Trans\hd \bQ  \left(\frac{\bx_i-\bx}{h}\right) \Big]^2\Big]
 + O\Big(\frac{1}{h^{d-2}}\Big)\\
=&\  \be_{\mu}^\Trans\be_{0}\Expectation\Big[ F(y|\bx_i)(1-F(y|\bx_i))\Big(\bc_2^\Trans\hd \bQ  \left(\frac{\bx_i-\bx}{h}\right) \Big)^2\Big]\\
&+ h \Big(\be_{\mu}^\Trans\bS_y^\inv \bT_y\bS_y^\inv \be_{\mu}\Big)
\Expectation\Big[ f(y|\bx_i) \Big(\bc_2^\Trans\hd \bQ  \left(\frac{\bx_i-\bx}{h}\right) \Big)^2\Big]
+ O\Big(\frac{1}{h^{d-2}}\Big).
\end{align*}}
To conclude the proof, we note that two scenarios can arise: $\mu=0$ and $\mu > 0$. In the second case,
\begin{align*}
\text{(*)} &= \frac{1}{h^{d-1}}\theta_{1,\mathbf{0}}(y|\bx) \left(\be_{\mu}^\Trans\bS_{y}^{-1}\bT_y \bS_{y}^{-1}\be_{\mu}\right)\left(\be_{\bnu}^\Trans \bS_{\bx}^{-1} \bT_{\bx} \bS_{\bx}^{-1} \be_{\bnu}\right) + O\Big(\frac{1}{h^{d-2}}\Big).
\end{align*}
The first case is more involved. If $\theta_{0,\mathbf{0}}(y|\bx)\neq 0,1$, then
\begin{align*}
\text{(*)} &= \frac{1}{h^d} \theta_{0,\mathbf{0}}(y|\bx)(1-\theta_{0,\mathbf{0}}(y|\bx))\left(\be_{\bnu}^\Trans \bS_{\bx}^{-1} \bT_{\bx} \bS_{\bx}^{-1} \be_{\bnu}\right) + \Op\Big( \frac{1}{h^{d-1}} \Big).
\end{align*}
If $\theta_{0,\mathbf{0}}(y|\bx)=0$ or $1$, then a further expansion is needed, which is why an extra $h$ will be present in the leading variance.

\subsection{Proof of Lemma~\ref{sa-lem:Variance estimation}}

\textbf{Consistency of $\check{\mathsf{V}}_{\mu,\bnu}(y,\bx)$}. For the purposes of this proof, let $\bc_1 = \bS_{y}^{-1}\be_{\mu}$, $\hat{\bc}_2 = \hat{\bS}_{\bx}^{-1}\be_{\bnu}$, and $\bc_2 = \bS_{\bx}^{-1}\be_{\bnu}$. To start, consider
\begin{align*}
&\frac{1}{n^2h^{2d+2\mu+2|\bnu|}}\sumIN \Big[ \int_{\mathcal{Y}}\Big(\Indicator\left(y_i\leq u\right) - \hat{F}(u|\bx_i)\Big)\frac{1}{h}\bc_1^\Trans\bP \left( \frac{u-y}{h} \right)\diff G(u) \bc_2^\Trans\bQ  \left(\frac{\bx_i-\bx}{h}\right) \Big]^2\\
\tag{$\text{I}$}=&\ \frac{1}{n^2h^{2d+2\mu+2|\bnu|}}\sumIN \Big[ \int_{\mathcal{Y}}\Big(\Indicator\left(y_i\leq u\right) - F(u|\bx_i)\Big)\frac{1}{h}\bc_1^\Trans\bP \left( \frac{u-y}{h} \right)\diff G(u) \bc_2^\Trans\bQ  \left(\frac{\bx_i-\bx}{h}\right) \Big]^2\\
-&\ \frac{2}{n^2h^{2d+2\mu+2|\bnu|}}\sumIN \left[\bc_2^\Trans\bQ  \left(\frac{\bx_i-\bx}{h}\right)\right]^2 \iint_{\frac{\Y-y}{h}} \Big(\Indicator\left(y_i\leq u_1\right) - F(u_1|\bx_i)\Big)\Big(\hat{F}(u_2|\bx_i)- F(u_2|\bx_i)\Big) \\
\tag{$\text{II}$}&\qquad\qquad\qquad\qquad \bc_1^\Trans\bP \left(u_1\right)\bc_1^\Trans\bP \left(u_2\right)g(y+hu_1)g(y + hu_2)\diff  u_1\diff  u_2\\
\tag{$\text{III}$}+&\ \frac{1}{n^2h^{2d+2\mu+2|\bnu|}}\sumIN \Big[ \int_{\mathcal{Y}}\Big(\hat{F}(u|\bx_i) - F(u|\bx_i)\Big)\frac{1}{h}\bc_1^\Trans\bP \left( \frac{u-y}{h} \right)\diff G(u) \bc_2^\Trans\bQ  \left(\frac{\bx_i-\bx}{h}\right) \Big]^2.
\end{align*}

First consider term (III). With the uniform convergence result for the estimated conditional distribution function, it is clear that
{\begin{align*}
\left|\text{(III)}\right| \precsim_\Prob
\begin{cases}
\mathsf{V}_{0,\bnu}(y,\bx)\left( h^{2\qOrder+2} + {\frac{\log (n)}{nh^{d}}}\right) & \text{if $\mu=0$, and $\theta_{0, \mathbf{0}}\neq 0$ or $1$}\\
\mathsf{V}_{\mu,\bnu}(y,\bx)\left( h^{2\qOrder+1} + {\frac{\log (n)}{nh^{d+1}}}\right) & \text{if $\mu > 0$, \ \ or $\theta_{0, \mathbf{0}}= 0$ or $1$}
\end{cases}.
\end{align*}}

Now we study term (I), which is clearly unbiased for $\mathsf{V}_{\mu,\bnu}(y,\bx)$. Therefore, we compute its variance.
\begin{align*}
&\Var\left[ \text{(I)} \right] = \frac{1}{n^3h^{4d+4\mu+4|\bnu|}}\Var\Big[ \left( \int_{\mathcal{Y}}\Big(\Indicator\left(y_i\leq u\right) - F(u|\bx_i)\Big)\frac{1}{h}\bc_1^\Trans\bP \left( \frac{u-y}{h} \right)\diff G(u) \bc_2^\Trans\bQ  \left(\frac{\bx_i-\bx}{h}\right) \right)^2 \Big]\\
&\leq \frac{1}{n^3h^{4d+4\mu+4|\bnu|}}\Expectation\Big[  \prod_{j=1}^4\Big[\int_{\mathcal{Y}}\Big(\Indicator\left(y_i\leq u_j\right) - F(u_j|\bx_i)\Big)\frac{1}{h}\bc_1^\Trans\bP \left( \frac{u_j-y}{h} \right)\diff G(u_j)\Big] \Big[\bc_2^\Trans\bQ  \left(\frac{\bx_i-\bx}{h}\right) \Big]^4 \Big].
\end{align*}
With iterative expectation (by conditioning on $\bx_i$), the above further reduces to
\begin{align*}
\Var\left[ \text{(I)} \right]
&= \frac{1}{n^3h^{3d+4\mu+4|\bnu|}}\theta_{0, \mathbf{0}}(1-\theta_{0, \mathbf{0}})(1-3\theta_{0, \mathbf{0}}(1-\theta_{0, \mathbf{0}}))\left[\bc_1^\Trans\bc_{y,0}\right]^4 \Expectation\left[\frac{1}{h^d}\left[\bc_2^\Trans\brx\right]^4\right]\\
&+ O\left(\frac{h}{n^3h^{4d+4\mu+4|\bnu|}}\right).
\end{align*}
In other words,
\begin{align*}
\Big|\text{(I)} - \mathsf{V}_{\mu,\bnu}(y,\bx)\Big| &\precsim_\Prob 
\begin{cases}
\mathsf{V}_{0,\bnu}(y,\bx)\sqrt{\frac{1}{nh^{d}}} & \text{if $\mu=0$, and $\theta_{0, \mathbf{0}}\neq 0$ or $1$}\\
\mathsf{V}_{\mu,\bnu}(y,\bx)\sqrt{\frac{1}{nh^{d+1}}} & \text{if $\mu > 0$, \ \ or $\theta_{0, \mathbf{0}}= 0$ or $1$}
\end{cases}.
\end{align*}

Finally, we consider (II). Using the Cauchy-Schwartz inequality, we have
\begin{align*}
\left|\text{(II)}\right|^2 &\leq \left|\text{(I)}\right| \cdot \left|\text{(III)}\right|.
\end{align*}
As a result,
\begin{align*}
\left|\text{(II)}\right| \precsim_\Prob
\begin{cases}
\mathsf{V}_{0,\bnu}(y,\bx)\sqrt{ h^{2\qOrder+2} + {\frac{\log (n)}{nh^{d}}}} & \text{if $\mu=0$, and $\theta_{0, \mathbf{0}}\neq 0$ or $1$}\\
\mathsf{V}_{\mu,\bnu}(y,\bx)\sqrt{ h^{2\qOrder+1} + {\frac{\log (n)}{nh^{d+1}}}} & \text{if $\mu > 0$, \ \ or $\theta_{0, \mathbf{0}}= 0$ or $1$}
\end{cases}.
\end{align*}
To conclude the proof for $\check{\mathsf{V}}_{\mu,\bnu}(y,\bx)$, we note that replacing $\hat{\bc}_2$ by $\bc_2$ only leads to an additional multiplicative factor $1 + \Op(1/\sqrt{nh^d})$. See Lemma~\ref{sa-lem:Taking hat off}.\bigskip

\noindent \textbf{Consistency of $\hat{\mathsf{V}}_{\mu,\bnu}(y,\bx)$}. For the purposes of this proof, let $\hat{\bc}_1 = \hat{\bS}_{y}^{-1}\be_{\mu}$, $\bc_1 = \bS_{y}^{-1}\be_{\mu}$, $\hat{\bc}_2 = \hat{\bS}_{\bx}^{-1}\be_{\bnu}$, and $\bc_2 = \bS_{\bx}^{-1}\be_{\bnu}$. We first consider the following decomposition
\begin{align*}
&\frac{1}{n^2h^{2d+2\mu+2|\bnu|}}\sumIN \Big[ \int_{\mathcal{Y}}\left(\Indicator\left(y_i\leq u\right) - \hat{F}(u|\bx_i)\right)\frac{1}{h}\bc_1^\Trans\bP \left( \frac{u-y}{h} \right)\diff \hat{F}_y(u) \bc_2^\Trans\bQ  \left(\frac{\bx_i-\bx}{h}\right) \Big]^2\\
=&\frac{1}{n^2h^{2d+2\mu+2|\bnu|}}\sumIN \Big[ \frac{1}{n^2}\sum_{j,k=1}^n\left(\Indicator\left(y_i\leq y_j\right) - F(y_j|\bx_i)\right)\left(\Indicator\left(y_i\leq y_k\right) - F(y_k|\bx_i)\right)\\
\tag{I}&\qquad\qquad\frac{1}{h^2}\bc_1^\Trans\bP \left( \frac{y_j-y}{h} \right)\bc_1^\Trans\bP \left( \frac{y_k-y}{h} \right)   \Big]\left[\bc_2^\Trans\bQ  \left(\frac{\bx_i-\bx}{h}\right)\right]^2\\
-&\frac{2}{n^2h^{2d+2\mu+2|\bnu|}}\sumIN \Big[ \frac{1}{n^2}\sum_{j,k=1}^n\left( \hat{F}(y_j|\bx_i) - F(y_j|\bx_i)\right)\left(\Indicator\left(y_i\leq y_k\right) - F(y_k|\bx_i)\right)\\
\tag{II}&\qquad\qquad\frac{1}{h^2}\bc_1^\Trans\bP \left( \frac{y_j-y}{h} \right)\bc_1^\Trans\bP \left( \frac{y_k-y}{h} \right)   \Big]\left[\bc_2^\Trans\bQ  \left(\frac{\bx_i-\bx}{h}\right)\right]^2\\
+&\frac{1}{n^2h^{2d+2\mu+2|\bnu|}}\sumIN \Big[ \frac{1}{n^2}\sum_{j,k=1}^n\left( \hat{F}(y_j|\bx_i) - F(y_j|\bx_i)\right)\left( \hat{F}(y_k|\bx_i) - F(y_k|\bx_i)\right)\\
\tag{III}&\qquad\qquad\frac{1}{h^2}\bc_1^\Trans\bP \left( \frac{y_j-y}{h} \right)\bc_1^\Trans\bP \left( \frac{y_k-y}{h} \right)   \Big]\left[\bc_2^\Trans\bQ  \left(\frac{\bx_i-\bx}{h}\right)\right]^2.
\end{align*}
By the uniform convergence rate of the estimated conditional distribution function, we have
\begin{align*}
\left|\text{(III)}\right| \precsim_\Prob
\begin{cases}
\mathsf{V}_{0,\bnu}(y,\bx)\left( h^{2\qOrder+2} + {\frac{\log (n)}{nh^{d}}}\right) & \text{if $\mu=0$, and $\theta_{0, \mathbf{0}}\neq 0$ or $1$}\\
\mathsf{V}_{\mu,\bnu}(y,\bx)\left( h^{2\qOrder+1} + {\frac{\log (n)}{nh^{d+1}}}\right) & \text{if $\mu > 0$, \ \ or $\theta_{0, \mathbf{0}}= 0$ or $1$}
\end{cases}.
\end{align*}
Next we consider (II). Using the Cauchy-Schwartz inequality, we have
\begin{align*}
\left|\text{(II)}\right|^2 &\leq \left|\text{(I)}\right| \cdot \left|\text{(III)}\right|.
\end{align*}

Finally, consider term (I), which has the expansion
\begin{align*}
\text{(I)} &= \frac{1}{n^4h^{2d+2\mu+2|\bnu|+2}}\sum_{\substack{i,j,k=1\\ \text{distinct}}}^n \Big[\left(\Indicator\left(y_i\leq y_j\right) - F(y_j|\bx_i)\right)\left(\Indicator\left(y_i\leq y_k\right) - F(y_k|\bx_i)\right)\\
\tag{I.1}&\qquad\qquad\bc_1^\Trans\bP \left( \frac{y_j-y}{h} \right)\bc_1^\Trans\bP \left( \frac{y_k-y}{h} \right)   \Big]\left[\bc_2^\Trans\bQ  \left(\frac{\bx_i-\bx}{h}\right)\right]^2\\
&+ \frac{2}{n^4h^{2d+2\mu+2|\bnu|+2}}\sum_{\substack{i,j=1\\ \text{distinct}}}^n \Big[\left(\Indicator\left(y_i\leq y_j\right) - F(y_j|\bx_i)\right)\left(1 - F(y_i|\bx_i)\right)\bc_1^\Trans\bP \left( \frac{y_j-y}{h} \right)\\
\tag{I.2}&\qquad\qquad\bc_1^\Trans\bP \left( \frac{y_i-y}{h} \right)   \Big]\left[\bc_2^\Trans\bQ  \left(\frac{\bx_i-\bx}{h}\right)\right]^2\\
\tag{I.3}&+ \frac{1}{n^4h^{2d+2\mu+2|\bnu|+2}}\sum_{\substack{i,j=1\\ \text{distinct}}}^n \left(\Indicator\left(y_i\leq y_j\right) - F(y_j|\bx_i)\right)^2\left[\bc_1^\Trans\bP \left( \frac{y_j-y}{h} \right)   \right]^2\left[\bc_2^\Trans\bQ  \left(\frac{\bx_i-\bx}{h}\right)\right]^2\\
\tag{I.4}&+ \frac{1}{n^4h^{2d+2\mu+2|\bnu|+2}}\sum_{i=1}^n \left(1 - F(y_i|\bx_i)\right)^2\left[\bc_1^\Trans\bP \left( \frac{y_i-y}{h} \right)\right]^2   \left[\bc_2^\Trans\bQ  \left(\frac{\bx_i-\bx}{h}\right)\right]^2.
\end{align*}
Then,
\begin{align*}
\left|\text{(I.2)}\right| &\leq \frac{2}{n^4h^{2d+2\mu+2|\bnu|+2}}\sum_{\substack{i,j=1\\ \text{distinct}}}^n \left|\bc_1^\Trans\bP \left( \frac{y_j-y}{h} \right)\right|\cdot \left|\bc_1^\Trans\bP \left( \frac{y_i-y}{h} \right)\right|  \left[\bc_2^\Trans\bQ  \left(\frac{\bx_i-\bx}{h}\right)\right]^2\\
&\leq \frac{2}{n^2h^{d+2\mu+2|\bnu|}}\left[ \frac{1}{nh}\sumIN \left|\bc_1^\Trans\bP \left( \frac{y_i-y}{h} \right)\right| \right]\left[ \frac{1}{nh^{d+1}}\sumIN \left|\bc_1^\Trans\bP \left( \frac{y_i-y}{h} \right)\right|  \left[\bc_2^\Trans\bQ  \left(\frac{\bx_i-\bx}{h}\right)\right]^2\right]\\
&\precsim_\Prob  \frac{1}{n^2h^{d+2\mu+2|\bnu|}}\left( 1 + \sqrt{\frac{1}{nh}} \right)\left( 1 + \sqrt{\frac{1}{nh^{d+1}}} \right)
\precsim \frac{1}{n^2h^{d+2\mu+2|\bnu|}}\\
&\precsim \begin{cases}
\mathsf{V}_{0,\bnu}(y,\bx)\frac{1}{n} & \text{if $\mu=0$, and $\theta_{0, \mathbf{0}}\neq 0$ or $1$}\\
\mathsf{V}_{\mu,\bnu}(y,\bx)\frac{1}{nh} & \text{if $\mu > 0$, \ \ or $\theta_{0, \mathbf{0}}= 0$ or $1$}
\end{cases}.
\end{align*}
Using similar techniques, one can show that
\begin{align*}
\left|\text{(I.3)}\right| &\precsim_\Prob \begin{cases}
\mathsf{V}_{0,\bnu}(y,\bx)\frac{1}{nh} & \text{if $\mu=0$, and $\theta_{0, \mathbf{0}}\neq 0$ or $1$}\\
\mathsf{V}_{\mu,\bnu}(y,\bx)\frac{1}{nh^2} & \text{if $\mu > 0$, \ \ or $\theta_{0, \mathbf{0}}= 0$ or $1$}
\end{cases}\\
\left|\text{(I.4)}\right| &\precsim_\Prob \begin{cases}
\mathsf{V}_{0,\bnu}(y,\bx)\frac{1}{n^2h} & \text{if $\mu=0$, and $\theta_{0, \mathbf{0}}\neq 0$ or $1$}\\
\mathsf{V}_{\mu,\bnu}(y,\bx)\frac{1}{n^2h^2} & \text{if $\mu > 0$, \ \ or $\theta_{0, \mathbf{0}}= 0$ or $1$}
\end{cases}.
\end{align*}

To streamline the remaining derivation, define
\begin{align*}
\phi_{j,i} &= \frac{1}{h}\left(\Indicator\left(y_i\leq y_j\right) - F(y_j|\bx_i)\right)\bc_1^\Trans\bP \left( \frac{y_j-y}{h} \right),\quad
\phi_i = \Expectation[\phi_{j,i}|y_i,\bx_i],\quad
\psi_i = \left[\bc_2^\Trans\bQ  \left(\frac{\bx_i-\bx}{h}\right)\right]^2.
\end{align*}
Then
\begin{align*}
&\text{(I.1)} = \frac{1}{n^4h^{2d+2\mu+2|\bnu|}}\sum_{\substack{i,j,k=1\\ \text{distinct}}}^n \phi_{j,i}\phi_{k,i}\psi_i\\
&= \underbrace{\frac{1}{n^4h^{2d+2\mu+2|\bnu|}}\sum_{\substack{i,j,k=1\\ \text{distinct}}}^n (\phi_{j,i} - \phi_i)(\phi_{k,i}-\phi_i)\psi_i}_{\textstyle \text{(I.1.1)}} + \left( 2 + O\left(\frac{1}{n}\right) \right)\underbrace{\frac{1}{n^3h^{2d+2\mu+2|\bnu|}}\sum_{\substack{i,j=1\\ \text{distinct}}}^n  (\phi_{j,i}-\phi_i)\phi_i\psi_i}_{\textstyle \text{(I.1.2)}}\\
&\qquad + \left( 1 + O\left(\frac{1}{n}\right) \right)\underbrace{\frac{1}{n^2h^{2d+2\mu+2|\bnu|}}\sum_{i=1}^n  \phi_i^2\psi_i}_{\textstyle \text{(I.1.3)}}.
\end{align*}

We have studied the term (I.1.3) in the proof for $\check{\mathsf{V}}_{\mu,\bnu}(y,\bx)$. In particular,
\begin{align*}
\Big|\text{(I.1.3)} - \mathsf{V}_{\mu,\bnu}(y,\bx)\Big| \precsim_\Prob \begin{cases}
\mathsf{V}_{0,\bnu}(y,\bx)\sqrt{\frac{1}{nh^{d}}} & \text{if $\mu=0$, and $\theta_{0, \mathbf{0}}\neq 0$ or $1$}\\
\mathsf{V}_{\mu,\bnu}(y,\bx)\sqrt{\frac{1}{nh^{d+1}}} & \text{if $\mu > 0$, \ \ or $\theta_{0, \mathbf{0}}= 0$ or $1$}
\end{cases}.
\end{align*}

Term (I.1.1) is a mean zero third order U-statistic. Consider its variance
\begin{align*}
& \left(\frac{1}{nh^{d+2\mu+2|\bnu|}}\right)^2\Expectation\Big[ \frac{1}{n^6h^{2d}}\sum_{\substack{i,j,k=1\\ \text{distinct}}}^n\sum_{\substack{i',j',k'=1\\ \text{distinct}}}^n (\phi_{j,i} - \phi_i)(\phi_{k,i}-\phi_i)(\phi_{j',i'} - \phi_{i'})(\phi_{k',i'}-\phi_{i'})\psi_i\psi_{i'} \Big].
\end{align*}
The above expectation is non-zero only in three scenarios: $(j=j',k=k',i\neq i')$, $(j=j',k=k',i=i')$ or $(j=i',k=k',i=j')$. Therefore,
\begin{align*}
\left|\text{(I.1.1)}\right| &\precsim_\Prob \left(\frac{1}{nh^{d+2\mu+2|\bnu|}}\right) \left( \frac{1}{nh} + \sqrt{\frac{1}{n^3h^{d+2}}}\right) \precsim \left(\frac{1}{nh^{d+2\mu+2|\bnu|}}\right)\frac{1}{nh} \\
&\precsim \begin{cases}
\mathsf{V}_{0,\bnu}(y,\bx)\frac{1}{nh} & \text{if $\mu=0$, and $\theta_{0, \mathbf{0}}\neq 0$ or $1$}\\
\mathsf{V}_{\mu,\bnu}(y,\bx)\frac{1}{nh^2} & \text{if $\mu > 0$, \ \ or $\theta_{0, \mathbf{0}}= 0$ or $1$}
\end{cases}.
\end{align*}
Finally consider (I.1.2), which has a mean of zero. Its variance is
\begin{align*}
\Var\left[ \text{(I.1.2)} \right] &= \Big(\frac{1}{nh^{d+2\mu+2|\bnu|}}\Big)^2\Expectation\Big[ \frac{1}{n^4h^{2d}}\sum_{\substack{i,j=1\\ \text{distinct}}}^n \sum_{\substack{i',j'=1\\ \text{distinct}}}^n (\phi_{j,i}-\phi_i)(\phi_{j',i'}-\phi_{i'})\phi_i\psi_i\phi_{i'}\psi_{i'} \Big]\\
&=\Big(\frac{1}{nh^{d+2\mu+2|\bnu|}}\Big)^2\Expectation\Big[ \frac{1}{n^4h^{2d}}\sum_{\substack{i,i',j=1\\ \text{distinct}}}^n  (\phi_{j,i}-\phi_i)(\phi_{j,i'}-\phi_{i'})\phi_i\psi_i\phi_{i'}\psi_{i'} \Big]\\
&+\Big(\frac{1}{nh^{d+2\mu+2|\bnu|}}\Big)^2\Expectation\Big[ \frac{1}{n^4h^{2d}}\sum_{\substack{i,j=1\\ \text{distinct}}}^n  (\phi_{j,i}-\phi_i)^2\phi_i^2\psi_i^2 \Big]\\
&+\Expectation\Big[ \frac{1}{n^4h^{2d}}\sum_{\substack{i,j=1\\ \text{distinct}}}^n (\phi_{j,i}-\phi_i)(\phi_{i,j}-\phi_{j})\phi_i\psi_i\phi_{j}\psi_{j} \Big]\\
&\precsim \left(\frac{1}{nh^{d+2\mu+2|\bnu|}}\right)^2\left( \frac{1}{nh} + \frac{1}{n^2h^{d+1}} \right).
\end{align*}
In addition, an extra $h$ factor emerges if $\mu > 0$, or if $\theta_{0, \mathbf{0}}= 0$ or $1$. As a result,
\begin{align*}
\left|\text{(I.1.2)}\right| \precsim_\Prob \begin{cases}
\mathsf{V}_{0,\bnu}(y,\bx)\sqrt{\frac{1}{nh} } & \text{if $\mu=0$, and $\theta_{0, \mathbf{0}}\neq 0$ or $1$}\\
\mathsf{V}_{\mu,\bnu}(y,\bx)\sqrt{\frac{1}{nh^2} } & \text{if $\mu > 0$, \ \ or $\theta_{0, \mathbf{0}}= 0$ or $1$}
\end{cases}.
\end{align*}
To conclude the proof for $\hat{\mathsf{V}}_{\mu,\bnu}(y,\bx)$, we note that replacing $\hat{\bc}_1$ by $\bc_1$ and $\hat{\bc}_2$ by $\bc_2$ only leads to an additional multiplicative factor $1 + \Op(1/\sqrt{nh^d})$. See Lemma~\ref{sa-lem:Taking hat off}.

\subsection{Proof of Theorem~\ref{sa-thm:Asymptotic normality}}

We will write
\begin{align*}
\bar{\SProc}_{\order}(\eval) &= \frac{1}{\sqrt{n}}\sumIN \frac{h^{d+\mu+|\bnu|}\K_{\mu,\bnu,h}^{\circ}\left( y_i,\bx_i; y,\bx \right)}{\sqrt{\Var\left[h^{d+\mu+|\bnu|}\K_{\mu,\bnu,h}^{\circ}\left( y_i,\bx_i; y,\bx \right)\right]}}.
\end{align*}
Define $\bc_1 = \bS_{y}^{-1}\be_{\mu}$ and $\bc_2 = \bS_{\bx}^{-1}\be_{\bnu}$.

To apply the Berry-Esseen theorem, we first compute the third moment
\begin{align*}
&\Expectation\left[ \left|h^{d+\mu+|\bnu|}\K_{\mu,\bnu,h}^{\circ}\left( y_i,\bx_i; y,\bx \right)\right|^3\right]\\
=& \Expectation\Big[\Big(\prod_{j=1}^3 \Big|\int_{\frac{\Y-y}{h}}\Big( \Indicator(y_i\leq y+hu_j) - F(y+h u_j|\bx_i) \Big)\bc_1^\Trans\bP \left(u_j\right) \diff G(u_j)\Big|\Big) \left|\bc_2^\Trans\bQ  \left(\frac{\bx_i-\bx}{h}\right)\right|^3\Big].
\end{align*}
The leading term in the above is simply
\begin{align*}
&\Expectation\Big[\Big(\prod_{j=1}^3 \Big|\int_{\frac{\Y-y}{h}}\Big( \Indicator(y_i\leq y) - F(y|\bx_i) \Big)\bc_1^\Trans\bP \left(u_j\right) \diff G(u_j)\Big|\Big) \Big|\bc_2^\Trans\bQ  \Big(\frac{\bx_i-\bx}{h}\Big)\Big|^3\Big]\\
=&  \left|\bc_1^\Trans\bc_{y,0}\right|^3\Expectation\Big[ \Big|\Big( \Indicator(y_i\leq y) - F(y|\bx_i) \Big)\Big|^3 \Big|\bc_2^\Trans\bQ  \left(\frac{\bx_i-\bx}{h}\right)\Big|^3\Big]\\
=& \left|\bc_1^\Trans\bc_{y,0}\right|^3\Expectation\Big[ \Big(\theta_{0, \mathbf{0}}(y|\bx)(1-\theta_{0, \mathbf{0}}(\eval))(2\theta_{0, \mathbf{0}}(y|\bx)^2-2\theta_{0, \mathbf{0}}(y|\bx)+1)\Big) \left|\bc_2^\Trans\bQ  \left(\frac{\bx_i-\bx}{h}\right)\right|^3\Big] = O(h^d).
\end{align*}
Note that the above will be exactly zero in cases (ii) and (iii) of Lemma~\ref{sa-lem:Variance}.

\subsection{Omitted details of Remark~\ref{remark:pointwise berry-esseen}}

\textbf{Approximation and coverage error of $\check{\SProc}_{\order}^\circ(\eval)$}. To start,
{\small\begin{align*}
&\check{\SProc}_{\order}^\circ(\eval) - \bar{\SProc}_{\order} (\eval) = \frac{1}{nh^{\mu+|\bnu|}\sqrt{\mathsf{V}_{\order}(\eval)} }\sumIN \be_{\mu}^\Trans\bS_y^\inv \Big[\int_\Y \Big[\Indicator(y_i\leq u) - {F}(u|\bx_i) \Big] \frac{1}{h}\bP \Big(\frac{u-y}{h}\Big)\diff G(u)\Big] \\
&\qquad\qquad\qquad\qquad\frac{1}{h^d}\bQ  \left(\frac{\bx_i-\bx}{h}\right)^\Trans \left(\hat{\bS}_{\bx}^\inv - {\bS}_{\bx}^\inv\right) \be_{\bnu}.
\end{align*}}
By allowing the constant $\mathfrak{c}_1$ to take possibly different values in each term, we have
\begin{align*}
&\Prob\Big[ \sup_{y \in \Y, \bx \in \X}\left| \check{\SProc}_{\order}^\circ(\eval) - \bar{\SProc}_{\order} (\eval) \right| \geq \mathfrak{c}_1\frac{\log (n)}{\sqrt{nh^d}} \Big] \\
&\leq \Prob\Big[\sup_{\bx\in \X} \left|\hat{\bS}_\bx-\bS_\bx\right| > \mathfrak{c}_1\sqrt{\frac{\log (n)}{nh^d}}\Big]
+ \Prob \Big[ \sup_{y \in \Y, \bx \in \X} \Big|\frac{\be_\mu^\Trans\bS_y^{-1}\left(\bar{\bR}_{\eval} - \Expectation \left[\bar{\bR}_{\eval}| \bX \right]\right)}{\sqrt{\mathsf{V}_{\order}(\eval)}}\Big|  > \mathfrak{c}_{1} \sqrt{\log (n)} \Big] \leq \mathfrak{c}_2 n^{-\mathfrak{c}_3},
\end{align*}
where the conclusions follow from the uniform rates established in Lemma~\ref{sa-lem:Taking hat off} and the variance calculations in Lemma~\ref{sa-lem:Variance}.
Next, we consider the normal approximation error. Note that
\begin{align*}
\Prob \left[ \bar{\SProc}_{\order} (\eval) \leq u - \mathfrak{c}_1\frac{\log (n)}{\sqrt{nh^d}} \right] - \mathfrak{c}_2 n^{-\mathfrak{c}_3} \leq \Prob \left[ \check{\SProc}_{\order}^\circ (\eval) \leq u \right] \leq \Prob \left[ \bar{\SProc}_{\order} (\eval) \leq u + \mathfrak{c}_1\frac{\log (n)}{\sqrt{nh^d}} \right] + \mathfrak{c}_2 n^{-\mathfrak{c}_3},
\end{align*}
which means
\begin{align*}
\sup_{u \in \mathbb{R}} \left| \Prob \left[ \check{\SProc}_{\order}^\circ (\eval) \leq u \right] - \Phi(u) \right| \precsim  \frac{\log (n)}{\sqrt{nh^d}} + \mathtt{r_{BE}},
\end{align*}
where $\mathtt{r_{BE}}$ is defined in Theorem~\ref{sa-thm:Asymptotic normality}.\bigskip

\noindent\textbf{Approximation and coverage error of $\hat{\SProc}_{\order}^\circ(\eval)$}. To begin with, we decompose the double sum into
\begin{align*}
&\frac{1}{n^2h^{d+1}}\sum_{i,j=1}^n \Big[\Indicator(y_i\leq y_j) - F(y_j|\bx_i) \Big] \bP \Big(\frac{y_j-y}{h}\Big)\bQ  \left(\frac{\bx_i-\bx}{h}\right)^\Trans\\
\tag{I}=& \frac{1}{n^2h^{d+1}}\sum_{i=1}^n \Bigg(\Big[1 - F(y_i|\bx_i) \Big] \bP \Big(\frac{y_i-y}{h}\Big) -\int_{\Y} \Big[\Indicator(y_i\leq u) - F(u|\bx_i) \Big]  \bP \Big(\frac{u-y}{h}\Big) \diff G(u)\Bigg)\bQ  \left(\frac{\bx_i-\bx}{h}\right)^\Trans \\
+& \frac{1}{n^2h^{d+1}}\sum_{\substack{i,j=1\\ \text{distinct}}}^n \Bigg(\Big[\Indicator(y_i\leq y_j) - F(y_j|\bx_i) \Big]  \bP \Big(\frac{y_j-y}{h}\Big)  -\int_{\Y} \Big[\Indicator(y_i\leq u) - F(u|\bx_i) \Big]  \bP \Big(\frac{u-y}{h}\Big) \diff G(u)\Bigg)\\
\tag{II}&\qquad\qquad\qquad\qquad\bQ  \left(\frac{\bx_i-\bx}{h}\right)^\Trans\\
\tag{III}+& \frac{1}{nh^{d+1}}\sumIN  \int_{\Y} \Big[\Indicator(y_i\leq u) - F(u|\bx_i) \Big]  \bP \Big(\frac{u-y}{h}\Big) \diff G(u)\bQ  \left(\frac{\bx_i-\bx}{h}\right)^\Trans,
\end{align*}
where we set $G = F_y$. Term (I) represents the leave-in bias, and it is straightforward to show that
\begin{align*}
\Prob\Big[ \sup_{y \in \Y, \bx \in \X} |\text{(I)}|  > \mathfrak{c}_{1} \frac{1}{n}\Big(1+\sqrt{\frac{\log (n)}{nh^{d+1}}}\Big) \Big] \leq \mathfrak{c}_{2}n^{-\mathfrak{c}_{3}},
\end{align*}
for some constants $\mathfrak{c}_{1}$, $\mathfrak{c}_{2}$, and $\mathfrak{c}_{3}$. See Lemma~\ref{sa-lem:Taking hat off} for the proof strategy.

Term (II) is a degenerate U-statistic. Define
\begin{align*}
u_{i,j} &= \bc_1^\Trans\Bigg(\Big[\Indicator(y_i\leq y_j) - F(y_j|\bx_i) \Big]  \bP \Big(\frac{y_j-y}{h}\Big)  -\int_{\Y} \Big[\Indicator(y_i\leq u) - F(u|\bx_i) \Big]  \bP \Big(\frac{u-y}{h}\Big) \diff G(u)\Bigg)\\
&\qquad\qquad\qquad\bQ  \left(\frac{\bx_i-\bx}{h}\right)^\Trans\bc_2,
\end{align*}
where $\bc_1$ and $\bc_2$ are arbitrary (fixed) vectors of conformable dimensions. Then we apply Equation (3.5) of \cite{Gine-Latala-Zinn_2000_Ustat} (Lemmas 7 and 8 in the main paper), which gives (the value of $C'$ may change for each line)
\begin{align*}
&\Prob\Big[\sup_{y \in \Y, \bx \in \X}\Big| \sum_{\substack{i,j=1,\  i\neq j}}^n u_{i,j} \Big| > t \Big] \leq C'\exp\left\{ -\frac{1}{C'}\min\left[ \frac{t}{\sqrt{n^2h^{d+1}}},\ \frac{t^{2/3}}{(nh)^{1/3}},\ t^{1/2} \right] + \log (n) \right\}\\
&= C'\exp\left\{ -\frac{\sqrt{C}}{C'}\min\left[ \log (n),\ \left((\log (n))^2nh^d\right)^{\frac{1}{3}},\ \left((\log (n))^2n^2h^{d+1}\right)^{\frac{1}{4}} \right] + \log (n) \right\}.
\end{align*}
As a result,
\begin{align*}
\Prob\Big[ \sup_{y \in \Y, \bx \in \X} |\text{(II)}|  > \mathfrak{c}_{1} \frac{\log (n)}{\sqrt{n^2h^{d+1}}} \Big] \leq \mathfrak{c}_{2}n^{-\mathfrak{c}_{3}},
\end{align*}
for some constants $\mathfrak{c}_{1}$, $\mathfrak{c}_{2}$, and $\mathfrak{c}_{3}$.

We now collect the pieces. The difference between $\hat{\SProc}_{\order}^\circ(\eval)$ and $\check{\SProc}_{\order}^\circ(\eval)$ is
\begin{align*}
&\hat{\SProc}_{\order}^\circ(\eval) - \check{\SProc}_{\order}^\circ(\eval)\\
=& \frac{1}{h^{\mu+|\bnu|}\sqrt{\mathsf{V}_{\order}(\eval)} }
\be_{\mu}^\Trans \hat{\bS}_y^\inv
\Big[ \text{(I)} + \text{(II)} \Big] \hat{\bS}_{\bx}^\inv  \be_{\bnu} + \frac{1}{\sqrt{\mathsf{V}_{\order}(\eval)} }
\be_{\mu}^\Trans \left(\hat{\bS}_y^\inv - {\bS}_y^\inv\right)
\Big[ \bar{\bR}_{\eval} - \Expectation \left[\bar{\bR}_{\eval}| \bX \right] \Big] \hat{\bS}_{\bx}^\inv  \be_{\bnu},
\end{align*}
and the conclusion follows from Lemmas \ref{sa-lem:Taking hat off} and \ref{sa-lem:Variance}.

\subsection{Omitted details of Remark \ref{SA-remark: local sample size}}

To show this result, we first partition the support $\mathcal{Y} \times \mathcal{X}$ into cubes with edge length $\mathfrak{c}_3h$, where the constant $\mathfrak{c}_3$ is chosen so that, for any $(y,\bx)$ in $\mathcal{Y} \times \mathcal{X}$, at least one of the cubes will be contained in the ball $\{ y':  |y'-y|\leq \mathfrak{c}_1h\}\times \{ \bx':  |\bx'-\bx|\leq \mathfrak{c}_1h\}$. The number of cubes in this partition is $\lceil1/(\mathfrak{c}_3h)^{d+1}\rceil$. Then the conclusion follows from Lemma \ref{SA-lem: smallest cell}. \bigskip

\noindent \textbf{Proof of Lemma \ref{SA-lem: smallest cell}}. For simplicity let $c_n = (1-\delta_n) \frac{n}{J_n}$. We first employ the union bound
\begin{align*}
\Prob\Big[\min_{1\leq j\leq J_n}z_j < c_n\Big] &\leq J_n\cdot \Prob\left[z_j < c_n\right].
\end{align*}
Note that $z_j\sim \mathrm{Binomial}(n;\frac{1}{J_n})$, and therefore
\begin{align*}
\Prob\left[z_j  < c_n\right] = \Prob\Big[z_j - \frac{n}{J_n}  < c_n - \frac{n}{J_n}\Big]
&\leq
\exp\Big( -\frac{1}{2}\frac{\left( \frac{n}{J_n} - c_n \right)^2}{\frac{n}{J_n}\left(1-\frac{1}{J_n}\right)+\frac{1}{3}\left( \frac{n}{J_n} - c_n \right)} \Big)\\
&\leq \exp\Big( -\frac{3}{8}\frac{J_n}{n}\Big( \frac{n}{J_n} - c_n \Big)^2 \Big).
\end{align*}
Then we have
\begin{align*}
\Prob\Big[\min_{1\leq j\leq J_n}z_j < c_n\Big] &\leq J_n\exp\left( -\frac{3\delta_n^2}{8}\frac{n}{J_n} \right)= \frac{1}{\pi_n\log n }\exp\left( -\left(\frac{3}{8}\delta_n^2\pi_n-1\right)\log n \right).
\end{align*}
Therefore, the above will vanish faster than any polynomial of $n$ provided that $\delta_n^2\pi_n\to\infty$.

\subsection{Proof of Lemma~\ref{sa-lem:Uniform rate of convergence}}

\textbf{Part (i) Convergence of $\check{\theta}_{\mu,\bnu} - \theta_{\mu,\bnu}$.} Recall that we have the following decomposition of our estimator
\begin{align*}
&\check{\theta}_\order  - \theta_\order\\
\tag{I}&= \frac{1}{nh^{1+d+\mu+|\bnu|}}\sumIN \left[  \int_{\Y}F(u|\bx_i)\be_{\mu}^\Trans\bS_y^\inv \bP \left(\frac{u-y}{h}\right) \diff G(u)\right] \bQ  \left(\frac{\bx_i-\bx}{h}\right)^\Trans \hat{\bS}_{\bx}^\inv \be_{\bnu} - \theta_\order\\
\tag{II}&+ \frac{1}{nh^{1+d+\mu+|\bnu|}}\sumIN \left[  \int_{\Y}\Big(\Indicator(y_i\leq u) - F(u|\bx_i)\Big)\be_{\mu}^\Trans\bS_y^\inv \bP \left(\frac{u-y}{h}\right) \diff G(u)\right] \bQ  \left(\frac{\bx_i-\bx}{h}\right)^\Trans \bS_{\bx}^\inv \be_{\bnu}\\
\tag{III}&+ \frac{1}{nh^{1+d+\mu+|\bnu|}}\sumIN \left[  \int_{\Y}\Big(\Indicator(y_i\leq u) - F(u|\bx_i)\Big)\be_{\mu}^\Trans\bS_y^\inv \bP \left(\frac{u-y}{h}\right) \diff G(u)\right] \bQ  \left(\frac{\bx_i-\bx}{h}\right)^\Trans \Big(\hat{\bS}_{\bx}^\inv - \bS_{\bx}^\inv\Big) \be_{\bnu}.
\end{align*}
(I) is simply the conditional bias, whose order is given in Lemma~\ref{sa-lem:Bias}. The convergence rate of (II) can be easily deduced from that of $\be_{\mu}^{\Trans}\bS_y^{-1}\left(\bar{\bR}_{\eval} - \Expectation \left[\bar{\bR}_{\eval}|\bX\right]\right)$ in Lemma~\ref{sa-lem:Taking hat off}. Finally, it should be clear that (III) is negligible relative to (II).

\bigskip\noindent
\textbf{Part (ii) Convergence of $\hat{\theta}_{\mu,\bnu} - \theta_{\mu,\bnu}$.} This part follows from Remark~\ref{remark:pointwise berry-esseen}.

\subsection{Proof of Lemma~\ref{sa-lem:Uniform variance estimation}}

\textbf{Uniform consistency of $\check{\mathsf{V}}_{\mu,\bnu}(y,\bx)$}. For the purposes of this proof, let $\bc_1 = \bS_{y}^{-1}\be_{\mu}$, $\hat{\bc}_2 = \hat{\bS}_{\bx}^{-1}\be_{\bnu}$, and $\bc_2 = \bS_{\bx}^{-1}\be_{\bnu}$. To start, consider
\begin{align*}
&\frac{1}{n^2h^{2d+2\mu+2|\bnu|}}\sumIN \left[ \int_{\mathcal{Y}}\Big(\Indicator\left(y_i\leq u\right) - \hat{F}(u|\bx_i)\Big)\frac{1}{h}\bc_1^\Trans\bP \left( \frac{u-y}{h} \right)\diff G(u) \bc_2^\Trans\bQ  \left(\frac{\bx_i-\bx}{h}\right) \right]^2\\
\tag{$\text{I}$}=&\ \frac{1}{n^2h^{2d+2\mu+2|\bnu|}}\sumIN \left[ \int_{\mathcal{Y}}\Big(\Indicator\left(y_i\leq u\right) - F(u|\bx_i)\Big)\frac{1}{h}\bc_1^\Trans\bP \left( \frac{u-y}{h} \right)\diff G(u) \bc_2^\Trans\bQ  \left(\frac{\bx_i-\bx}{h}\right) \right]^2\\
-&\ \frac{2}{n^2h^{2d+2\mu+2|\bnu|}}\sumIN \iint_{\frac{\Y-y}{h}} \Big(\Indicator\left(y_i\leq u_1\right) - F(u_1|\bx_i)\Big)\Big(\hat{F}(u_2|\bx_i)- F(u_2|\bx_i)\Big)\bc_1^\Trans\bP \left(u_1\right)\bc_1^\Trans\bP \left(u_2\right) \\
\tag{$\text{II}$}&\qquad\qquad\qquad\qquad g(y+hu_1)g(y + hu_2)\diff  u_1\diff  u_2\left[\bc_2^\Trans\bQ  \left(\frac{\bx_i-\bx}{h}\right)\right]^2\\
\tag{$\text{III}$}+&\ \frac{1}{n^2h^{2d+2\mu+2|\bnu|}}\sumIN \left[ \int_{\mathcal{Y}}\Big(\hat{F}(u|\bx_i) - F(u|\bx_i)\Big)\frac{1}{h}\bc_1^\Trans\bP \left( \frac{u-y}{h} \right)\diff G(u) \bc_2^\Trans\bQ  \left(\frac{\bx_i-\bx}{h}\right) \right]^2.
\end{align*}

First consider term (I). Clearly this term is unbiased for $\mathsf{V}_{\mu,\bnu}(y,\bx)$. In the proof of Lemma~\ref{sa-lem:Variance estimation}, we showed that
\begin{align*}
\Var\left[ \left( \int_{\mathcal{Y}}\Big(\Indicator\left(y_i\leq u\right) - F(u|\bx_i)\Big)\frac{1}{h}\bc_1^\Trans\bP \left( \frac{u-y}{h} \right)\diff G(u) \bc_2^\Trans\bQ  \left(\frac{\bx_i-\bx}{h}\right) \right)^2 \right]
&\leq C_1 \begin{cases}
h^d & \text{if }\mu=0\\
h^{d+1} & \text{if }\mu>0
\end{cases}.
\end{align*}
Also note that
\begin{align*}
\left( \int_{\mathcal{Y}}\Big(\Indicator\left(y_i\leq u\right) - F(u|\bx_i)\Big)\frac{1}{h}\bc_1^\Trans\bP \left( \frac{u-y}{h} \right)\diff G(u) \bc_2^\Trans\bQ  \left(\frac{\bx_i-\bx}{h}\right) \right)^2 \leq C_2.
\end{align*}
In the above, the constants, $C_1$ and $C_2$, can be chosen to be independent of the evaluation point, the sample size, and the bandwidth. Then by a proper discretization of $\mathcal{Y}\times \mathcal{X}$, and applying the union bound and Bernstein's inequality, one has
\begin{align*}
\Prob\Big[\sup_{y\in\mathcal{Y},\bx\in\mathcal{X}}\Big| \frac{\text{(I)} -\mathsf{V}_{\mu,\bnu}(y,\bx) }{\mathsf{V}_{\mu,\bnu}(y,\bx)} \Big| > \mathfrak{c}_1\mathtt{r_1}\Big] \leq \mathfrak{c}_2n^{-\mathfrak{c}_3},\qquad \mathtt{r_1}=\begin{cases}
\sqrt{\frac{\log (n)}{nh^{d}}} & \text{if }\mu=0\\
\sqrt{\frac{\log (n)}{nh^{d+1}}} & \text{if }\mu>0
\end{cases},
\end{align*}
for some constants $\mathfrak{c}_1$, $\mathfrak{c}_2$, and $\mathfrak{c}_3$. In addition, $\mathfrak{c}_3$ can be made arbitrarily large by appropriate choices of $\mathfrak{c}_1$. See the proof of Lemma~\ref{sa-lem:Taking hat off} for an example of this proof strategy.

Next consider term (III). With the uniform convergence result for the estimated conditional distribution function, it is clear that
\begin{align*}
\Prob\Big[\sup_{y\in\mathcal{Y},\bx\in\mathcal{X}}\Big| \frac{\text{(III)} }{\mathsf{V}_{\mu,\bnu}(y,\bx)} \Big| > \mathfrak{c}_1\mathtt{r_3}\Big] \leq \mathfrak{c}_2n^{-\mathfrak{c}_3},\qquad \mathtt{r_3}= h^{2\qOrder+1} + {\frac{\log (n)}{nh^{d+1}}}.
\end{align*}

Finally, we consider (II). Using the Cauchy-Schwartz inequality, we have
\begin{align*}
\left|\text{(II)}\right|^2 &\leq \left|\text{(I)}\right| \cdot \left|\text{(III)}\right|.
\end{align*}
As a result,
\begin{align*}
\Prob\Big[\sup_{y\in\mathcal{Y},\bx\in\mathcal{X}}\Big| \frac{\text{(II)} }{\mathsf{V}_{\mu,\bnu}(y,\bx)} \Big| > \mathfrak{c}_1\mathtt{r_2}\Big] \leq \mathfrak{c}_2n^{-\mathfrak{c}_3},\qquad \mathtt{r_2}=\sqrt{1 + \mathtt{r_1}}\sqrt{\mathtt{r_3}}.
\end{align*}
To conclude the proof for $\check{\mathsf{V}}_{\mu,\bnu}(y,\bx)$, we note that replacing $\hat{\bc}_2$ by $\bc_2$ only leads to an additional multiplicative factor $1 + \Op(\sqrt{\log (n)/(nh^d)})$. See Lemma~\ref{sa-lem:Taking hat off}.\bigskip

\noindent\textbf{Uniform consistency of $\hat{\mathsf{V}}_{\mu,\bnu}(y,\bx)$}. For the purposes of this proof, let $\hat{\bc}_1 = \hat{\bS}_{y}^{-1}\be_{\mu}$, $\bc_1 = \bS_{y}^{-1}\be_{\mu}$, $\hat{\bc}_2 = \hat{\bS}_{\bx}^{-1}\be_{\bnu}$, and $\bc_2 = \bS_{\bx}^{-1}\be_{\bnu}$. We consider the same decomposition used in the proof of Lemma~\ref{sa-lem:Variance estimation}:
\begin{align*}
&\frac{1}{n^2h^{2d+2\mu+2|\bnu|}}\sumIN \Big[ \int_{\mathcal{Y}}\left(\Indicator\left(y_i\leq u\right) - \hat{F}(u|\bx_i)\right)\frac{1}{h}\bc_1^\Trans\bP \left( \frac{u-y}{h} \right)\diff \hat{F}_y(u) \bc_2^\Trans\bQ  \left(\frac{\bx_i-\bx}{h}\right) \Big]^2\\
=&\frac{1}{n^2h^{2d+2\mu+2|\bnu|}}\sumIN \Big[ \frac{1}{n^2}\sum_{j,k=1}^n\left(\Indicator\left(y_i\leq y_j\right) - F(y_j|\bx_i)\right)\left(\Indicator\left(y_i\leq y_k\right) - F(y_k|\bx_i)\right)\\
\tag{I}&\qquad\qquad\frac{1}{h^2}\bc_1^\Trans\bP \left( \frac{y_j-y}{h} \right)\bc_1^\Trans\bP \left( \frac{y_k-y}{h} \right)   \Big]\left[\bc_2^\Trans\bQ  \left(\frac{\bx_i-\bx}{h}\right)\right]^2\\
-&\frac{2}{n^2h^{2d+2\mu+2|\bnu|}}\sumIN \Big[ \frac{1}{n^2}\sum_{j,k=1}^n\left( \hat{F}(y_j|\bx_i) - F(y_j|\bx_i)\right)\left(\Indicator\left(y_i\leq y_k\right) - F(y_k|\bx_i)\right)\\
\tag{II}&\qquad\qquad\frac{1}{h^2}\bc_1^\Trans\bP \left( \frac{y_j-y}{h} \right)\bc_1^\Trans\bP \left( \frac{y_k-y}{h} \right)   \Big]\left[\bc_2^\Trans\bQ  \left(\frac{\bx_i-\bx}{h}\right)\right]^2\\
+&\frac{1}{n^2h^{2d+2\mu+2|\bnu|}}\sumIN \Big[ \frac{1}{n^2}\sum_{j,k=1}^n\left( \hat{F}(y_j|\bx_i) - F(y_j|\bx_i)\right)\left( \hat{F}(y_k|\bx_i) - F(y_k|\bx_i)\right)\\
\tag{III}&\qquad\qquad\frac{1}{h^2}\bc_1^\Trans\bP \left( \frac{y_j-y}{h} \right)\bc_1^\Trans\bP \left( \frac{y_k-y}{h} \right)   \Big]\left[\bc_2^\Trans\bQ  \left(\frac{\bx_i-\bx}{h}\right)\right]^2.
\end{align*}
By the uniform convergence rate for the estimated conditional distribution function, we have
\begin{align*}
\Prob\Big[\sup_{y\in\mathcal{Y},\bx\in\mathcal{X}}\Big| \frac{\text{(III)} }{\mathsf{V}_{\mu,\bnu}(y,\bx)} \Big| > \mathfrak{c}_1\mathtt{r_3}\Big] \leq \mathfrak{c}_2n^{-\mathfrak{c}_3},\qquad \mathtt{r_3}= h^{2\qOrder+1} + {\frac{\log (n)}{nh^{d+1}}}.
\end{align*}

Employing the Cauchy-Schwartz inequality gives
\begin{align*}
\left|\text{(II)}\right|^2 &\leq \left|\text{(I)}\right| \cdot \left|\text{(III)}\right|.
\end{align*}
As a result, a probabilistic order for term (II) follows that of terms (I) and (III).

Finally, consider term (I), which has the expansion
\begin{align*}
\text{(I)}&= \frac{1}{n^4h^{2d+2\mu+2|\bnu|+2}}\sum_{\substack{i,j,k=1\\ \text{distinct}}}^n \Big[\left(\Indicator\left(y_i\leq y_j\right) - F(y_j|\bx_i)\right)\left(\Indicator\left(y_i\leq y_k\right) - F(y_k|\bx_i)\right)\\
\tag{I.1}&\qquad\qquad\bc_1^\Trans\bP \left( \frac{y_j-y}{h} \right)\bc_1^\Trans\bP \left( \frac{y_k-y}{h} \right)   \Big]\Big[\bc_2^\Trans\bQ  \left(\frac{\bx_i-\bx}{h}\right)\Big]^2\\
&+ \frac{2}{n^4h^{2d+2\mu+2|\bnu|+2}}\sum_{\substack{i,j=1\\ \text{distinct}}}^n \Big[\left(\Indicator\left(y_i\leq y_j\right) - F(y_j|\bx_i)\right)\left(1 - F(y_i|\bx_i)\right)\\
\tag{I.2}&\qquad\qquad\bc_1^\Trans\bP \left( \frac{y_j-y}{h} \right)\bc_1^\Trans\bP \left( \frac{y_i-y}{h} \right)   \Big]\Big[\bc_2^\Trans\bQ  \left(\frac{\bx_i-\bx}{h}\right)\Big]^2\\
\tag{I.3}&+ \frac{1}{n^4h^{2d+2\mu+2|\bnu|+2}}\sum_{\substack{i,j=1\\ \text{distinct}}}^n \left(\Indicator\left(y_i\leq y_j\right) - F(y_j|\bx_i)\right)^2\left[\bc_1^\Trans\bP \left( \frac{y_j-y}{h} \right)   \right]^2\left[\bc_2^\Trans\bQ  \left(\frac{\bx_i-\bx}{h}\right)\right]^2\\
\tag{I.4}&+ \frac{1}{n^4h^{2d+2\mu+2|\bnu|+2}}\sum_{i=1}^n \left(1 - F(y_i|\bx_i)\right)^2\left[\bc_1^\Trans\bP \left( \frac{y_i-y}{h} \right)\right]^2   \left[\bc_2^\Trans\bQ  \left(\frac{\bx_i-\bx}{h}\right)\right]^2.
\end{align*}
Then,
\begin{align*}
\left|\text{(I.2)}\right| &\leq \frac{1}{n^4h^{2d+2\mu+2|\bnu|+2}}\sum_{\substack{i,j=1\\ \text{distinct}}}^n \left|\bc_1^\Trans\bP \left( \frac{y_j-y}{h} \right)\right|\cdot \left|\bc_1^\Trans\bP \left( \frac{y_i-y}{h} \right)\right|  \left[\bc_2^\Trans\bQ  \left(\frac{\bx_i-\bx}{h}\right)\right]^2\\
&\leq \left(\frac{1}{nh^{d+2\mu+2|\bnu|}}\right)\frac{2}{n}\left[ \frac{1}{nh}\sumIN \left|\bc_1^\Trans\bP \left( \frac{y_i-y}{h} \right)\right| \right]\left[ \frac{1}{nh^{d+1}}\sumIN \left|\bc_1^\Trans\bP \left( \frac{y_i-y}{h} \right)\right|  \left[\bc_2^\Trans\bQ  \left(\frac{\bx_i-\bx}{h}\right)\right]^2\right],
\end{align*}
which means
\begin{align*}
\Prob\Big[\sup_{y \in \Y, \bx \in \X} \Big|\frac{\text{(I.2)}}{\mathsf{V}_{\mu,\bnu}(y,\bx)}\Big| \geq \mathfrak{c}_1\frac{1}{nh} \Big] &\leq \mathfrak{c}_2n^{-\mathfrak{c}_3}.
\end{align*}

Using similar techniques, one can show that
\begin{align*}
\Prob\Big[\sup_{y \in \Y, \bx \in \X} \Big|\frac{\text{(I.3)}}{\mathsf{V}_{\mu,\bnu}(y,\bx)}\Big| \geq \mathfrak{c}_1\frac{1}{nh^2} \Big] &\leq \mathfrak{c}_2n^{-\mathfrak{c}_3},\quad \Prob\Big[\sup_{y \in \Y, \bx \in \X} \Big|\frac{\text{(I.4)}}{\mathsf{V}_{\mu,\bnu}(y,\bx)}\Big| \geq \mathfrak{c}_1\frac{1}{n^2h^2} \Big] \leq \mathfrak{c}_2n^{-\mathfrak{c}_3}.
\end{align*}

To streamline the remaining derivation, define
\begin{align*}
\phi_{j,i} &= \frac{1}{h}\left(\Indicator\left(y_i\leq y_j\right) - F(y_j|\bx_i)\right)\bc_1^\Trans\bP \left( \frac{y_j-y}{h} \right),\quad
\phi_i = \Expectation[\phi_{j,i}|y_i,\bx_i],\quad
\psi_i = \left[\bc_2^\Trans\bQ  \left(\frac{\bx_i-\bx}{h}\right)\right]^2.
\end{align*}
Then
\begin{align*}
&\text{(I.1)} = \frac{1}{n^4h^{2d+2\mu+2|\bnu|}}\sum_{\substack{i,j,k=1\\ \text{distinct}}}^n \phi_{j,i}\phi_{k,i}\psi_i\\
&= \underbrace{\frac{1}{n^4h^{2d+2\mu+2|\bnu|}}\sum_{\substack{i,j,k=1\\ \text{distinct}}}^n (\phi_{j,i} - \phi_i)(\phi_{k,i}-\phi_i)\psi_i}_{\textstyle \text{(I.1.1)}} + \left( 2 + O\left(\frac{1}{n}\right) \right)\underbrace{\frac{1}{n^3h^{2d+2\mu+2|\bnu|}}\sum_{\substack{i,j=1\\ \text{distinct}}}^n  (\phi_{j,i}-\phi_i)\phi_i\psi_i}_{\textstyle \text{(I.1.2)}} \\
&\qquad+ \left( 1 + O\left(\frac{1}{n}\right) \right)\underbrace{\frac{1}{n^2h^{2d+2\mu+2|\bnu|}}\sum_{i=1}^n  \phi_i^2\psi_i}_{\textstyle \text{(I.1.3)}}.
\end{align*}
By employing the same techniques in the proof for $\check{\mathsf{V}}_{\mu,\bnu}(y,\bx)$, we have that
\begin{align*}
\Prob\Big[\sup_{y\in\mathcal{Y},\bx\in\mathcal{X}}\Big| \frac{\text{(I.1.3)} -\mathsf{V}_{\mu,\bnu}(y,\bx) }{\mathsf{V}_{\mu,\bnu}(y,\bx)} \Big| > \mathfrak{c}_1\mathtt{r_1}\Big] \leq \mathfrak{c}_2n^{-\mathfrak{c}_3},\qquad \mathtt{r_1}=\begin{cases}
\sqrt{\frac{\log (n)}{nh^{d}}} & \text{if }\mu=0\\
\sqrt{\frac{\log (n)}{nh^{d+1}}} & \text{if }\mu>0
\end{cases}.
\end{align*}

Term (I.1.2) admits the following decomposition:
\begin{align*}
\text{(I.1.2)} &= \underbrace{\frac{1}{nh^{d+2\mu+2|\bnu|}}\frac{n-1}{n^2}\sum_{j=1}^n  \Expectation\Big[\frac{1}{h^d}(\phi_{j,i}-\phi_i)\phi_i\psi_i\Big| y_j,\bx_j\Big]}_{\textstyle \text{(I.1.2.1)}} \\
&\qquad + \underbrace{\frac{1}{nh^{d+2\mu+2|\bnu|}}\frac{1}{n^2h^{d}}\sum_{\substack{i,j=1\\ \text{distinct}}}^n  \left((\phi_{j,i}-\phi_i)\phi_i\psi_i - \Expectation\Big[\frac{1}{h^d}(\phi_{j,i}-\phi_i)\phi_i\psi_i\Big| y_j,\bx_j\Big]\right)}_{\textstyle \text{(I.1.2.2)}}.
\end{align*}
Using the same techniques of Lemmas \ref{sa-lem:Taking hat off} and \ref{sa-lem:Variance estimation}, we have
\begin{align*}
&\mu = 0\qquad \Prob\Big[\sup_{y\in\mathcal{Y},\bx\in\mathcal{X}}\Big| \frac{\text{(I.1.2.1)}  }{\mathsf{V}_{\mu,\bnu}(y,\bx)} \Big| > \mathfrak{c}_1\sqrt{\frac{\log (n)}{nh}}\Big] \leq \mathfrak{c}_2n^{-\mathfrak{c}_3},\\
&\mu > 0\qquad \Prob\Big[\sup_{y\in\mathcal{Y},\bx\in\mathcal{X}}\Big| \frac{\text{(I.1.2.1)}  }{\mathsf{V}_{\mu,\bnu}(y,\bx)} \Big| > \mathfrak{c}_1\sqrt{\frac{\log (n)}{nh^2}}\Big] \leq \mathfrak{c}_2n^{-\mathfrak{c}_3}.
\end{align*}
Term (I.1.2.2) is a degenerate second order U-statistic. We adopt Equation (3.5) of \cite{Gine-Latala-Zinn_2000_Ustat} (Lemmas 7 and 8 in the main paper), which implies (see Remark~\ref{remark:pointwise berry-esseen} and its proof for an example)
\begin{align*}
\Prob\Big[\sup_{y\in\mathcal{Y},\bx\in\mathcal{X}}\Big| \frac{\text{(I.1.2.2)}  }{\mathsf{V}_{\mu,\bnu}(y,\bx)} \Big| > \mathfrak{c}_1\sqrt{\frac{\log (n)}{n^2h^{d+3}}}\Big] \leq \mathfrak{c}_2n^{-\mathfrak{c}_3}.
\end{align*}
To handle term (I.1.1), first consider the quantity $\phi_{j,i} - \phi_i$, which takes the form
\begin{align*}
&\max_i\Big|\frac{1}{n}\sum_{j=1}^n\left(\phi_{j,i} - \phi_i\right) \Big| \leq \sup_{y' \in \Y, \bx' \in \X}\Big|\frac{1}{nh}\sum_{j=1}^n\Big[\left(\Indicator\left(y'\leq y_j\right) - F(y_j|\bx')\right)\bc_1^\Trans\bP \left( \frac{y_j-y}{h} \right) \\
&\qquad\qquad\qquad - \int \left(\Indicator\left(y'\leq u\right)- F(u|\bx')\right)\bc_1^\Trans\bP \Big( \frac{u-y}{h} \Big)\diff G(u)\Big]\Big|.
\end{align*}
Then it is straightforward to show that
\begin{align*}
\Prob\Big[\max_i\sup_{ y \in \Y, \bx \in \X} \Big|\frac{1}{n}\sum_{j=1}^n\left(\phi_{j,i} - \phi_i\right) \Big| \geq \mathfrak{c}_1\sqrt{\frac{\log (n)}{nh}} \Big] \leq \mathfrak{c}_2n^{-\mathfrak{c}_3}.
\end{align*}
As a result,
\begin{align*}
\Prob\Big[\sup_{y\in\mathcal{Y},\bx\in\mathcal{X}}\Big| \frac{\text{(I.1.1)}  }{\mathsf{V}_{\mu,\bnu}(y,\bx)} \Big| > \mathfrak{c}_1\frac{\log (n)}{nh^2}\Big] \leq \mathfrak{c}_2n^{-\mathfrak{c}_3}.
\end{align*}

To conclude the proof for $\hat{\mathsf{V}}_{\mu,\bnu}(y,\bx)$, we note that replacing $\hat{\bc}_1$ by $\bc_1$ and $\hat{\bc}_2$ by $\bc_2$ only leads to an additional multiplicative factor $1 + \Op(\sqrt{\log (n)/nh^d})$. See Lemma~\ref{sa-lem:Taking hat off}.

\subsection{Proof of Lemma~\ref{sa-lem:uniform Studentization error}}

First consider $\check{\TProc}_{\order}^\circ(\eval)$. The difference between $\check{\TProc}_{\order}^\circ(\eval)$ and $\check{\SProc}_{\order}^\circ(\eval)$ is
\begin{align*}
\check{\TProc}_{\order}^\circ(\eval) - \check{\SProc}_{\order}^\circ(\eval) &= \left( \sqrt{\frac{\mathsf{V}_{\order}(\eval)}{\check{\mathsf{V}}_{\order}(\eval)}} - 1 \right)\check{\SProc}_{\order}^\circ(\eval).
\end{align*}
From Lemma~\ref{sa-lem:Uniform variance estimation}, we have
\begin{align*}
\Prob\Big[\sup_{y \in \Y, \bx \in \X}\Big| \sqrt{\frac{\mathsf{V}_{\order}(\eval)}{\check{\mathsf{V}}_{\order}(\eval)}} - 1 \Big| > \mathfrak{c}_1 \mathtt{r_{VE}}\Big] \leq \mathfrak{c}_2n^{-\mathfrak{c}_3}.
\end{align*}
To close the proof, it is straightforward to verify that
\begin{align*}
\Prob\Big[\sup_{y \in \Y, \bx \in \X}\left|\check{\SProc}_{\order}^\circ(\eval) \right| > \mathfrak{c}_1 \sqrt{\log (n)}\Big] \leq \mathfrak{c}_2n^{-\mathfrak{c}_3},
\end{align*}
which follows from the uniform convergence rate in Lemma~\ref{sa-lem:Uniform rate of convergence}. The same technique applies to the analysis of $\hat{\TProc}_{\order}^\circ(\eval)-\hat{\SProc}_{\order}^\circ(\eval)$.

\subsection{Proof of Theorem~\ref{sa-lem:strong approximation}}

See the proof of Theorem 2 in the main paper (i.e., Appendix A.6).

\subsection{Proof of Lemma~\ref{sa-lem:uniform covariance estimation}}

Consider $\check{\rho}_{\order}(y,\bx,y',\bx')$. Note that we can decompose the difference into
\begin{align*}
&\ \check{\rho}_{\order}(y,\bx,y',\bx') - \rho_{\order}(y,\bx,y',\bx')
= \frac{\check{\mathsf{C}}_{\order}(y,\bx,y',\bx')}{ \sqrt{\check{\mathsf{V}}_{\order}(\eval)\check{\mathsf{V}}_{\order}(y',\bx')}} - \frac{{\mathsf{C}}_{\order}(y,\bx,y',\bx')}{ \sqrt{{\mathsf{V}}_{\order}(\eval){\mathsf{V}}_{\order}(y',\bx')}}\\
=&\ \underbrace{\frac{\check{\mathsf{C}}_{\order}(y,\bx,y',\bx')-{\mathsf{C}}_{\order}(y,\bx,y',\bx')}{ \sqrt{\check{\mathsf{V}}_{\order}(\eval)\check{\mathsf{V}}_{\order}(y',\bx')}}}_{\textstyle \text{(I)}} + \underbrace{\rho_{\order}(y,\bx,y',\bx')\left( \sqrt{\frac{{\mathsf{V}}_{\order}(\eval){\mathsf{V}}_{\order}(y',\bx')}{\check{\mathsf{V}}_{\order}(\eval)\check{\mathsf{V}}_{\order}(y',\bx')}} - 1 \right)}_{\textstyle \text{(II)}}.
\end{align*}
The probabilistic order of the second term is given in Lemma~\ref{sa-lem:Uniform variance estimation}.

Using similar techniques as in the proof of Lemma \ref{sa-lem:Taking hat off} or \ref{sa-lem:Uniform variance estimation}, it is also straightforward to verify that term (I) has the same order. That is,
\begin{align*}
\Prob\Big[\sup_{y,y' \in \Y, \bx,\bx' \in \X}\left|\text{(I)} \right| > \mathfrak{c}_1 \mathtt{r_{VE}}\Big] \leq \mathfrak{c}_2n^{-\mathfrak{c}_3}.
\end{align*}

\subsection{Proof of Lemma~\ref{sa-lem:feasible Gaussian approx}}

Consider an $\epsilon$ discretization of $\mathcal{Y}\times \mathcal{X}$, which is denoted by $\mathcal{A}_{\epsilon} = \{(y_\ell,\bx_{\ell}^\Trans)^\Trans:\ 1\leq \ell \leq L\}$. Then one can define two Gaussian vectors, $\bz,\check{\bz}\in\mathbb{R}^L$, such that
\begin{align*}
\Cov[z_{\ell}, z_{\ell'}] = \rho(y_{\ell},\bx_{\ell},y_{\ell'},\bx_{\ell'}),\quad \Cov[\check{z}_{\ell}, \check{z}_{\ell'}|\text{Data}] = \check{\rho}(y_{\ell},\bx_{\ell},y_{\ell'},\bx_{\ell'}).
\end{align*}
Then we apply the Gaussian comparison result in Corollary 5.1 of \cite{chernozhukov2022central} (Lemma 11 in the main paper) and the error rate in Lemma~\ref{sa-lem:uniform covariance estimation}, which lead to
\begin{align*}
&\ \sup_{u \in \mathbb{R}}\Big| \Prob\Big[ \sup_{1 \leq \ell \leq L}|\check{z}_{\ell}|\leq u \Big| \bY,\bX\Big] - \Prob\Big[ \sup_{1 \leq \ell \leq L}|z_{\ell}|\leq u  \Big] \Big|\\
=&\ \sup_{u \in \mathbb{R}}\Big| \Prob\Big[ \sup_{1 \leq \ell \leq L}|\check{\GProc}_{\order}(y_{\ell},\bx_{\ell})|\leq u \Big| \bY,\bX \Big] - \Prob\Big[ \sup_{1 \leq \ell \leq L}|\GProc_{\order}(y_{\ell},\bx_{\ell})|\leq u  \Big] \Big|\\
\precsim_\Prob&\  \Big[ h^{\frac{\qOrder}{2}} +  \Big(\frac{\log (n)}{nh^{d+2}}\Big)^{\frac{1}{4}} \Big]\log(\frac{1}{\epsilon}).
\end{align*}
Since $\epsilon$ only enters the above error bound logarithmically, one can choose $\epsilon = n^{-c}$ for some $c$ large enough, so that the error that arises from discretization becomes negligible. The same applies to $\hat{\GProc}_{\order}(y_{\ell},\bx_{\ell})$.

\subsection{Proof of Theorem~\ref{sa-thm:feasible Gaussian approx}}

First consider $\check{\TProc}_{\order}^\circ(\eval)$. Since
\begin{align*}
&\ \sup_{y \in \Y, \bx \in \X}| \bar{\SProc}_{\order}(\eval)| - \sup_{y \in \Y, \bx \in \X}|\check{\TProc}_{\order}^\circ(\eval) - \bar{\SProc}_{\order}(\eval)|
\leq \sup_{y \in \Y, \bx \in \X}|\check{\TProc}_{\order}^\circ(\eval) | \\
&\leq \sup_{y \in \Y, \bx \in \X}| \bar{\SProc}_{\order}(\eval)| + \sup_{y \in \Y, \bx \in \X}|\check{\TProc}_{\order}^\circ(\eval) - \bar{\SProc}_{\order}(\eval)|,
\end{align*}
then with Lemma~\ref{sa-lem:uniform Studentization error},
\begin{align*}
\Prob\Big[\sup_{y \in \Y, \bx \in \X}|\bar{\SProc}_{\order}(\eval)| \leq u - \mathfrak{c}_1\mathtt{r_{SE}}\Big] - \mathfrak{c}_2n^{-\mathfrak{c}_3}
&\leq \Prob\Big[ \sup_{y \in \Y, \bx \in \X}|\check{\TProc}_{\order}^\circ(\eval) | \leq u \Big] \\
&\leq \Prob\Big[\sup_{y \in \Y, \bx \in \X}| \bar{\SProc}_{\order}(\eval)| \leq u + \mathfrak{c}_1\mathtt{r_{SE}}\Big] + \mathfrak{c}_2n^{-\mathfrak{c}_3}.
\end{align*}
In the above, we also used the fact that the difference  $\check{\SProc}_{\order}^\circ(\eval)-\bar{\SProc}_{\order}(\eval)$ is negligible compared to $\mathtt{r_{SE}}$ (see Remark~\ref{remark:pointwise berry-esseen}).

By applying Lemma~\ref{sa-lem:strong approximation},
\begin{align*}
&\Prob\Big[\sup_{y \in \Y, \bx \in \X}| \GProc_{\order}(\eval)| \leq u - \mathfrak{c}_1(\mathtt{r_{SE}}+\mathtt{r_{SA}})\Big] - \mathfrak{c}_2n^{-\mathfrak{c}_3}
\leq \Prob\Big[ \sup_{y \in \Y, \bx \in \X}|\check{\TProc}_{\order}^\circ(\eval) | \leq u \Big] \\
&\leq \Prob\Big[\sup_{y \in \Y, \bx \in \X}| \GProc_{\order}(\eval)| \leq u + \mathfrak{c}_1(\mathtt{r_{SE}}+\mathtt{r_{SA}})\Big] + \mathfrak{c}_2n^{-\mathfrak{c}_3}.
\end{align*}
Finally, we apply the Gaussian comparison result in Lemma~\ref{sa-lem:feasible Gaussian approx}, which implies that

\begin{align*}
&\ \sup_{u \in \mathbb{R}}\Big| \Prob\Big[ \sup_{y \in \Y, \bx \in \X}|\check{\TProc}_{\order}^\circ(\eval)|\leq u  \Big] -
\Prob\Big[ \sup_{y \in \Y, \bx \in \X}|\check{\GProc}_{\order}\eval|\leq u \Big| \bY,\bX \Big]  \Big|\\
\precsim_\Prob&\  \mathfrak{c}_2n^{-\mathfrak{c}_3} + \log (n)\sqrt{\mathtt{r_{VE}}} + \sup_{u\in\mathbb{R}}\Prob\Big[ \sup_{y \in \Y, \bx \in \X}| \GProc_{\order}(\eval)| \in [u, u+\mathfrak{c}_1(\mathtt{r_{SE}}+\mathtt{r_{SA}})] \Big].
\end{align*}
Finally, due to Theorem 2.1 of \cite{Chernozhukov-Chetverikov-Kato_2014b_AoS} (Lemma 12 in the main paper), we have
\begin{align*}
\sup_{u\in\mathbb{R}}\Prob\Big[ \sup_{y \in \Y, \bx \in \X}| \GProc_{\order}(\eval)| \in [u, u+\mathfrak{c}_1(\mathtt{r_{SE}}+\mathtt{r_{SA}})] \Big] \precsim \sqrt{\log (n)}(\mathtt{r_{SE}}+\mathtt{r_{SA}}).
\end{align*}

\subsection{Proof of Theorem~\ref{sa-thm:confidence bands}}

Note that $\theta_{\mu,\bnu}(y|\bx)$ falls into the confidence band $\check{\mathcal{C}}_{\order}(1-\alpha)$ if and only if
\begin{align*}
\sup_{ y \in \Y, \bx \in \X}\left|\frac{\check{\theta}_{\mu,\bnu}(y|\bx) - \theta_{\mu,\bnu}(y|\bx)}{\sqrt{\check{\mathsf{V}}_{\order}(\eval)}}\right| \leq \check{\mathtt{cv}}_{\order}(\alpha).
\end{align*}
A sufficient condition would then be
\begin{align*}
\sup_{ y \in \Y, \bx \in \X}\left|\check{\TProc}_{\order}^\circ(\eval)\right| + \sup_{ y \in \Y, \bx \in \X}\left| \frac{\Expectation[\check{\theta}_{\mu,\bnu}(y|\bx)|\bX] - \theta_{\mu,\bnu}(y|\bx)}{\sqrt{\check{\mathsf{V}}_{\order}(\eval)}} \right| \leq \check{\mathtt{cv}}_{\order}(\alpha).
\end{align*}
The conclusion then follows from Theorem~\ref{sa-thm:feasible Gaussian approx} and the bias calculation in Lemma~\ref{sa-lem:Bias}. The same analysis applies to $\hat{\mathcal{C}}_{\order}(1-\alpha)$.

\subsection{Proof of Theorem~\ref{sa-thm:HT-specification}}

To start, we decompose the test statistic into
\begin{align*}
\check{\TProc}_{\mathtt{PS}}(\eval) =
\check{\TProc}_{\order}^\circ(\eval)
+\frac{\Expectation[\check{\theta}_{\order}(y|\bx)] - \theta_{\order}(y|\bx) }{\sqrt{\check{\mathsf{V}}_{\order}(\eval)}}
+ \frac{\theta_{\order}(y|\bx) - \theta_{\order}(y|\bx;\hat{\bgamma}) }{\sqrt{\check{\mathsf{V}}_{\order}(\eval)}}.
\end{align*}
Then by the leading bias order in Lemma~\ref{sa-lem:Bias} and the leading variance order in Lemma~\ref{sa-lem:Variance}, we have that
\begin{align*}
\Prob\Big[ \sup_{ y \in \Y, \bx \in \X}\Big|\frac{\Expectation[\check{\theta}_{\order}(y|\bx)] - \theta_{\order}(y|\bx) }{\sqrt{\check{\mathsf{V}}_{\order}(\eval)}}\Big| > \mathfrak{c}_1 \frac{\mathtt{r_{B}}}{\mathtt{r_{V}}}\left( 1 + \mathtt{r_{VE}} \right) \Big] \leq \mathfrak{c}_2n^{-\mathfrak{c}_3}.
\end{align*}
Similarly, under the null hypothesis,
\begin{align*}
\Prob\Big[ \sup_{ y \in \Y, \bx \in \X}\Big| \frac{\theta_{\order}(y|\bx) - \theta_{\order}(y|\bx;\hat{\bgamma}) }{\sqrt{\check{\mathsf{V}}_{\order}(\eval)}} \Big| > \mathfrak{c}_1\frac{\mathtt{r_{PS}}}{\mathtt{r_{V}}}\left( 1 + \mathtt{r_{VE}} \right)  \Big]\leq \mathfrak{c}_2n^{-\mathfrak{c}_3}.
\end{align*}
Then we have the following error bound
\begin{align*}
&\sup_{u \in \mathbb{R}}\Big| \Prob\Big[ \sup_{y \in \Y, \bx \in \X}|\check{\TProc}_{\mathtt{PS}}(\eval) |\leq u  \Big] -
\Prob\Big[ \sup_{y \in \Y, \bx \in \X}|\check{\GProc}_{\order}(\eval)|\leq u \Big| \bY,\bX \Big]  \Big| \\
&\precsim_\Prob \sqrt{\log (n)}\left(\mathtt{r_{SE}}+\mathtt{r_{SA}} + \frac{\mathtt{r_{B}}+\mathtt{r_{PS}}}{\mathtt{r_{V}}} \right) + (\log (n))\sqrt{\mathtt{r_{VE}}}.
\end{align*}
As a result,
\begin{align*}
\Prob\Big[ \sup_{y \in \Y, \bx \in \X}|\check{\TProc}_{\mathtt{PS}}(\eval) | > \mathtt{cv}_{\order}(\alpha) \Big] \leq \alpha  + \mathfrak{c}\left(\sqrt{\log (n)}\left(\mathtt{r_{SE}}+\mathtt{r_{SA}} + \frac{\mathtt{r_{B}}+\mathtt{r_{PS}}}{\mathtt{r_{V}}} \right) + (\log (n))\sqrt{\mathtt{r_{VE}}}\right).
\end{align*}
The same strategy can be employed to establish results for $\hat{\TProc}_{\mathtt{PS}}(\eval)$.

\subsection{Proof of Theorem~\ref{sa-thm:HT-restriction}}

The conclusion follows directly from Theorem~\ref{sa-thm:confidence bands}.

\begin{funding}
Cattaneo gratefully acknowledges financial support from the
National Science Foundation through grants SES-1947805 and DMS-2210561, and
from the National Institute of Health (R01 GM072611-16).

Jansson gratefully acknowledges financial support from the
National Science Foundation through grant SES-1947662 and
the research support of CREATES.
\end{funding}

\bibliographystyle{imsart-nameyear} 
\bibliography{CCJM_2023_Bernoulli--bib}       